\documentclass[a4paper]{amsart}
\usepackage{amsmath}
\usepackage{amsthm}
\usepackage{amssymb}
\usepackage{amsfonts}
\usepackage{mathrsfs}
\usepackage{stmaryrd}
\usepackage{url}
\numberwithin{equation}{section}
\numberwithin{figure}{section}
\numberwithin{table}{section}
\usepackage{enumitem}

\newenvironment{renum}{
  \begin{enumerate}[label=\textup{(\roman*)}]
}{
  \end{enumerate}
}

\theoremstyle{plain}

\newtheorem*{thm*}{Theorem}

\newtheorem{pro}[subsection]{Proposition}
\theoremstyle{definition}

\newcommand\no{n\textsuperscript{o}}   % Bourbaki numeros

\newcommand{\bb}[1]{\mathbb{#1}}
\newcommand\NN{\bb{N}}
\newcommand\ZZ{\bb{Z}}
\newcommand\RR{\bb{R}}
\newcommand{\ca}[1]{\mathcal{#1}}
\DeclareMathAlphabet{\mathpzc}{OT1}{pzc}{m}{it}
\let\epsilon\varepsilon
\let\phi\varphi
\let\theta\vartheta
\newcommand\eps\epsilon
\newcommand{\dbl}{\llbracket}
\newcommand{\dbr}{\rrbracket}
\newcommand{\CX}{\ensuremath{\mathpzc{X}}} % the class X
\newcommand{\CM}{\ensuremath{\mathpzc{M}}} % the class M
\newcommand{\GL}{\operatorname{GL}}
\usepackage[all]{xy}
\newcommand{\CoxBullet}{\bullet}

\newcommand{\CircuitN}{%
\begin{xy} 0;<0.7cm,0cm>: 
(0,0)*=0{\CoxBullet}="1" , 
(1,0)*=0{\CoxBullet}="2" , 
(2,0)*=0{\CoxBullet}="3" ,
(3,0)*=0{\CoxBullet}="4" , 
(1.5,1)*=0{\CoxBullet}="5" , 
\ar@{-} "1";"2",
\ar@{..} "2";"3"
\ar@{-} "3";"4"
\ar@{-} "4";"5"
\ar@{-}"5";"1"
\end{xy}
}
\newcommand{\Linn}[2]{%
\begin{xy} 0;<0.7cm,0cm>: 
(0,0)*=0{\CoxBullet}="1" , 
(1,0)*=0{\CoxBullet}="2" , 
(2,0)*=0{\CoxBullet}="3" ,
(3,0)*=0{\CoxBullet}="4" , 
(4,0)*=0{\CoxBullet}="5" , 
(5,0)*=0{\CoxBullet}="6" ,
\ar@{-}^{#1} "1";"2",
\ar@{-} "2";"3"
\ar@{..} "3";"4"
\ar@{-} "4";"5"
\ar@{-}^{#2} "5";"6"
\end{xy}
}

\newcommand{\Lintwo}[1]{%
\begin{xy} 0;<0.7cm,0cm>: 
(0,0)*=0{\CoxBullet}="1" , 
(1,0)*=0{\CoxBullet}="2" , 
\ar@{-}^{#1} "1";"2",
\end{xy}
}

\newcommand{\LinT}[2]{%
\begin{xy} 0;<0.7cm,0cm>: 
(0,0)*=0{\CoxBullet}="1" , 
(1,0)*=0{\CoxBullet}="2" , 
(2,0)*=0{\CoxBullet}="3" ,
\ar@{-}^{#1} "1";"2",
\ar@{-}^{#2} "2";"3"
\end{xy}
}

\newcommand{\Tri}[3]{%
\begin{xy} 0;<0.7cm,0cm>: 
(0,0)*=0{\CoxBullet}="1" , 
(1,0)*=0{\CoxBullet}="2" , 
(0.5,0.866)*=0{\CoxBullet}="3" ,
\ar@{-}^{#1} "1";"2",
\ar@{-}_{#2} "2";"3",
\ar@{-}_{#3} "3";"1",
\end{xy}
}

\newcommand{\LinQ}[3] 
{\begin{xy} 0;<0.7cm,0cm>:
(0,0)*=0{\CoxBullet}="1" , 
(1,0)*=0{\CoxBullet}="2" , 
(2,0)*=0{\CoxBullet}="3" ,
(3,0)*=0{\CoxBullet}="4",
\ar@{-}^{#1} "1";"2",
\ar@{-}^{#2} "2";"3",
\ar@{-}^{#3} "3";"4"
\end{xy}}

\newcommand{\LinP}[4] 
{\begin{xy} 0;<0.7cm,0cm>:
(0,0)*=0{\CoxBullet}="1" , 
(1,0)*=0{\CoxBullet}="2" , 
(2,0)*=0{\CoxBullet}="3" ,
(3,0)*=0{\CoxBullet}="4",
(4,0)*=0{\CoxBullet}="5",
\ar@{-}^{#1} "1";"2",
\ar@{-}^{#2} "2";"3",
\ar@{-}^{#3} "3";"4",
\ar@{-}^{#4} "4";"5"
\end{xy}}

\newcommand{\LinS}[5] 
{\begin{xy} 0;<0.7cm,0cm>:
(0,0)*=0{\CoxBullet}="1" , 
(1,0)*=0{\CoxBullet}="2" , 
(2,0)*=0{\CoxBullet}="3" ,
(3,0)*=0{\CoxBullet}="4",
(4,0)*=0{\CoxBullet}="5",
(5,0)*=0{\CoxBullet}="6",
\ar@{-}^{#1} "1";"2",
\ar@{-}^{#2} "2";"3",
\ar@{-}^{#3} "3";"4",
\ar@{-}^{#4} "4";"5",
\ar@{-}^{#5} "5";"6"
\end{xy}}

\newcommand{\QFork}[3] 
{\begin{xy} 0;<0.7cm,0cm>:
(0,0)*=0{\CoxBullet}="1" , 
(1,0)*=0{\CoxBullet}="2" , 
(1.866,0.5)*=0{\CoxBullet}="3" ,
(1.866,-0.5)*=0{\CoxBullet}="4",
\ar@{-}^{#1} "1";"2",
\ar@{-}^{#2} "2";"3",
\ar@{-}_{#3} "2";"4"
\end{xy}}

\newcommand{\PFork}[3] 
{\begin{xy} 0;<0.7cm,0cm>:
(0,0)*=0{\CoxBullet}="1" , 
(1,0)*=0{\CoxBullet}="2" , 
(2,0)*=0{\CoxBullet}="3" , 
(2.866,0.5)*=0{\CoxBullet}="4" ,
(2.866,-0.5)*=0{\CoxBullet}="5",
\ar@{-}^{#1} "1";"2",
\ar@{-}^{#2} "2";"3",
\ar@{-}^{#3} "3";"4",
\ar@{-} "3";"5"
\end{xy}}

\newcommand{\SFork}[1] 
{\begin{xy} 0;<0.7cm,0cm>:
(0,0)*=0{\CoxBullet}="1" , 
(1,0)*=0{\CoxBullet}="2" , 
(2,0)*=0{\CoxBullet}="3" , 
(3,0)*=0{\CoxBullet}="4" , 
(3.866,0.5)*=0{\CoxBullet}="5" ,
(3.866,-0.5)*=0{\CoxBullet}="6",
\ar@{-} "1";"2",
\ar@{-}^{#1} "2";"3",
\ar@{-} "3";"4",
\ar@{-} "4";"5",
\ar@{-} "4";"6"
\end{xy}}

\newcommand{\nFork}[1]{%
\begin{xy} 0;<0.7cm,0cm>:
(0,0)*=0{\CoxBullet}="1" , 
(1,0)*=0{\CoxBullet}="2" , 
(2,0)*=0{\CoxBullet}="3" , 
(3,0)*=0{\CoxBullet}="4" , 
(4,0)*=0{\CoxBullet}="5" , 
(4.866,0.5)*=0{\CoxBullet}="6" ,
(4.866,-0.5)*=0{\CoxBullet}="7",
\ar@{-}^{#1} "1";"2",
\ar@{-} "2";"3",
\ar@{.} "3";"4",
\ar@{-} "4";"5",
\ar@{-} "5";"6",
\ar@{-} "5";"7"
\end{xy}}

\newcommand{\nDoubleFork}{%
\begin{xy} 0;<0.7cm,0cm>:
(0,-0.5)*=0{\CoxBullet}="1" , 
(0,0.5)*=0{\CoxBullet}="2" , 
(0.866,0)*=0{\CoxBullet}="3" , 
(1.866,0)*=0{\CoxBullet}="4" , 
(2.866,0)*=0{\CoxBullet}="5" , 
(3.866,0)*=0{\CoxBullet}="6" , 
(4.732,0.5)*=0{\CoxBullet}="7" ,
(4.732,-0.5)*=0{\CoxBullet}="8",
\ar@{-} "1";"3",
\ar@{-} "2";"3",
\ar@{-} "3";"4",
\ar@{.} "4";"5",
\ar@{-} "5";"6",
\ar@{-} "6";"7"
\ar@{-} "6";"8"
\end{xy}}

\newcommand{\Quad}[4]
{\begin{xy} 0;<0.7cm,0cm>:
(0,0)*=0{\CoxBullet}="1" , 
(1,0)*=0{\CoxBullet}="2" , 
(1,1)*=0{\CoxBullet}="3" ,
(0,1)*=0{\CoxBullet}="4",
\ar@{-}^{#1} "1";"2",
\ar@{-}_{#2} "2";"3",
\ar@{-}_{#3} "3";"4",
\ar@{-}_{#4} "4";"1"
\end{xy}}

\newcommand{\Pent}[5]
{\begin{xy} 0;<0.7cm,0cm>:
(0,0)*=0{\CoxBullet}="1" , 
(1,0)*=0{\CoxBullet}="2" , 
(1.5,0.866)*=0{\CoxBullet}="3" ,
(0.5,0.866)*=0{\CoxBullet}="4",
(-0.5,0.866)*=0{\CoxBullet}="5",
\ar@{-}^{#1} "1";"2",
\ar@{-}_{#2} "2";"3",
\ar@{-}_{#3} "3";"4",
\ar@{-}_{#4} "4";"5",
\ar@{-}_{#5} "5";"1",
\end{xy}}

\newcommand{\QELoop}[1] 
{\begin{xy} 0;<0.7cm,0cm>:
(0,0)*=0{\CoxBullet}="1" , 
(1,0)*=0{\CoxBullet}="2" , 
(1.866,0.5)*=0{\CoxBullet}="3" ,
(1.866,-0.5)*=0{\CoxBullet}="4",
\ar@{-}^{#1} "1";"2",
\ar@{-} "2";"3",
\ar@{-} "2";"4",
\ar@{-} "3";"4"
\end{xy}}

\newcommand{\PELoop}[1] 
{\begin{xy} 0;<0.7cm,0cm>:
(0,0)*=0{\CoxBullet}="1" , 
(1,0)*=0{\CoxBullet}="2" , 
(1.866,0.5)*=0{\CoxBullet}="3" ,
(1.866,-0.5)*=0{\CoxBullet}="4",
(2.732,0)*=0{\CoxBullet}="5",
\ar@{-}^{#1} "1";"2",
\ar@{-} "2";"3",
\ar@{-} "2";"4",
\ar@{-} "4";"5",
\ar@{-} "3";"5"
\end{xy}}

\newcommand{\SELoop}
{\begin{xy} 0;<0.7cm,0cm>:
(0,0)*=0{\CoxBullet}="1" , 
(1,0)*=0{\CoxBullet}="2" , 
(1.866,0.5)*=0{\CoxBullet}="3" ,
(1.866,-0.5)*=0{\CoxBullet}="4",
(2.866,0.5)*=0{\CoxBullet}="5",
(2.866,-0.5)*=0{\CoxBullet}="6",
\ar@{-} "1";"2",
\ar@{-} "2";"3",
\ar@{-} "3";"5",
\ar@{-} "5";"6",
\ar@{-} "6";"4",
\ar@{-} "4";"2"
\end{xy}}

\newcommand{\HELoop}
{\begin{xy} 0;<0.7cm,0cm>:
(0,0)*=0{\CoxBullet}="1" , 
(1,0)*=0{\CoxBullet}="2" , 
(1.866,0.5)*=0{\CoxBullet}="3" ,
(1.866,-0.5)*=0{\CoxBullet}="4",
(2.866,0.5)*=0{\CoxBullet}="5",
(2.866,-0.5)*=0{\CoxBullet}="6",
(3.732,0)*=0{\CoxBullet}="7",
\ar@{-} "1";"2",
\ar@{-} "2";"3",
\ar@{-} "3";"5",
\ar@{-} "5";"7",
\ar@{-} "7";"6",
\ar@{-} "6";"4",
\ar@{-} "4";"2",
\end{xy}}

\newcommand{\OELoop}
{\begin{xy} 0;<0.7cm,0cm>:
(0,0)*=0{\CoxBullet}="1" , 
(1,0)*=0{\CoxBullet}="2" , 
(1.866,0.5)*=0{\CoxBullet}="3" ,
(1.866,-0.5)*=0{\CoxBullet}="4",
(2.866,0.5)*=0{\CoxBullet}="5",
(2.866,-0.5)*=0{\CoxBullet}="6",
(3.866,0.5)*=0{\CoxBullet}="7",
(3.866,-0.5)*=0{\CoxBullet}="8",
\ar@{-} "1";"2",
\ar@{-} "2";"3",
\ar@{-} "3";"5",
\ar@{-} "5";"7",
\ar@{-} "7";"8"
\ar@{-} "8";"6",
\ar@{-} "6";"4",
\ar@{-} "4";"2"
\end{xy}}

\newcommand{\NELoop}
{\begin{xy} 0;<0.7cm,0cm>:
(0,0)*=0{\CoxBullet}="1" , 
(1,0)*=0{\CoxBullet}="2" , 
(1.866,0.5)*=0{\CoxBullet}="3" ,
(1.866,-0.5)*=0{\CoxBullet}="4",
(2.866,0.5)*=0{\CoxBullet}="5",
(2.866,-0.5)*=0{\CoxBullet}="6",
(3.866,0.5)*=0{\CoxBullet}="7",
(3.866,-0.5)*=0{\CoxBullet}="8",
(4.732,0)*=0{\CoxBullet}="9",
\ar@{-} "1";"2",
\ar@{-} "2";"3",
\ar@{-} "3";"5",
\ar@{-} "5";"7",
\ar@{-} "7";"9",
\ar@{-} "9";"8",
\ar@{-} "8";"6",
\ar@{-} "6";"4",
\ar@{-} "4";"2",
\end{xy}}

\newcommand{\Tetra}
{\begin{xy} 0;<0.7cm,0cm>:
(0,0)*=0{\CoxBullet}="1" , 
(0,1)*=0{\CoxBullet}="2" , 
(-0.866,-0.5)*=0{\CoxBullet}="3" ,
(0.866,-0.5)*=0{\CoxBullet}="4",
\ar@{-} "1";"2",
\ar@{-} "2";"3",
\ar@{-} "2";"4",
\ar@{-} "3";"4",
\ar@{-} "1";"4",
\ar@{-} "3";"1"
\end{xy}}

\newcommand{\BigQ}
{\begin{xy} 0;<0.7cm,0cm>:
(0,0)*=0{\CoxBullet}="1" , 
(0.707,0.707)*=0{\CoxBullet}="2" , 
(-0.707,0.707)*=0{\CoxBullet}="3" ,
(-0.707,-0.707)*=0{\CoxBullet}="4",
(0.707,-0.707)*=0{\CoxBullet}="5",
\ar@{-} "1";"2",
\ar@{-} "2";"3",
\ar@{-} "3";"4",
\ar@{-} "4";"5",
\ar@{-} "5";"2",
\ar@{-} "4";"1"
\end{xy}}

\newcommand{\EsixL}[2]{%
\begin{xy} 0;<0.7cm,0cm>:
(0,0)*=0{\CoxBullet}="1" , 
(1,0)*=0{\CoxBullet}="2" , 
(2,0)*=0{\CoxBullet}="3" ,
(3,0)*=0{\CoxBullet}="4",
(4,0)*=0{\CoxBullet}="5",
(2,1)*=0{\CoxBullet}="6",
\ar@{-}^{#1} "1";"2",
\ar@{-} "2";"3",
\ar@{-} "3";"4",
\ar@{-}^{#2} "4";"5",
\ar@{-} "3";"6"
\end{xy}}

\newcommand{\EsevenL}[1]{%
\begin{xy} 0;<0.7cm,0cm>:
(0,0)*=0{\CoxBullet}="1" , 
(1,0)*=0{\CoxBullet}="2" , 
(2,0)*=0{\CoxBullet}="3" ,
(3,0)*=0{\CoxBullet}="4",
(4,0)*=0{\CoxBullet}="5",
(2,1)*=0{\CoxBullet}="6",
(5,0)*=0{\CoxBullet}="7",
\ar@{-} "1";"2",
\ar@{-} "2";"3",
\ar@{-} "3";"4",
\ar@{-} "4";"5",
\ar@{-} "3";"6"
\ar@{-}^{#1} "5";"7",
\end{xy}}

\newcommand{\EeightL}[1]{%
\begin{xy} 0;<0.7cm,0cm>:
(0,0)*=0{\CoxBullet}="1" , 
(1,0)*=0{\CoxBullet}="2" , 
(2,0)*=0{\CoxBullet}="3" ,
(3,0)*=0{\CoxBullet}="4",
(4,0)*=0{\CoxBullet}="5",
(2,1)*=0{\CoxBullet}="6",
(5,0)*=0{\CoxBullet}="7",
(6,0)*=0{\CoxBullet}="8",
\ar@{-} "1";"2",
\ar@{-} "2";"3",
\ar@{-} "3";"4",
\ar@{-} "4";"5",
\ar@{-} "3";"6",
\ar@{-} "5";"7",
\ar@{-}^{#1} "7";"8",
\end{xy}}

\newcommand{\TEeightL}[1]{%
\begin{xy} 0;<0.7cm,0cm>:
(0,0)*=0{\CoxBullet}="1" , 
(1,0)*=0{\CoxBullet}="2" , 
(2,0)*=0{\CoxBullet}="3" ,
(3,0)*=0{\CoxBullet}="4",
(4,0)*=0{\CoxBullet}="5",
(5,0)*=0{\CoxBullet}="6",
(6,0)*=0{\CoxBullet}="7",
(7,0)*=0{\CoxBullet}="8",
(2,1)*=0{\CoxBullet}="9",
\ar@{-} "1";"2",
\ar@{-} "2";"3",
\ar@{-} "3";"4",
\ar@{-} "4";"5",
\ar@{-} "5";"6",
\ar@{-} "6";"7",
\ar@{-}^{#1} "7";"8"
\ar@{-} "3";"9",
\end{xy}}

\newcommand{\TTEsix}{%
\begin{xy} 0;<0.7cm,0cm>:
(0,0)*=0{\CoxBullet}="1" , 
(1,0)*=0{\CoxBullet}="2" , 
(2,0)*=0{\CoxBullet}="3" ,
(3,0)*=0{\CoxBullet}="4",
(3.866,0.5)*=0{\CoxBullet}="5",
(4.866,0.5)*=0{\CoxBullet}="6",
(3.866,-0.5)*=0{\CoxBullet}="7",
(4.866,-0.5)*=0{\CoxBullet}="8",
\ar@{-} "1";"2",
\ar@{-} "2";"3",
\ar@{-} "3";"4",
\ar@{-} "4";"5",
\ar@{-} "5";"6",
\ar@{-} "4";"7",
\ar@{-} "7";"8",
\end{xy}}

\newcommand{\TEsix}{%
\begin{xy} 0;<0.7cm,0cm>:
(0,0)*=0{\CoxBullet}="1" , 
(1,0)*=0{\CoxBullet}="2" , 
(2,0)*=0{\CoxBullet}="3" ,
(2.866,0.5)*=0{\CoxBullet}="4",
(3.866,0.5)*=0{\CoxBullet}="5",
(2.866,-0.5)*=0{\CoxBullet}="6",
(3.866,-0.5)*=0{\CoxBullet}="7",
\ar@{-} "1";"2",
\ar@{-} "2";"3",
\ar@{-} "3";"4",
\ar@{-} "4";"5",
\ar@{-} "3";"6",
\ar@{-} "6";"7",
\end{xy}}

\newcommand{\TTEseven}{%
\begin{xy} 0;<0.7cm,0cm>:
(0,0)*=0{\CoxBullet}="1" , 
(1,0)*=0{\CoxBullet}="2" , 
(2,0)*=0{\CoxBullet}="3" ,
(3,0)*=0{\CoxBullet}="4",
(4,0)*=0{\CoxBullet}="5",
(5,0)*=0{\CoxBullet}="6",
(6,0)*=0{\CoxBullet}="7",
(7,0)*=0{\CoxBullet}="8",
(3,1)*=0{\CoxBullet}="9",
\ar@{-} "1";"2",
\ar@{-} "2";"3",
\ar@{-} "3";"4",
\ar@{-} "4";"5",
\ar@{-} "5";"6",
\ar@{-} "6";"7",
\ar@{-} "7";"8",
\ar@{-} "4";"9",
\end{xy}}

\newcommand{\TEseven}{%
\begin{xy} 0;<0.7cm,0cm>:
(0,0)*=0{\CoxBullet}="1" , 
(1,0)*=0{\CoxBullet}="2" , 
(2,0)*=0{\CoxBullet}="3" ,
(3,0)*=0{\CoxBullet}="4",
(4,0)*=0{\CoxBullet}="5",
(5,0)*=0{\CoxBullet}="6",
(6,0)*=0{\CoxBullet}="7",
(2,1)*=0{\CoxBullet}="8",
\ar@{-} "1";"2",
\ar@{-} "2";"3",
\ar@{-} "3";"4",
\ar@{-} "4";"5",
\ar@{-} "5";"6",
\ar@{-} "6";"7",
\ar@{-} "3";"8",
\end{xy}}

\newcommand{\TTEeightL}[1]{%
\begin{xy} 0;<0.7cm,0cm>:
(0,0)*=0{\CoxBullet}="1" , 
(1,0)*=0{\CoxBullet}="2" , 
(2,0)*=0{\CoxBullet}="3" ,
(3,0)*=0{\CoxBullet}="4",
(4,0)*=0{\CoxBullet}="5",
(5,0)*=0{\CoxBullet}="6",
(6,0)*=0{\CoxBullet}="7",
(7,0)*=0{\CoxBullet}="8",
(8,0)*=0{\CoxBullet}="9",
(2,1)*=0{\CoxBullet}="10",
\ar@{-} "1";"2",
\ar@{-} "2";"3",
\ar@{-} "3";"4",
\ar@{-} "4";"5",
\ar@{-} "5";"6",
\ar@{-} "6";"7",
\ar@{-} "7";"8",
\ar@{-}^{#1} "8";"9"
\ar@{-} "3";"10",
\end{xy}}

\newcommand{\TDfour}[1]
{\begin{xy} 0;<0.7cm,0cm>:
(0,0)*=0{\CoxBullet}="1" , 
(1,0)*=0{\CoxBullet}="2" , 
(0,1)*=0{\CoxBullet}="3" ,
(-1,0)*=0{\CoxBullet}="4",
(0,-1)*=0{\CoxBullet}="5",
\ar@{-}^{#1} "4";"1",
\ar@{-} "1";"3",
\ar@{-} "1";"2",
\ar@{-} "1";"5",
\end{xy}}

\newcommand{\TTDfour}[1]{%
\begin{xy} 0;<0.7cm,0cm>:
(0,0)*=0{\CoxBullet}="1" , 
(1,0)*=0{\CoxBullet}="2" , 
(2,0)*=0{\CoxBullet}="3" ,
(3,0)*=0{\CoxBullet}="4",
(2,-1)*=0{\CoxBullet}="5",
(2,1)*=0{\CoxBullet}="6",
\ar@{-}^{#1} "1";"2",
\ar@{-} "2";"3",
\ar@{-} "3";"4",
\ar@{-} "3";"5",
\ar@{-} "3";"6"
\end{xy}}

\newcommand{\TTDfive}{%
\begin{xy} 0;<0.7cm,0cm>:
(0,0)*=0{\CoxBullet}="1" , 
(1,0)*=0{\CoxBullet}="2" , 
(2,0)*=0{\CoxBullet}="3" ,
(3,0)*=0{\CoxBullet}="4",
(4,0)*=0{\CoxBullet}="5",
(1,1)*=0{\CoxBullet}="6",
(2,1)*=0{\CoxBullet}="7",
\ar@{-} "1";"2",
\ar@{-} "2";"3",
\ar@{-} "3";"4",
\ar@{-} "4";"5",
\ar@{-} "2";"6",
\ar@{-} "3";"7",
\end{xy}}

\newcommand{\TTDsix}{%
\begin{xy} 0;<0.7cm,0cm>:
(0,0)*=0{\CoxBullet}="1" , 
(1,0)*=0{\CoxBullet}="2" , 
(2,0)*=0{\CoxBullet}="3" ,
(3,0)*=0{\CoxBullet}="4",
(4,0)*=0{\CoxBullet}="5",
(5,0)*=0{\CoxBullet}="6",
(1,1)*=0{\CoxBullet}="7",
(3,1)*=0{\CoxBullet}="8",
\ar@{-} "1";"2",
\ar@{-} "2";"3",
\ar@{-} "3";"4",
\ar@{-} "4";"5",
\ar@{-} "5";"6",
\ar@{-} "2";"7",
\ar@{-} "4";"8",
\end{xy}}

\newcommand{\TTDseven}{%
\begin{xy} 0;<0.7cm,0cm>:
(0,0)*=0{\CoxBullet}="1" , 
(1,0)*=0{\CoxBullet}="2" , 
(2,0)*=0{\CoxBullet}="3" ,
(3,0)*=0{\CoxBullet}="4",
(4,0)*=0{\CoxBullet}="5",
(5,0)*=0{\CoxBullet}="6",
(6,0)*=0{\CoxBullet}="7",
(1,1)*=0{\CoxBullet}="8",
(4,1)*=0{\CoxBullet}="9",
\ar@{-} "1";"2",
\ar@{-} "2";"3",
\ar@{-} "3";"4",
\ar@{-} "4";"5",
\ar@{-} "5";"6",
\ar@{-} "6";"7",
\ar@{-} "2";"8",
\ar@{-} "5";"9",
\end{xy}}

\newcommand{\TTDeight}{%
\begin{xy} 0;<0.7cm,0cm>:
(0,0)*=0{\CoxBullet}="1" , 
(1,0)*=0{\CoxBullet}="2" , 
(2,0)*=0{\CoxBullet}="3" ,
(3,0)*=0{\CoxBullet}="4",
(4,0)*=0{\CoxBullet}="5",
(5,0)*=0{\CoxBullet}="6",
(6,0)*=0{\CoxBullet}="7",
(7,0)*=0{\CoxBullet}="8",
(1,1)*=0{\CoxBullet}="9",
(5,1)*=0{\CoxBullet}="10",
\ar@{-} "1";"2",
\ar@{-} "2";"3",
\ar@{-} "3";"4",
\ar@{-} "4";"5",
\ar@{-} "5";"6",
\ar@{-} "6";"7",
\ar@{-} "7";"8",
\ar@{-} "2";"9",
\ar@{-} "6";"10",
\end{xy}}

\newcommand{\PStar}
{\begin{xy} 0;<0.7cm,0cm>:
(0,0)*=0{\CoxBullet}="1" , 
(0,1)*=0{\CoxBullet}="2" , 
a(162)*=0{\CoxBullet}="3" ,
a(234)*=0{\CoxBullet}="4",
a(306)*=0{\CoxBullet}="5",
a(18)*=0{\CoxBullet}="6",
\ar@{-} "1";"2",
\ar@{-} "1";"3",
\ar@{-} "1";"4",
\ar@{-} "1";"5",
\ar@{-} "1";"6"
\end{xy}}

\newcommand{\SixLoop}[2]
{\begin{xy} 0;<0.7cm,0cm>:
(0,0)*=0{\CoxBullet}="1" , 
(0.866,0.5)*=0{\CoxBullet}="2" , 
(0.866,-0.5)*=0{\CoxBullet}="3" ,
(1.732,0.5)*=0{\CoxBullet}="4",
(1.732,-0.5)*=0{\CoxBullet}="5",
(2.732,0)*=0{\CoxBullet}="6",
\ar@{-} "1";"2",
\ar@{-}^{#2} "2";"4",
\ar@{-} "4";"6",
\ar@{-} "6";"5",
\ar@{-}_{#1} "5";"3",
\ar@{-} "3";"1",
\end{xy}}

\newcommand{\CQuad}[5]{%
\begin{xy} 0;<0.7cm,0cm>:
(0,0)*=0{\CoxBullet}="1" , 
(1,0)*=0{\CoxBullet}="2" , 
(1,1)*=0{\CoxBullet}="3" ,
(0,1)*=0{\CoxBullet}="4",
\ar@{-}^{#1} "1";"2",
\ar@{-}^{#2} "2";"3",
\ar@{-}^{#3} "3";"4",
\ar@{-}^{#4} "4";"1",
\ar@{-}|{#5} "1";"3"
\end{xy}}

\newcommand{\Diagrambox}[3]{\raisebox{0pt}[#1][#2]{#3}}

\newcommand{\CoxGrUIV}{%
  \begin{xy} 0;<0.7cm,0cm>:
    (0,0)*=0{\CoxBullet}="1" , 
    (0,1)*=0{\CoxBullet}="2" , 
    (-0.866,-0.5)*=0{\CoxBullet}="3" ,
    (0.866,-0.5)*=0{\CoxBullet}="4",
    \ar@{-}|{\infty} "1";"2",
    \ar@{-}|{\infty} "2";"3",
    \ar@{-}|{\infty} "2";"4",
    \ar@{-}|{\infty} "3";"4",
    \ar@{-}|{\infty} "1";"4",
    \ar@{-}|{\infty} "3";"1"
  \end{xy}
}

\newcommand{\CoxGrSItwo}[1]{\Lintwo{m}}

\newcommand{\CoxGrHCI}[2]{\LinT{#1}{#2}}
\newcommand{\CoxGrHCII}[3]{\Diagrambox{24pt}{3pt}{\Tri{#1}{#2}{#3}}}

\usepackage[margin=2cm]{geometry}
\usepackage{microtype}
\usepackage{multicol}

\title{Data about hyperbolic Coxeter systems}
\author{T. Terragni} 
\email{tom.terragni@gmail.com}
\urladdr{https://sites.google.com/site/tomterragni}

\date{\today}

\usepackage{tikz}
\usetikzlibrary{matrix,arrows}
\DeclareMathOperator{\cleq}{\preceq}
\newcommand{\WS}{(W,S)}
\newcommand{\WSp}{(W',S')}
\newcommand{\pWS}{p_{\WS}(t)}
\newcommand{\WI}{(W_I,I)}
\newcommand{\pWI}{p_{\WI}(t)}
\newcommand{\spr}{\ca{F}}
\DeclareMathOperator{\degree}{deg}

\newcommand{\SH}[3]{\ensuremath{\langle #1,#2,#3\rangle}}
\newcommand{\EHC}[1]{\ensuremath{\mathbf{hc}(#1)}}%
\newcommand{\EHNC}[1]{\ensuremath{\mathbf{hnc}(#1)}}
\newcounter{mycnt}

\stepcounter{tocdepth}

\usepackage{breqn}
\newcommand{\breakingcomma}{%
  \begingroup\lccode`~=`,
  \lowercase{\endgroup\expandafter\def\expandafter~\expandafter{~\penalty0 }}}

\begin{document}

%%%%%%%%%%%%%%%%%%%%%%%%%%%%%%%%%%% 
%% Abstract
\begin{abstract}
We collect several data about Coxeter systems (cf.~\cite{bourbaki--gal46,humphreys--rgcg}), with particular emphasis on the hyperbolic ones.

For each ($\cleq$-minimal) hyperbolic Coxeter system $(W,S)$ the Poincar\'e series 
\[p_{(W,S)}(t)=\sum_{w\in W} t^{\ell(w)}\]
and the growth rate 
\[ \omega(W,S)=\limsup_n \sqrt[n]{a_n}\]
are explicitly computed using Magma (cf.~\cite{MAGMA}). 

These computations were performed in connection to the proof of \cite[Thm.~B]{terragni--gcg}.

Since the Poincar\'e series represents a rational function, one may recover the sequence $(a_k)_{k\in \NN_0}$ through a linear recurrence relation on the coefficients, provided that enough terms at the beginning of the sequence are known. 
For each Coxeter system the initial coefficients $(a_k)_{k=0}^N$ are computed, where $N$ is the degree of the numerator of $p_{(W,S)}(t)$.

\medskip Consider the environment. Do not bulk print.


\end{abstract}

\maketitle
\addtocontents{toc}{%
\protect\begin{multicols}{2}%
}
\tableofcontents

%%%%%%%%%%%%%%%%%%%%%%%%%%%%%%%%%%% 
%% Body
\section{Coxeter groups}\label{s:cox-gps}
\subsection{Coxeter systems}\label{ss:cox-syst}
Let $S$ be a finite set, and let $M$ be an $(S\times S)$-matrix such that 
$m_{s,s}=1$, and $m_{s,r}=m_{r,s}\in \ZZ_{\geq 2}\cup\{\infty\}$ for all $s,r\in S$, $s\neq r$. Then $M$ is a \emph{Coxeter matrix} over $S$. 

The \emph{Coxeter system} associated with a Coxeter matrix $M$ over $S$ is the pair $\WS$ where $W$ is the group 
\begin{equation}\label{eq:cox-pres}
W=W(M)=\langle S\mid (sr)^{m_{s,r}} \,\text{ if }m_{s,r}<\infty\rangle.
\end{equation}

The Coxeter matrix $M$ (or, equivalently, the presentation \eqref{eq:cox-pres}) is often encoded in the Coxeter graph $\Gamma(M)$ (cf.~\cite[Ch.~IV \no 1.9]{bourbaki--gal46}).
Either datum is called the \emph{type} of $\WS$, and a Coxeter system is called \emph{irreducible} if its Coxeter graph $\Gamma$ is connected.

If $I\subseteq S$ let $W_I=\langle I\rangle \leq W$. The \emph{parabolic subsystem} $(W_I,I)$ is a Coxeter system in its own right, with Coxeter matrix $M_I=(m_{s,r})_{s,r\in I}$.
Its Coxeter graph is the graph induced from $\Gamma$ by the vertices in $I$.
The finite set $\spr=\spr(W,S)=\{I\subseteq S\mid |W_I|<\infty\}$ is called the set of \emph{spherical residues}.

A \emph{Coxeter-isomorphism} $\phi\colon \WS\to (W',S')$ of Coxeter systems of types $M$ and $M'$, respectively, is a bijection $\phi\colon S\to S'$ such that $m'_{\phi(s),\phi(r)}=m_{s,r}$ for all $s,r\in S$.

Any Coxeter group $\WS$ is linear via the Tits' reflection representation $\rho\colon W\to \GL(\RR^S)$ (cf.~\cite[Ch.~V, \S4]{bourbaki--gal46}). The representation $\rho$ is determined by the symmetric matrix\footnote{For short, we put $\frac\pi\infty=0$.} $B=B_M=\left(-\cos\frac{\pi}{m_{s,r}}\right)_{s,r\in S}$, and the signature of $B$ induces the following tetrachotomy on irreducible Coxeter systems.

\begin{renum}
\item If $B$ is positive definite, then $\WS$ is \emph{spherical},
\item if $B$ is positive semidefinite with $0$ a simple eigenvalue, then $\WS$ is \emph{affine}, 
\item if $B$ has $|S|-1$ positive and $1$ negative eigenvalue, then $\WS$ is \emph{hyperbolic}\footnote{There are several non-compatible notions of hyperbolicity (cf.~\cite[Note 6.9]{davis--gtcg}). In the present work ``hyperbolic'' coincides with Bourbaki's notion (cf.~\cite[Ch. V, \S4, Ex.13]{bourbaki--gal46}).}, or
\item none of the above conditions applies.
\end{renum}
The irreducible Coxeter system $\WS$ is spherical if, and only if, $W$ is a finite group. The classification of spherical and affine systems is classical (cf.~\cite[Ch. VI]{bourbaki--gal46}). 
For a characterisation of hyperbolic Coxeter systems see \S\ref{ss:min}.

\subsection{Poincar\'e series}\label{s:poincare-ser}
The \emph{Poincar\'e series} $\pWS$ of $\WS$ is defined as
\begin{equation}\label{eq:pws}
\pWS = \sum_{w\in W} t^{\ell(w)}\in \ZZ\dbl t\dbr.
\end{equation} 
If $\WS$ is spherical then $\pWS$ is a polynomial, which can be explicitly computed in terms of the degrees of the polynomial invariants of $\WS$, simply known as the \emph{degrees} of $\WS$ (cf.~\cite{solomon--ofcg}, \cite[Ch.~3]{humphreys--rgcg}):
\[\pWS = \prod_{d \text{ degree}} \frac{t^d -1}{t-1}.\]

If $\widetilde{X}_n$ is an affine Coxeter system and $X_n$ is the corresponding spherical system, then a theorem of Bott (cf.~\cite{bott--amttlg}) provides the formula
\[p_{\widetilde{X}_n}(t)= p_{X_n}(t) \prod_{d\text{ degree}} \frac{1}{1-t^{d-1}}.\]
 
For arbitrary Coxeter systems, the Poincar\'e series can be computed using the following property. 

\begin{pro}[{\cite{steinberg--elag}}]\label{pro:ps}
Let $\WS$ be a Coxeter system with Poincar\'e series $\pWS$. Then
\begin{equation}\label{eq:sph-steinberg}
  \frac{1}{p_{(W,S)}(t^{-1})}=\sum_{I\in \spr}\frac{(-1)^{|I|}}{\pWI},
\end{equation}
where $\spr=\spr\WS$. In particular, the Poincar\'e series $\pWS$ is a rational function.
\end{pro}

\subsection{Growth rate}
If $a_n=a_n^{\WS}=|\{w\in W \mid \ell(w)=n\}|$, then \eqref{eq:pws} can be rewritten as 
\[\pWS = \sum_{n\in \NN_0} a_nt^n.\]
Therefore the number $\omega\WS=\limsup_n \sqrt[n]{a_n}$ is the inverse of the convergency radius of the (complex) power series $\pWS$. In the context of the growth of finitely generated groups, $\omega\WS$ is called the \emph{(exponential) growth rate} of $\WS$.

%% \begin{lem} Let $(W_1,S_1)$ and $(W_2,S_2)$ be Coxeter systems, and let $\WS=(W_1\times W_2,S_1\sqcup S_2)$ be their product. 
%% Then $\omega(W,S)=\max\{\omega{(W_1,S_1)}, \omega{(W_2,S_2)} \}$.
%% \end{lem}

\subsection{The partial order $\cleq$ on the class of Coxeter systems}\label{s:orders}

Let $\WS$ and $\WSp$ be Coxeter systems with Coxeter matrices $M$, $M'$ respectively. 
Define $\WS\cleq\WSp$ whenever there exists an injective map $\phi\colon S\to S'$ such that $m_{s,r}\leq m'_{\phi(s),\phi(r)}$ for all $s,r\in S$ (cf.~\cite[\S6]{mcmullen--cgsnhm}).

In particular, if $\WS$ and $\WSp$ are Coxeter-isomorphic (cf.~\S\ref{ss:cox-syst}) then $\WS\cleq\WSp$ and $\WSp\cleq\WS$. 
Therefore the preorder $\cleq$ descends to a partial order on the set of Coxeter-isomorphism classes of Coxeter systems. 
With a mild abuse of notation we will avoid the distinction between a Coxeter system and its Coxeter-isomorphism class.

Let $\CX$ be the set of (Coxeter-isomorphism classes of) non-affine, non-spherical, irreducible Coxeter systems, and let  $\CM =\min_{\cleq} \CX$ be the set of $\cleq$-minimal elements of $\CX$.

\subsection{Hyperbolic Coxeter systems and $\CM$}\label{ss:min}

It is well known that hyperbolic Coxeter systems are characterised as those systems such that every proper irreducible parabolic subsystem is either of spherical or affine type (cf.~\cite[Ch.~V, \S4, Ex.~13]{bourbaki--gal46}). 
By minimality, $\CM$ consists of hyperbolic Coxeter systems, which are classified in an infinite family of rank-three systems, and $72$ exceptions of rank $|S|\geq 4$ (cf.~\cite[\S\S6.8--6.9]{humphreys--rgcg}).
The infinite family consists of the $\langle a,b,c\rangle$-triangle groups with $\frac1a+\frac1b+\frac1c<1$, and among those only the $\langle2,3,7\rangle$, $\langle3,3,4\rangle$ and $\langle2,4,5\rangle$-triangle groups are $\cleq$-minimal. 
Among the $72$ exceptions, $35$ are in $\CM$ (cf.~\cite[Thm. 6.6, Table 5]{mcmullen--cgsnhm}).

\subsection{Linear recurrence relations}

It is well-known that a power series represents a rational function if, and only if, there is a linear reccurrence formula for the coefficients of the series.

Explicitly, if 
\[ \sum_{k=0}^\infty a_k t^k =\frac{n(t)}{d(t)},\]
with
\[n(t)= \sum_{k=0}^N \nu_k t^k \qquad\text{ and } \qquad d(t)=\sum_{k=0}^D \delta_k t^k,\]
one has
\[\sum_{r=0}^D \delta_r a_{k-r} =\nu_k \quad\text{for all }k.\]
In particular, for all $k> N$ one has the required linear recurrence relation
\begin{equation}\label{eq:rec}
  a_k = -\frac{1}{\delta_0} \sum_{r=1}^D\delta_r a_{k-r}.
\end{equation}

Then, the relation \eqref{eq:rec} together with the values $a_k$ for $k\leq N=\degree(n(t))$, completely determines the sequence $(a_k)_{k\geq 0}$. 
The values $(a_k)_{k=0}^N$ are displayed for each system in \S\ref{s:hyp}.

\newpage
\section{Data about hyperbolic Coxeter systems}\label{s:hyp}

\subsection{The infinite family}
For $a,b,c\in \NN\cup\{\infty\}$ satisfying\footnote{For short, put $\frac1\infty=0$.}  
$\frac1a+\frac1b+\frac1c<1$ let us denote 
\[\SH{a}{b}{c}=\CoxGrHCII{a}{b}{c},\]
and then on has, for $a,b,c<\infty$,
\[\begin{split}
p_{\SH{a}{b}{c}}&= \frac{(1+t) (1-t^a) (1-t^b) (1-t^c)}{%
    1-2t+t^{a+1}+t^{b+1}+t^{c+1}-t^{a+b}-t^{a+c}-t^{b+c}+2 t^{a+b+c}-t^{a+b+c+1}} ,\\
p_{\SH{a}{b}{\infty}}&= \frac{(1+t) (1-t^a) (1-t^b)}{1-2t+t^{a+1}+t^{b+1}-t^{a+b}},\\%  {t^{a+b}-t^{a+1}-t^{b+1}+2t-1},\\
p_{\SH{a}{\infty}{\infty}}&= \frac{(1+t) (1-t^a)}{1 -2t+  t^{a+1}} ,\text{ and}\\
p_{\SH{\infty}{\infty}{\infty}}&= \frac{1+t}{1-2t} .\\
\end{split}\]

For the reader's convenience the computations for the three $\cleq$-minimal Coxeter systems in the series follow.

\newpage
\subsubsection{\SH{2}{3}{7}}
\[\Gamma(W,S) = \CoxGrHCI{}{7}\]

\vskip12pt 
\[ M=\begin{pmatrix} 
1 & 3 & 2 \\
3 & 1 & 7  \\
2 & 7 & 1  
\end{pmatrix}\]
Numerator of $p_{(W,S)}(t)$:
\begin{dmath*} n_{(W,S)}(t) = 
1+4t+8t^{2}+11t^{3}+12t^{4}+12t^{5}+12t^{6}+11t^{7}+8t^{8}+4t^{9}+t^{10}
\end{dmath*}
Denominator of $p_{(W,S)}(t)$:
\begin{dmath*} d_{(W,S)}(t) = 
1+t-t^{3}-t^{4}-t^{5}-t^{6}-t^{7}+t^{9}+t^{10}
\end{dmath*}
Initial values of $a_k$:
\begin{dmath*}  (a_k)_{k=0}^{10}=(1,3, 5, 7, 9, 12, 16, 20, 24, 28, 33).
\end{dmath*}
Exponential growth rate:
\[\omega(W,S) =\lambda_{\textrm{Lehmer}}=
1.17628081825991750654407033847
\dots\]
Cocompact? Yes\\
$(W,S)\in \CM$? Yes
\newpage

\subsubsection{\SH{2}{4}{5}}
\[\Gamma(W,S) = \CoxGrHCI{4}{5}\]

\vskip12pt 
\[ M=\begin{pmatrix} 
1 & 4 & 2 \\
4 & 1 & 5  \\
2 & 5 & 1  
\end{pmatrix}\]
Numerator of $p_{(W,S)}(t)$:
\begin{dmath*} n_{(W,S)}(t) = 
1+3t+5t^{2}+7t^{3}+8t^{4}+7t^{5}+5t^{6}+3t^{7}+t^{8}
\end{dmath*}
Denominator of $p_{(W,S)}(t)$:
\begin{dmath*} d_{(W,S)}(t) = 
1-t^{3}-t^{4}-t^{5}+t^{8}
\end{dmath*}
Initial values of $a_k$:
\begin{dmath*}  (a_k)_{k=0}^{8}=(1, 3, 5, 8, 12, 16, 21, 28, 36).
\end{dmath*}
Exponential growth rate:
\[\omega(W,S) =
1.28063815626775759670190253271
\dots\]
Cocompact? Yes\\
$(W,S)\in \CM$? Yes
\newpage

\subsubsection{\SH{3}{3}{4}}
\[\Gamma(W,S) = \CoxGrHCII{}{}{4}\]

\vskip12pt 
\[ M=\begin{pmatrix} 
1 & 3 & 4 \\
3 & 1 & 3  \\
4 & 3 & 1  
\end{pmatrix}\]
Numerator of $p_{(W,S)}(t)$:
\begin{dmath*} n_{(W,S)}(t) = 
1+3t+5t^{2}+6t^{3}+5t^{4}+3t^{5}+t^{6}
\end{dmath*}
Denominator of $p_{(W,S)}(t)$:
\begin{dmath*} d_{(W,S)}(t) = 
1-t^{2}-t^{3}-t^{4}+t^{6}
\end{dmath*}
Initial values of $a_k$:
\begin{dmath*}  (a_k)_{k=0}^{6}=(1, 3, 6, 10, 15, 22, 31).
\end{dmath*}
Exponential growth rate:
\[\omega(W,S) =
1.40126836793985491510176409562
\dots\]
Cocompact? Yes\\
$(W,S)\in \CM$? Yes
\newpage

\subsection{The 72 exceptional types}
%% \input{tables-1-10}
%% \input{tables-11-20}
%% \input{tables-21-46}
%% \input{tables-47-60}
%% \input{tables-61-65}
%% \input{tables-66}
%% \input{tables-67}
%% \input{tables-68-72}
% produced with RUNfortables.magma
\subsubsection{\EHC{1}}
\setcounter{mycnt}{3}
\[\Gamma(W,S) = \quad\csname CoxGrHC\Roman{mycnt}\endcsname\]

\vskip12pt 
\[ M=\begin{pmatrix} 
1 & 5 & 2 & 2 \\
5 & 1 & 3 & 2 \\
2 & 3 & 1 & 4 \\
2 & 2 & 4 & 1
\end{pmatrix}\]
Numerator of $p_{(W,S)}(t)$:
\begin{dmath*} n_{(W,S)}(t) = 
-1-2t-3t^{2}-5t^{3}-6t^{4}-7t^{5}-8t^{6}-8t^{7}-8t^{8}-8t^{9}-7t^{10}-6t^{11}-5t^{12}-3t^{13}-2t^{14}-t^{15}
\end{dmath*}
Denominator of $p_{(W,S)}(t)$:
\begin{dmath*} d_{(W,S)}(t) = 
-1+2t-2t^{2}+2t^{3}-t^{4}+t^{5}-t^{6}+t^{7}-t^{8}+t^{9}-t^{10}+t^{11}-2t^{12}+2t^{13}-2t^{14}+t^{15}
\end{dmath*}
Initial values of $a_k$:
\begin{dmath*}  (a_k)_{k=0}^{15}=(1
, 4
, 9
, 17
, 29
, 46
, 70
, 103
, 148
, 210
, 295
, 411
, 569
, 783
, 1074
, 1470
).
 \end{dmath*}
Exponential growth rate:
\[\omega(W,S) =
1.35999971171150086545103049530
\dots\]
Cocompact? Yes\\
$(W,S)\in \CM$? Yes\\
Magma command: \texttt{HyperbolicCoxeterGraph(1)} or \texttt{HyperbolicCoxeterMatrix(1)} \\
\newpage
\subsubsection{\EHC{2}}
\setcounter{mycnt}{4}
\[\Gamma(W,S) = \quad\csname CoxGrHC\Roman{mycnt}\endcsname\]

\vskip12pt 
\[ M=\begin{pmatrix} 
1 & 5 & 2 & 3 \\
5 & 1 & 3 & 2 \\
2 & 3 & 1 & 2 \\
3 & 2 & 2 & 1
\end{pmatrix}\]
Numerator of $p_{(W,S)}(t)$:
\begin{dmath*} n_{(W,S)}(t) = 
-1-2t-2t^{2}-2t^{3}-2t^{4}-3t^{5}-3t^{6}-2t^{7}-2t^{8}-2t^{9}-2t^{10}-t^{11}
\end{dmath*}
Denominator of $p_{(W,S)}(t)$:
\begin{dmath*} d_{(W,S)}(t) = 
-1+2t-t^{2}+t^{4}-2t^{5}+2t^{6}-t^{7}+t^{9}-2t^{10}+t^{11}
\end{dmath*}
Initial values of $a_k$:
\begin{dmath*}  (a_k)_{k=0}^{11}=(1
, 4
, 9
, 16
, 26
, 41
, 62
, 90
, 128
, 181
, 254
, 352
).
 \end{dmath*}
Exponential growth rate:
\[\omega(W,S) =
1.35098033771623731021140357306
\dots\]
Cocompact? Yes\\
$(W,S)\in \CM$? Yes\\
Magma command: \texttt{HyperbolicCoxeterGraph(2)} or \texttt{HyperbolicCoxeterMatrix(2)} \\
\newpage
\subsubsection{\EHC{3}}
\setcounter{mycnt}{5}
\[\Gamma(W,S) = \quad\csname CoxGrHC\Roman{mycnt}\endcsname\]

\vskip12pt 
\[ M=\begin{pmatrix} 
1 & 5 & 2 & 2 \\
5 & 1 & 3 & 2 \\
2 & 3 & 1 & 5 \\
2 & 2 & 5 & 1
\end{pmatrix}\]
Numerator of $p_{(W,S)}(t)$:
\begin{dmath*} n_{(W,S)}(t) = 
-1-2t-2t^{2}-3t^{3}-4t^{4}-4t^{5}-4t^{6}-4t^{7}-4t^{8}-4t^{9}-3t^{10}-2t^{11}-2t^{12}-t^{13}
\end{dmath*}
Denominator of $p_{(W,S)}(t)$:
\begin{dmath*} d_{(W,S)}(t) = 
-1+2t-t^{2}+t^{4}-t^{5}+t^{6}-t^{7}+t^{8}-t^{9}+t^{11}-2t^{12}+t^{13}
\end{dmath*}
Initial values of $a_k$:
\begin{dmath*}  (a_k)_{k=0}^{13}=(1
, 4
, 9
, 17
, 30
, 50
, 80
, 125
, 193
, 296
, 450
, 680
, 1025
, 1541
).
 \end{dmath*}
Exponential growth rate:
\[\omega(W,S) =
1.49671107560954952105387691751
\dots\]
Cocompact? Yes\\
$(W,S)\in \CM$? No\\
Magma command: \texttt{HyperbolicCoxeterGraph(3)} or \texttt{HyperbolicCoxeterMatrix(3)} \\
\newpage
\subsubsection{\EHC{4}}
\setcounter{mycnt}{6}
\[\Gamma(W,S) = \quad\csname CoxGrHC\Roman{mycnt}\endcsname\]

\vskip12pt 
\[ M=\begin{pmatrix} 
1 & 5 & 2 & 2 \\
5 & 1 & 3 & 3 \\
2 & 3 & 1 & 2 \\
2 & 3 & 2 & 1
\end{pmatrix}\]
Numerator of $p_{(W,S)}(t)$:
\begin{dmath*} n_{(W,S)}(t) = 
-1-t-t^{2}-2t^{3}-t^{4}-2t^{5}-2t^{6}-t^{7}-2t^{8}-t^{9}-t^{10}-t^{11}
\end{dmath*}
Denominator of $p_{(W,S)}(t)$:
\begin{dmath*} d_{(W,S)}(t) = 
-1+3t-4t^{2}+4t^{3}-3t^{4}+2t^{5}-2t^{6}+3t^{7}-4t^{8}+4t^{9}-3t^{10}+t^{11}
\end{dmath*}
Initial values of $a_k$:
\begin{dmath*}  (a_k)_{k=0}^{11}=(1
, 4
, 9
, 17
, 29
, 47
, 74
, 113
, 170
, 253
, 374
, 550
).
 \end{dmath*}
Exponential growth rate:
\[\omega(W,S) =
1.44842304024420580152689399820
\dots\]
Cocompact? Yes\\
$(W,S)\in \CM$? Yes\\
Magma command: \texttt{HyperbolicCoxeterGraph(4)} or \texttt{HyperbolicCoxeterMatrix(4)} \\
\newpage
\subsubsection{\EHC{5}}
\setcounter{mycnt}{7}
\[\Gamma(W,S) = \quad\csname CoxGrHC\Roman{mycnt}\endcsname\]

\vskip12pt 
\[ M=\begin{pmatrix} 
1 & 4 & 2 & 3 \\
4 & 1 & 3 & 2 \\
2 & 3 & 1 & 3 \\
3 & 2 & 3 & 1
\end{pmatrix}\]
Numerator of $p_{(W,S)}(t)$:
\begin{dmath*} n_{(W,S)}(t) = 
-1-2t-2t^{2}-3t^{3}-3t^{4}-2t^{5}-2t^{6}-t^{7}
\end{dmath*}
Denominator of $p_{(W,S)}(t)$:
\begin{dmath*} d_{(W,S)}(t) = 
-1+2t-2t^{3}+2t^{4}-2t^{6}+t^{7}
\end{dmath*}
Initial values of $a_k$:
\begin{dmath*}  (a_k)_{k=0}^{7}=(1
, 4
, 10
, 21
, 39
, 68
, 114
, 186
).
 \end{dmath*}
Exponential growth rate:
\[\omega(W,S) =
1.55603019132268226053042892660
\dots\]
Cocompact? Yes\\
$(W,S)\in \CM$? Yes\\
Magma command: \texttt{HyperbolicCoxeterGraph(5)} or \texttt{HyperbolicCoxeterMatrix(5)} \\
\newpage
\subsubsection{\EHC{6}}
\setcounter{mycnt}{8}
\[\Gamma(W,S) = \quad\csname CoxGrHC\Roman{mycnt}\endcsname\]

\vskip12pt 
\[ M=\begin{pmatrix} 
1 & 4 & 2 & 3 \\
4 & 1 & 3 & 2 \\
2 & 3 & 1 & 4 \\
3 & 2 & 4 & 1
\end{pmatrix}\]
Numerator of $p_{(W,S)}(t)$:
\begin{dmath*} n_{(W,S)}(t) = 
-1-2t-2t^{2}-3t^{3}-3t^{4}-2t^{5}-2t^{6}-t^{7}
\end{dmath*}
Denominator of $p_{(W,S)}(t)$:
\begin{dmath*} d_{(W,S)}(t) = 
-1+2t-t^{3}+t^{4}-2t^{6}+t^{7}
\end{dmath*}
Initial values of $a_k$:
\begin{dmath*}  (a_k)_{k=0}^{7}=(1
, 4
, 10
, 22
, 44
, 84
, 156
, 284
).
 \end{dmath*}
Exponential growth rate:
\[\omega(W,S) =
1.78164359860800194739266335197
\dots\]
Cocompact? Yes\\
$(W,S)\in \CM$? No\\
Magma command: \texttt{HyperbolicCoxeterGraph(6)} or \texttt{HyperbolicCoxeterMatrix(6)} \\
\newpage
\subsubsection{\EHC{7}}
\setcounter{mycnt}{9}
\[\Gamma(W,S) = \quad\csname CoxGrHC\Roman{mycnt}\endcsname\]

\vskip12pt 
\[ M=\begin{pmatrix} 
1 & 5 & 2 & 3 \\
5 & 1 & 3 & 2 \\
2 & 3 & 1 & 4 \\
3 & 2 & 4 & 1
\end{pmatrix}\]
Numerator of $p_{(W,S)}(t)$:
\begin{dmath*} n_{(W,S)}(t) = 
-1-t-t^{2}-2t^{3}-t^{4}-2t^{5}-2t^{6}-t^{7}-2t^{8}-t^{9}-t^{10}-t^{11}
\end{dmath*}
Denominator of $p_{(W,S)}(t)$:
\begin{dmath*} d_{(W,S)}(t) = 
-1+3t-3t^{2}+2t^{3}-2t^{5}+2t^{6}-2t^{8}+3t^{9}-3t^{10}+t^{11}
\end{dmath*}
Initial values of $a_k$:
\begin{dmath*}  (a_k)_{k=0}^{11}=(1
, 4
, 10
, 22
, 45
, 89
, 172
, 328
, 622
, 1176
, 2220
, 4186
).
 \end{dmath*}
Exponential growth rate:
\[\omega(W,S) =
1.88320350591352586416894746536
\dots\]
Cocompact? Yes\\
$(W,S)\in \CM$? No\\
Magma command: \texttt{HyperbolicCoxeterGraph(7)} or \texttt{HyperbolicCoxeterMatrix(7)} \\
\newpage
\subsubsection{\EHC{8}}
\setcounter{mycnt}{10}
\[\Gamma(W,S) = \quad\csname CoxGrHC\Roman{mycnt}\endcsname\]

\vskip12pt 
\[ M=\begin{pmatrix} 
1 & 5 & 2 & 3 \\
5 & 1 & 3 & 2 \\
2 & 3 & 1 & 3 \\
3 & 2 & 3 & 1
\end{pmatrix}\]
Numerator of $p_{(W,S)}(t)$:
\begin{dmath*} n_{(W,S)}(t) = 
-1-t-t^{2}-2t^{3}-t^{4}-2t^{5}-2t^{6}-t^{7}-2t^{8}-t^{9}-t^{10}-t^{11}
\end{dmath*}
Denominator of $p_{(W,S)}(t)$:
\begin{dmath*} d_{(W,S)}(t) = 
-1+3t-3t^{2}+t^{3}+2t^{4}-4t^{5}+4t^{6}-2t^{7}-t^{8}+3t^{9}-3t^{10}+t^{11}
\end{dmath*}
Initial values of $a_k$:
\begin{dmath*}\breakingcomma
 (a_k)_{k=0}^{11}=(1
, 4
, 10
, 21
, 40
, 73
, 130
, 228
, 396
, 684
, 1178
, 2025
).
 \end{dmath*}
Exponential growth rate:
\[\omega(W,S) =
1.71336036067247148570537838116
\dots\]
Cocompact? Yes\\
$(W,S)\in \CM$? No\\
Magma command: \texttt{HyperbolicCoxeterGraph(8)} or \texttt{HyperbolicCoxeterMatrix(8)} \\
\newpage
\subsubsection{\EHC{9}}
\setcounter{mycnt}{11}
\[\Gamma(W,S) = \quad\csname CoxGrHC\Roman{mycnt}\endcsname\]

\vskip12pt 
\[ M=\begin{pmatrix} 
1 & 5 & 2 & 3 \\
5 & 1 & 3 & 2 \\
2 & 3 & 1 & 5 \\
3 & 2 & 5 & 1
\end{pmatrix}\]
Numerator of $p_{(W,S)}(t)$:
\begin{dmath*} n_{(W,S)}(t) = 
-1-t-t^{3}-t^{4}-t^{5}-t^{6}-t^{8}-t^{9}
\end{dmath*}
Denominator of $p_{(W,S)}(t)$:
\begin{dmath*} d_{(W,S)}(t) = 
-1+3t-2t^{2}-t^{3}+3t^{4}-3t^{5}+t^{6}+2t^{7}-3t^{8}+t^{9}
\end{dmath*}
Initial values of $a_k$:
\begin{dmath*}\breakingcomma
 (a_k)_{k=0}^{9}=(1
, 4
, 10
, 22
, 46
, 94
, 188
, 372
, 734
, 1446
).
 \end{dmath*}
Exponential growth rate:
\[\omega(W,S) =
1.96355303898882461408647248762
\dots\]
Cocompact? Yes\\
$(W,S)\in \CM$? No\\
Magma command: \texttt{HyperbolicCoxeterGraph(9)} or \texttt{HyperbolicCoxeterMatrix(9)} \\
\newpage
\subsubsection{\EHC{10}}
\setcounter{mycnt}{12}
\[\Gamma(W,S) = \quad\csname CoxGrHC\Roman{mycnt}\endcsname\]

\vskip12pt 
\[ M=\begin{pmatrix} 
1 & 5 & 2 & 2 & 2 \\
5 & 1 & 3 & 2 & 2 \\
2 & 3 & 1 & 3 & 2 \\
2 & 2 & 3 & 1 & 4 \\
2 & 2 & 2 & 4 & 1
\end{pmatrix}\]
Numerator of $p_{(W,S)}(t)$:
\begin{dmath*} n_{(W,S)}(t) = 
1+4t+9t^{2}+16t^{3}+26t^{4}+40t^{5}+58t^{6}+80t^{7}+106t^{8}+136t^{9}+170t^{10}+208t^{11}+249t^{12}+292t^{13}+337t^{14}+384t^{15}+432t^{16}+480t^{17}+528t^{18}+576t^{19}+623t^{20}+668t^{21}+711t^{22}+752t^{23}+790t^{24}+824t^{25}+854t^{26}+880t^{27}+902t^{28}+920t^{29}+933t^{30}+940t^{31}+942t^{32}+940t^{33}+933t^{34}+920t^{35}+902t^{36}+880t^{37}+854t^{38}+824t^{39}+790t^{40}+752t^{41}+711t^{42}+668t^{43}+623t^{44}+576t^{45}+528t^{46}+480t^{47}+432t^{48}+384t^{49}+337t^{50}+292t^{51}+249t^{52}+208t^{53}+170t^{54}+136t^{55}+106t^{56}+80t^{57}+58t^{58}+40t^{59}+26t^{60}+16t^{61}+9t^{62}+4t^{63}+t^{64}
\end{dmath*}
Denominator of $p_{(W,S)}(t)$:
\begin{dmath*} d_{(W,S)}(t) = 
1-t-t^{3}+2t^{4}-2t^{5}+t^{6}-3t^{7}+3t^{8}-3t^{9}+3t^{10}-5t^{11}+5t^{12}-5t^{13}+6t^{14}-7t^{15}+8t^{16}-8t^{17}+9t^{18}-9t^{19}+11t^{20}-11t^{21}+12t^{22}-11t^{23}+14t^{24}-13t^{25}+14t^{26}-13t^{27}+16t^{28}-14t^{29}+15t^{30}-14t^{31}+17t^{32}-14t^{33}+15t^{34}-14t^{35}+16t^{36}-13t^{37}+14t^{38}-13t^{39}+14t^{40}-11t^{41}+12t^{42}-11t^{43}+11t^{44}-9t^{45}+9t^{46}-8t^{47}+8t^{48}-7t^{49}+6t^{50}-5t^{51}+5t^{52}-5t^{53}+3t^{54}-3t^{55}+3t^{56}-3t^{57}+t^{58}-2t^{59}+2t^{60}-t^{61}-t^{63}+t^{64}
\end{dmath*}
Initial values of $a_k$:
\begin{dmath*}\breakingcomma
 (a_k)_{k=0}^{64}=(1
, 5
, 14
, 31
, 60
, 106
, 176
, 280
, 432
, 651
, 964
, 1409
, 2038
, 2923
, 4165
, 5904
, 8334
, 11725
, 16453
, 23040
, 32210
, 44970
, 62719
, 87399
, 121707
, 169389
, 235647
, 327705
, 455597
, 633259
, 880042
, 1222819
, 1698911
, 2360144
, 3278488
, 4553889
, 6325143
, 8784994
, 12201104
, 16945176
, 23533389
, 32682552
, 45388091
, 63032342
, 87534935
, 121561635
, 168814355
, 234433890
, 325559185
, 452103992
, 627835396
, 871871540
, 1210761505
, 1681373607
, 2334906324
, 3242458212
, 4502762768
, 6252929715
, 8683361252
, 12058466031
, 16745424351
, 23254133928
, 32292681987
, 44844378522
, 62274731206
).
 \end{dmath*}
Exponential growth rate:
\[\omega(W,S) =
1.38868480369313087111252071358
\dots\]
Cocompact? Yes\\
$(W,S)\in \CM$? No\\
Magma command: \texttt{HyperbolicCoxeterGraph(10)} or \texttt{HyperbolicCoxeterMatrix(10)} \\
\newpage
% produced with RUNfortables.magma
\subsubsection{\EHC{11}}
\setcounter{mycnt}{13}
\[\Gamma(W,S) = \quad\csname CoxGrHC\Roman{mycnt}\endcsname\]

\vskip12pt 
\[ M=\begin{pmatrix} 
1 & 5 & 2 & 2 & 2 \\
5 & 1 & 3 & 2 & 2 \\
2 & 3 & 1 & 3 & 2 \\
2 & 2 & 3 & 1 & 3 \\
2 & 2 & 2 & 3 & 1
\end{pmatrix}\]
Numerator of $p_{(W,S)}(t)$:
\begin{dmath*} n_{(W,S)}(t) = 
1+4t+9t^{2}+16t^{3}+25t^{4}+36t^{5}+49t^{6}+64t^{7}+81t^{8}+100t^{9}+121t^{10}+144t^{11}+168t^{12}+192t^{13}+216t^{14}+240t^{15}+264t^{16}+288t^{17}+312t^{18}+336t^{19}+359t^{20}+380t^{21}+399t^{22}+416t^{23}+431t^{24}+444t^{25}+455t^{26}+464t^{27}+471t^{28}+476t^{29}+478t^{30}+476t^{31}+471t^{32}+464t^{33}+455t^{34}+444t^{35}+431t^{36}+416t^{37}+399t^{38}+380t^{39}+359t^{40}+336t^{41}+312t^{42}+288t^{43}+264t^{44}+240t^{45}+216t^{46}+192t^{47}+168t^{48}+144t^{49}+121t^{50}+100t^{51}+81t^{52}+64t^{53}+49t^{54}+36t^{55}+25t^{56}+16t^{57}+9t^{58}+4t^{59}+t^{60}
\end{dmath*}
Denominator of $p_{(W,S)}(t)$:
\begin{dmath*} d_{(W,S)}(t) = 
1-t-t^{7}+t^{8}-t^{9}+t^{10}-t^{11}+t^{14}-t^{15}+t^{16}-2t^{17}+2t^{18}-t^{19}+t^{20}-t^{21}+t^{22}-t^{23}+2t^{24}-2t^{25}+2t^{26}-2t^{27}+2t^{28}-t^{29}+t^{30}-t^{31}+2t^{32}-2t^{33}+2t^{34}-2t^{35}+2t^{36}-t^{37}+t^{38}-t^{39}+t^{40}-t^{41}+2t^{42}-2t^{43}+t^{44}-t^{45}+t^{46}-t^{49}+t^{50}-t^{51}+t^{52}-t^{53}-t^{59}+t^{60}
\end{dmath*}
Initial values of $a_k$:
\begin{dmath*}\breakingcomma
 (a_k)_{k=0}^{60}=(1
, 5
, 14
, 30
, 55
, 91
, 140
, 205
, 290
, 400
, 541
, 720
, 945
, 1225
, 1571
, 1996
, 2516
, 3151
, 3926
, 4871
, 6021
, 7417
, 9108
, 11153
, 13623
, 16604
, 20200
, 24536
, 29762
, 36058
, 43639
, 52762
, 63735
, 76928
, 92786
, 111845
, 134749
, 162270
, 195333
, 235047
, 282742
, 340014
, 408779
, 491337
, 590449
, 709429
, 852254
, 1023695
, 1229475
, 1476459
, 1772884
, 2128635
, 2555576
, 3067945
, 3682827
, 4420721
, 5306222
, 6368840
, 7643983
, 9174137
, 11010281
).
 \end{dmath*}
Exponential growth rate:
\[\omega(W,S) =
1.19988303671613064177882753906
\dots\]
Cocompact? Yes\\
$(W,S)\in \CM$? Yes\\
Magma command: \texttt{HyperbolicCoxeterGraph(11)} or \texttt{HyperbolicCoxeterMatrix(11)} \\
\newpage
\subsubsection{\EHC{12}}
\setcounter{mycnt}{14}
\[\Gamma(W,S) = \quad\csname CoxGrHC\Roman{mycnt}\endcsname\]

\vskip12pt 
\[ M=\begin{pmatrix} 
1 & 5 & 2 & 2 & 2 \\
5 & 1 & 3 & 2 & 2 \\
2 & 3 & 1 & 3 & 2 \\
2 & 2 & 3 & 1 & 5 \\
2 & 2 & 2 & 5 & 1
\end{pmatrix}\]
Numerator of $p_{(W,S)}(t)$:
\begin{dmath*} n_{(W,S)}(t) = 
1+4t+9t^{2}+16t^{3}+25t^{4}+36t^{5}+49t^{6}+64t^{7}+81t^{8}+100t^{9}+121t^{10}+144t^{11}+168t^{12}+192t^{13}+216t^{14}+240t^{15}+264t^{16}+288t^{17}+312t^{18}+336t^{19}+359t^{20}+380t^{21}+399t^{22}+416t^{23}+431t^{24}+444t^{25}+455t^{26}+464t^{27}+471t^{28}+476t^{29}+478t^{30}+476t^{31}+471t^{32}+464t^{33}+455t^{34}+444t^{35}+431t^{36}+416t^{37}+399t^{38}+380t^{39}+359t^{40}+336t^{41}+312t^{42}+288t^{43}+264t^{44}+240t^{45}+216t^{46}+192t^{47}+168t^{48}+144t^{49}+121t^{50}+100t^{51}+81t^{52}+64t^{53}+49t^{54}+36t^{55}+25t^{56}+16t^{57}+9t^{58}+4t^{59}+t^{60}
\end{dmath*}
Denominator of $p_{(W,S)}(t)$:
\begin{dmath*} d_{(W,S)}(t) = 
1-t-t^{3}-t^{7}-t^{9}+2t^{10}-2t^{11}+2t^{12}-2t^{13}+2t^{14}+2t^{16}-2t^{17}+2t^{18}-2t^{19}+6t^{20}-3t^{21}+4t^{22}-3t^{23}+4t^{24}+4t^{26}-3t^{27}+4t^{28}-3t^{29}+8t^{30}-3t^{31}+4t^{32}-3t^{33}+4t^{34}+4t^{36}-3t^{37}+4t^{38}-3t^{39}+6t^{40}-2t^{41}+2t^{42}-2t^{43}+2t^{44}+2t^{46}-2t^{47}+2t^{48}-2t^{49}+2t^{50}-t^{51}-t^{53}-t^{57}-t^{59}+t^{60}
\end{dmath*}
Initial values of $a_k$:
\begin{dmath*}\breakingcomma
 (a_k)_{k=0}^{60}=(1
, 5
, 14
, 31
, 61
, 111
, 191
, 317
, 514
, 820
, 1292
, 2017
, 3127
, 4821
, 7402
, 11330
, 17302
, 26376
, 40159
, 61088
, 92857
, 141068
, 214221
, 325207
, 493574
, 748973
, 1136377
, 1723992
, 2615261
, 3967072
, 6017366
, 9127010
, 13843307
, 20996318
, 31844934
, 48298426
, 73252457
, 111098669
, 168497540
, 255550538
, 387577835
, 587814418
, 891499094
, 1352076001
, 2050599979
, 3110000980
, 4716717890
, 7153508967
, 10849213011
, 16454218500
, 24954922142
, 37847319774
, 57400278978
, 87054825446
, 132029711446
, 200239837221
, 303689153595
, 460583175654
, 698532877831
, 1059413807404
, 1606735546889
).
 \end{dmath*}
Exponential growth rate:
\[\omega(W,S) =
1.51662692404935293158277836808
\dots\]
Cocompact? Yes\\
$(W,S)\in \CM$? No\\
Magma command: \texttt{HyperbolicCoxeterGraph(12)} or \texttt{HyperbolicCoxeterMatrix(12)} \\
\newpage
\subsubsection{\EHC{13}}
\setcounter{mycnt}{15}
\[\Gamma(W,S) = \quad\csname CoxGrHC\Roman{mycnt}\endcsname\]

\vskip12pt 
\[ M=\begin{pmatrix} 
1 & 5 & 2 & 2 & 2 \\
5 & 1 & 3 & 2 & 2 \\
2 & 3 & 1 & 3 & 3 \\
2 & 2 & 3 & 1 & 2 \\
2 & 2 & 3 & 2 & 1
\end{pmatrix}\]
Numerator of $p_{(W,S)}(t)$:
\begin{dmath*} n_{(W,S)}(t) = 
1+4t+9t^{2}+16t^{3}+25t^{4}+36t^{5}+49t^{6}+64t^{7}+81t^{8}+100t^{9}+121t^{10}+144t^{11}+168t^{12}+192t^{13}+216t^{14}+240t^{15}+264t^{16}+288t^{17}+312t^{18}+336t^{19}+359t^{20}+380t^{21}+399t^{22}+416t^{23}+431t^{24}+444t^{25}+455t^{26}+464t^{27}+471t^{28}+476t^{29}+478t^{30}+476t^{31}+471t^{32}+464t^{33}+455t^{34}+444t^{35}+431t^{36}+416t^{37}+399t^{38}+380t^{39}+359t^{40}+336t^{41}+312t^{42}+288t^{43}+264t^{44}+240t^{45}+216t^{46}+192t^{47}+168t^{48}+144t^{49}+121t^{50}+100t^{51}+81t^{52}+64t^{53}+49t^{54}+36t^{55}+25t^{56}+16t^{57}+9t^{58}+4t^{59}+t^{60}
\end{dmath*}
Denominator of $p_{(W,S)}(t)$:
\begin{dmath*} d_{(W,S)}(t) = 
1-t-t^{3}+t^{4}-t^{5}-t^{7}+t^{8}-t^{9}-t^{11}+2t^{12}-t^{13}+t^{14}-t^{15}+2t^{16}-t^{17}+t^{18}-t^{19}+3t^{20}-t^{21}+2t^{22}-t^{23}+3t^{24}-t^{25}+2t^{26}-t^{27}+3t^{28}-t^{29}+3t^{30}-t^{31}+3t^{32}-t^{33}+2t^{34}-t^{35}+3t^{36}-t^{37}+2t^{38}-t^{39}+3t^{40}-t^{41}+t^{42}-t^{43}+2t^{44}-t^{45}+t^{46}-t^{47}+2t^{48}-t^{49}-t^{51}+t^{52}-t^{53}-t^{55}+t^{56}-t^{57}-t^{59}+t^{60}
\end{dmath*}
Initial values of $a_k$:
\begin{dmath*}\breakingcomma
 (a_k)_{k=0}^{60}=(1
, 5
, 14
, 31
, 60
, 106
, 177
, 285
, 447
, 688
, 1045
, 1572
, 2346
, 3478
, 5130
, 7536
, 11036
, 16125
, 23520
, 34261
, 49856
, 72490
, 105331
, 152971
, 222068
, 322276
, 467590
, 678301
, 983827
, 1426815
, 2069090
, 3000278
, 4350312
, 6307558
, 9145091
, 13258795
, 19222585
, 27868466
, 40402604
, 58573580
, 84916313
, 123105683
, 178469174
, 258730056
, 375084773
, 543764840
, 788301123
, 1142806360
, 1656733858
, 2401776174
, 3481866138
, 5047675378
, 7317632842
, 10608395387
, 15379021056
, 22295007143
, 32321126305
, 46856011231
, 67927262280
, 98474295402
, 142758387259
).
 \end{dmath*}
Exponential growth rate:
\[\omega(W,S) =
1.44970181390854096032366536730
\dots\]
Cocompact? Yes\\
$(W,S)\in \CM$? Yes\\
Magma command: \texttt{HyperbolicCoxeterGraph(13)} or \texttt{HyperbolicCoxeterMatrix(13)} \\
\newpage
\subsubsection{\EHC{14}}
\setcounter{mycnt}{16}
\[\Gamma(W,S) = \quad\csname CoxGrHC\Roman{mycnt}\endcsname\]

\vskip12pt 
\[ M=\begin{pmatrix} 
1 & 4 & 2 & 2 & 3 \\
4 & 1 & 3 & 2 & 2 \\
2 & 3 & 1 & 3 & 2 \\
2 & 2 & 3 & 1 & 3 \\
3 & 2 & 2 & 3 & 1
\end{pmatrix}\]
Numerator of $p_{(W,S)}(t)$:
\begin{dmath*} n_{(W,S)}(t) = 
1+5t+14t^{2}+30t^{3}+55t^{4}+90t^{5}+134t^{6}+185t^{7}+240t^{8}+295t^{9}+346t^{10}+390t^{11}+424t^{12}+445t^{13}+452t^{14}+445t^{15}+424t^{16}+390t^{17}+346t^{18}+295t^{19}+240t^{20}+185t^{21}+134t^{22}+90t^{23}+55t^{24}+30t^{25}+14t^{26}+5t^{27}+t^{28}
\end{dmath*}
Denominator of $p_{(W,S)}(t)$:
\begin{dmath*} d_{(W,S)}(t) = 
1-t^{2}-t^{3}-t^{4}-2t^{5}-2t^{6}-t^{7}+t^{8}+t^{9}+2t^{10}+2t^{11}+3t^{12}+2t^{13}+3t^{14}+2t^{15}+3t^{16}+2t^{17}+2t^{18}+t^{19}+t^{20}-t^{21}-2t^{22}-2t^{23}-t^{24}-t^{25}-t^{26}+t^{28}
\end{dmath*}
Initial values of $a_k$:
\begin{dmath*}\breakingcomma
 (a_k)_{k=0}^{28}=(1
, 5
, 15
, 36
, 76
, 148
, 273
, 486
, 843
, 1435
, 2410
, 4009
, 6623
, 10887
, 17833
, 29135
, 47511
, 77372
, 125879
, 204652
, 332551
, 540183
, 877221
, 1424278
, 2312177
, 3753224
, 6091955
, 9887499
, 16047226
).
 \end{dmath*}
Exponential growth rate:
\[\omega(W,S) =
1.62282456264156210393713091832
\dots\]
Cocompact? Yes\\
$(W,S)\in \CM$? Yes\\
Magma command: \texttt{HyperbolicCoxeterGraph(14)} or \texttt{HyperbolicCoxeterMatrix(14)} \\
\newpage
\subsubsection{\EHNC{1}}
\setcounter{mycnt}{1}
\[\Gamma(W,S) = \quad\csname CoxGrHNC\Roman{mycnt}\endcsname\]

\vskip12pt 
\[ M=\begin{pmatrix} 
1 & 4 & 2 & 3 \\
4 & 1 & 4 & 2 \\
2 & 4 & 1 & 3 \\
3 & 2 & 3 & 1
\end{pmatrix}\]
Numerator of $p_{(W,S)}(t)$:
\begin{dmath*} n_{(W,S)}(t) = 
1+2t+2t^{2}+3t^{3}+3t^{4}+2t^{5}+2t^{6}+t^{7}
\end{dmath*}
Denominator of $p_{(W,S)}(t)$:
\begin{dmath*} d_{(W,S)}(t) = 
1-2t+t^{3}-t^{4}+t^{6}
\end{dmath*}
Initial values of $a_k$:
\begin{dmath*}\breakingcomma
 (a_k)_{k=0}^{7}=(1
, 4
, 10
, 22
, 44
, 84
, 157
, 289
).
 \end{dmath*}
Exponential growth rate:
\[\omega(W,S) =
1.81240361926804266078969297255
\dots\]
Cocompact? No\\
$(W,S)\in \CM$? No\\
Magma command: \texttt{HyperbolicCoxeterGraph(15)} or \texttt{HyperbolicCoxeterMatrix(15)} \\
\newpage
\subsubsection{\EHNC{2}}
\setcounter{mycnt}{2}
\[\Gamma(W,S) = \quad\csname CoxGrHNC\Roman{mycnt}\endcsname\]

\vskip12pt 
\[ M=\begin{pmatrix} 
1 & 4 & 2 & 3 \\
4 & 1 & 4 & 2 \\
2 & 4 & 1 & 4 \\
3 & 2 & 4 & 1
\end{pmatrix}\]
Numerator of $p_{(W,S)}(t)$:
\begin{dmath*} n_{(W,S)}(t) = 
1+2t+2t^{2}+3t^{3}+3t^{4}+2t^{5}+2t^{6}+t^{7}
\end{dmath*}
Denominator of $p_{(W,S)}(t)$:
\begin{dmath*} d_{(W,S)}(t) = 
1-2t+t^{7}
\end{dmath*}
Initial values of $a_k$:
\begin{dmath*}\breakingcomma
 (a_k)_{k=0}^{7}=(1
, 4
, 10
, 23
, 49
, 100
, 202
, 404
).
 \end{dmath*}
Exponential growth rate:
\[\omega(W,S) =
1.98358284342432633038562929339
\dots\]
Cocompact? No\\
$(W,S)\in \CM$? No\\
Magma command: \texttt{HyperbolicCoxeterGraph(16)} or \texttt{HyperbolicCoxeterMatrix(16)} \\
\newpage
\subsubsection{\EHNC{3}}
\setcounter{mycnt}{3}
\[\Gamma(W,S) = \quad\csname CoxGrHNC\Roman{mycnt}\endcsname\]

\vskip12pt 
\[ M=\begin{pmatrix} 
1 & 4 & 2 & 4 \\
4 & 1 & 4 & 2 \\
2 & 4 & 1 & 4 \\
4 & 2 & 4 & 1
\end{pmatrix}\]
Numerator of $p_{(W,S)}(t)$:
\begin{dmath*} n_{(W,S)}(t) = 
1+2t+2t^{2}+2t^{3}+t^{4}
\end{dmath*}
Denominator of $p_{(W,S)}(t)$:
\begin{dmath*} d_{(W,S)}(t) = 
1-2t-2t^{3}+3t^{4}
\end{dmath*}
Initial values of $a_k$:
\begin{dmath*}\breakingcomma
 (a_k)_{k=0}^{4}=(1
, 4
, 10
, 24
, 54
).
 \end{dmath*}
Exponential growth rate:
\[\omega(W,S) =
2.13039543476727879287505602649
\dots\]
Cocompact? No\\
$(W,S)\in \CM$? No\\
Magma command: \texttt{HyperbolicCoxeterGraph(17)} or \texttt{HyperbolicCoxeterMatrix(17)} \\
\newpage
\subsubsection{\EHNC{4}}
\setcounter{mycnt}{4}
\[\Gamma(W,S) = \quad\csname CoxGrHNC\Roman{mycnt}\endcsname\]

\vskip12pt 
\[ M=\begin{pmatrix} 
1 & 6 & 2 & 3 \\
6 & 1 & 3 & 2 \\
2 & 3 & 1 & 3 \\
3 & 2 & 3 & 1
\end{pmatrix}\]
Numerator of $p_{(W,S)}(t)$:
\begin{dmath*} n_{(W,S)}(t) = 
1+2t+3t^{2}+4t^{3}+4t^{4}+4t^{5}+3t^{6}+2t^{7}+t^{8}
\end{dmath*}
Denominator of $p_{(W,S)}(t)$:
\begin{dmath*} d_{(W,S)}(t) = 
1-2t+t^{2}-t^{3}-t^{5}+t^{6}+t^{8}
\end{dmath*}
Initial values of $a_k$:
\begin{dmath*}\breakingcomma
 (a_k)_{k=0}^{8}=(1
, 4
, 10
, 21
, 40
, 74
, 135
, 244
, 438
).
 \end{dmath*}
Exponential growth rate:
\[\omega(W,S) =
1.77981392094460364656787163406
\dots\]
Cocompact? No\\
$(W,S)\in \CM$? No\\
Magma command: \texttt{HyperbolicCoxeterGraph(18)} or \texttt{HyperbolicCoxeterMatrix(18)} \\
\newpage
\subsubsection{\EHNC{5}}
\setcounter{mycnt}{5}
\[\Gamma(W,S) = \quad\csname CoxGrHNC\Roman{mycnt}\endcsname\]

\vskip12pt 
\[ M=\begin{pmatrix} 
1 & 6 & 2 & 3 \\
6 & 1 & 3 & 2 \\
2 & 3 & 1 & 4 \\
3 & 2 & 4 & 1
\end{pmatrix}\]
Numerator of $p_{(W,S)}(t)$:
\begin{dmath*} n_{(W,S)}(t) = 
1+3t+5t^{2}+7t^{3}+8t^{4}+8t^{5}+7t^{6}+5t^{7}+3t^{8}+t^{9}
\end{dmath*}
Denominator of $p_{(W,S)}(t)$:
\begin{dmath*} d_{(W,S)}(t) = 
1-t-t^{2}-t^{3}-t^{4}-t^{5}+t^{6}+t^{7}+t^{8}+t^{9}
\end{dmath*}
Initial values of $a_k$:
\begin{dmath*}\breakingcomma
 (a_k)_{k=0}^{9}=(1
, 4
, 10
, 22
, 45
, 90
, 177
, 344
, 666
, 1286
).
 \end{dmath*}
Exponential growth rate:
\[\omega(W,S) =
1.92756197548292530426190586173
\dots\]
Cocompact? No\\
$(W,S)\in \CM$? No\\
Magma command: \texttt{HyperbolicCoxeterGraph(19)} or \texttt{HyperbolicCoxeterMatrix(19)} \\
\newpage
\subsubsection{\EHNC{6}}
\setcounter{mycnt}{6}
\[\Gamma(W,S) = \quad\csname CoxGrHNC\Roman{mycnt}\endcsname\]

\vskip12pt 
\[ M=\begin{pmatrix} 
1 & 6 & 2 & 3 \\
6 & 1 & 3 & 2 \\
2 & 3 & 1 & 5 \\
3 & 2 & 5 & 1
\end{pmatrix}\]
Numerator of $p_{(W,S)}(t)$:
\begin{dmath*} n_{(W,S)}(t) = 
1+2t+2t^{2}+2t^{3}+2t^{4}+3t^{5}+3t^{6}+2t^{7}+2t^{8}+2t^{9}+2t^{10}+t^{11}
\end{dmath*}
Denominator of $p_{(W,S)}(t)$:
\begin{dmath*} d_{(W,S)}(t) = 
1-2t+t^{11}
\end{dmath*}
Initial values of $a_k$:
\begin{dmath*}\breakingcomma
 (a_k)_{k=0}^{11}=(1
, 4
, 10
, 22
, 46
, 95
, 193
, 388
, 778
, 1558
, 3118
, 6236
).
 \end{dmath*}
Exponential growth rate:
\[\omega(W,S) =
1.99901863271010113866340923913
\dots\]
Cocompact? No\\
$(W,S)\in \CM$? No\\
Magma command: \texttt{HyperbolicCoxeterGraph(20)} or \texttt{HyperbolicCoxeterMatrix(20)} \\
\newpage
% produced with RUNfortables.magma
\subsubsection{\EHNC{7}}
\setcounter{mycnt}{7}
\[\Gamma(W,S) = \quad\csname CoxGrHNC\Roman{mycnt}\endcsname\]

\vskip12pt 
\[ M=\begin{pmatrix} 
1 & 6 & 2 & 3 \\
6 & 1 & 3 & 2 \\
2 & 3 & 1 & 6 \\
3 & 2 & 6 & 1
\end{pmatrix}\]
Numerator of $p_{(W,S)}(t)$:
\begin{dmath*} n_{(W,S)}(t) = 
1+2t+2t^{2}+2t^{3}+2t^{4}+2t^{5}+t^{6}
\end{dmath*}
Denominator of $p_{(W,S)}(t)$:
\begin{dmath*} d_{(W,S)}(t) = 
1-2t-2t^{5}+3t^{6}
\end{dmath*}
Initial values of $a_k$:
\begin{dmath*}\breakingcomma
 (a_k)_{k=0}^{6}=(1
, 4
, 10
, 22
, 46
, 96
, 198
).
 \end{dmath*}
Exponential growth rate:
\[\omega(W,S) =
2.03073557275235254628366194872
\dots\]
Cocompact? No\\
$(W,S)\in \CM$? No\\
Magma command: \texttt{HyperbolicCoxeterGraph(21)} or \texttt{HyperbolicCoxeterMatrix(21)} \\
\newpage
\subsubsection{\EHNC{8}}
\setcounter{mycnt}{8}
\[\Gamma(W,S) = \quad\csname CoxGrHNC\Roman{mycnt}\endcsname\]

\vskip12pt 
\[ M=\begin{pmatrix} 
1 & 3 & 3 & 3 \\
3 & 1 & 3 & 2 \\
3 & 3 & 1 & 3 \\
3 & 2 & 3 & 1
\end{pmatrix}\]
Numerator of $p_{(W,S)}(t)$:
\begin{dmath*} n_{(W,S)}(t) = 
1+2t+3t^{2}+3t^{3}+2t^{4}+t^{5}
\end{dmath*}
Denominator of $p_{(W,S)}(t)$:
\begin{dmath*} d_{(W,S)}(t) = 
1-2t+t^{5}
\end{dmath*}
Initial values of $a_k$:
\begin{dmath*}\breakingcomma
 (a_k)_{k=0}^{5}=(1
, 4
, 11
, 25
, 52
, 104
).
 \end{dmath*}
Exponential growth rate:
\[\omega(W,S) =
1.92756197548292530426190586173
\dots\]
Cocompact? No\\
$(W,S)\in \CM$? No\\
Magma command: \texttt{HyperbolicCoxeterGraph(22)} or \texttt{HyperbolicCoxeterMatrix(22)} \\
\newpage
\subsubsection{\EHNC{9}}
\setcounter{mycnt}{9}
\[\Gamma(W,S) = \quad\csname CoxGrHNC\Roman{mycnt}\endcsname\]

\vskip12pt 
\[ M=\begin{pmatrix} 
1 & 3 & 3 & 2 \\
3 & 1 & 3 & 2 \\
3 & 3 & 1 & 3 \\
2 & 2 & 3 & 1
\end{pmatrix}\]
Numerator of $p_{(W,S)}(t)$:
\begin{dmath*} n_{(W,S)}(t) = 
1+2t+3t^{2}+3t^{3}+2t^{4}+t^{5}
\end{dmath*}
Denominator of $p_{(W,S)}(t)$:
\begin{dmath*} d_{(W,S)}(t) = 
1-2t+t^{2}-t^{3}+t^{4}
\end{dmath*}
Initial values of $a_k$:
\begin{dmath*}\breakingcomma
 (a_k)_{k=0}^{5}=(1
, 4
, 10
, 20
, 35
, 57
).
 \end{dmath*}
Exponential growth rate:
\[\omega(W,S) =
1.46557123187676802665673122522
\dots\]
Cocompact? No\\
$(W,S)\in \CM$? Yes\\
Magma command: \texttt{HyperbolicCoxeterGraph(23)} or \texttt{HyperbolicCoxeterMatrix(23)} \\
\newpage
\subsubsection{\EHNC{10}}
\setcounter{mycnt}{10}
\[\Gamma(W,S) = \quad\csname CoxGrHNC\Roman{mycnt}\endcsname\]

\vskip12pt 
\[ M=\begin{pmatrix} 
1 & 3 & 3 & 2 \\
3 & 1 & 3 & 2 \\
3 & 3 & 1 & 4 \\
2 & 2 & 4 & 1
\end{pmatrix}\]
Numerator of $p_{(W,S)}(t)$:
\begin{dmath*} n_{(W,S)}(t) = 
1+3t+5t^{2}+7t^{3}+8t^{4}+8t^{5}+7t^{6}+5t^{7}+3t^{8}+t^{9}
\end{dmath*}
Denominator of $p_{(W,S)}(t)$:
\begin{dmath*} d_{(W,S)}(t) = 
1-t-t^{2}-t^{4}+t^{6}+t^{8}
\end{dmath*}
Initial values of $a_k$:
\begin{dmath*}\breakingcomma
 (a_k)_{k=0}^{9}=(1
, 4
, 10
, 21
, 40
, 73
, 129
, 224
, 385
, 658
).
 \end{dmath*}
Exponential growth rate:
\[\omega(W,S) =
1.69783739529014608500737078192
\dots\]
Cocompact? No\\
$(W,S)\in \CM$? No\\
Magma command: \texttt{HyperbolicCoxeterGraph(24)} or \texttt{HyperbolicCoxeterMatrix(24)} \\
\newpage
\subsubsection{\EHNC{11}}
\setcounter{mycnt}{11}
\[\Gamma(W,S) = \quad\csname CoxGrHNC\Roman{mycnt}\endcsname\]

\vskip12pt 
\[ M=\begin{pmatrix} 
1 & 3 & 3 & 2 \\
3 & 1 & 3 & 2 \\
3 & 3 & 1 & 5 \\
2 & 2 & 5 & 1
\end{pmatrix}\]
Numerator of $p_{(W,S)}(t)$:
\begin{dmath*} n_{(W,S)}(t) = 
1+2t+2t^{2}+2t^{3}+2t^{4}+3t^{5}+3t^{6}+2t^{7}+2t^{8}+2t^{9}+2t^{10}+t^{11}
\end{dmath*}
Denominator of $p_{(W,S)}(t)$:
\begin{dmath*} d_{(W,S)}(t) = 
1-2t+t^{3}-t^{4}+t^{5}-t^{6}+t^{8}-t^{9}+t^{10}
\end{dmath*}
Initial values of $a_k$:
\begin{dmath*}\breakingcomma
 (a_k)_{k=0}^{11}=(1
, 4
, 10
, 21
, 41
, 78
, 145
, 266
, 485
, 882
, 1601
, 2902
).
 \end{dmath*}
Exponential growth rate:
\[\omega(W,S) =
1.80985205464161390068500549702
\dots\]
Cocompact? No\\
$(W,S)\in \CM$? No\\
Magma command: \texttt{HyperbolicCoxeterGraph(25)} or \texttt{HyperbolicCoxeterMatrix(25)} \\
\newpage
\subsubsection{\EHNC{12}}
\setcounter{mycnt}{12}
\[\Gamma(W,S) = \quad\csname CoxGrHNC\Roman{mycnt}\endcsname\]

\vskip12pt 
\[ M=\begin{pmatrix} 
1 & 3 & 3 & 2 \\
3 & 1 & 3 & 2 \\
3 & 3 & 1 & 6 \\
2 & 2 & 6 & 1
\end{pmatrix}\]
Numerator of $p_{(W,S)}(t)$:
\begin{dmath*} n_{(W,S)}(t) = 
1+2t+2t^{2}+2t^{3}+2t^{4}+2t^{5}+t^{6}
\end{dmath*}
Denominator of $p_{(W,S)}(t)$:
\begin{dmath*} d_{(W,S)}(t) = 
1-2t+t^{3}-t^{4}-t^{5}+2t^{6}
\end{dmath*}
Initial values of $a_k$:
\begin{dmath*}\breakingcomma
 (a_k)_{k=0}^{6}=(1
, 4
, 10
, 21
, 41
, 79
, 150
).
 \end{dmath*}
Exponential growth rate:
\[\omega(W,S) =
1.86007305043413707421477985317
\dots\]
Cocompact? No\\
$(W,S)\in \CM$? No\\
Magma command: \texttt{HyperbolicCoxeterGraph(26)} or \texttt{HyperbolicCoxeterMatrix(26)} \\
\newpage
\subsubsection{\EHNC{13}}
\setcounter{mycnt}{13}
\[\Gamma(W,S) = \quad\csname CoxGrHNC\Roman{mycnt}\endcsname\]

\vskip12pt 
\[ M=\begin{pmatrix} 
1 & 4 & 2 & 2 \\
4 & 1 & 4 & 2 \\
2 & 4 & 1 & 3 \\
2 & 2 & 3 & 1
\end{pmatrix}\]
Numerator of $p_{(W,S)}(t)$:
\begin{dmath*} n_{(W,S)}(t) = 
1+3t+5t^{2}+7t^{3}+8t^{4}+8t^{5}+7t^{6}+5t^{7}+3t^{8}+t^{9}
\end{dmath*}
Denominator of $p_{(W,S)}(t)$:
\begin{dmath*} d_{(W,S)}(t) = 
1-t-t^{3}+t^{8}
\end{dmath*}
Initial values of $a_k$:
\begin{dmath*}\breakingcomma
 (a_k)_{k=0}^{9}=(1
, 4
, 9
, 17
, 29
, 46
, 70
, 104
, 152
, 219
).
 \end{dmath*}
Exponential growth rate:
\[\omega(W,S) =
1.41216255165664341543382830907
\dots\]
Cocompact? No\\
$(W,S)\in \CM$? Yes\\
Magma command: \texttt{HyperbolicCoxeterGraph(27)} or \texttt{HyperbolicCoxeterMatrix(27)} \\
\newpage
\subsubsection{\EHNC{14}}
\setcounter{mycnt}{14}
\[\Gamma(W,S) = \quad\csname CoxGrHNC\Roman{mycnt}\endcsname\]

\vskip12pt 
\[ M=\begin{pmatrix} 
1 & 4 & 2 & 2 \\
4 & 1 & 4 & 2 \\
2 & 4 & 1 & 4 \\
2 & 2 & 4 & 1
\end{pmatrix}\]
Numerator of $p_{(W,S)}(t)$:
\begin{dmath*} n_{(W,S)}(t) = 
1+3t+4t^{2}+4t^{3}+3t^{4}+t^{5}
\end{dmath*}
Denominator of $p_{(W,S)}(t)$:
\begin{dmath*} d_{(W,S)}(t) = 
1-t-t^{2}-t^{3}+t^{4}+t^{5}
\end{dmath*}
Initial values of $a_k$:
\begin{dmath*}\breakingcomma
 (a_k)_{k=0}^{5}=(1
, 4
, 9
, 18
, 33
, 56
).
 \end{dmath*}
Exponential growth rate:
\[\omega(W,S) =
1.61803398874989484820458683437
\dots\]
Cocompact? No\\
$(W,S)\in \CM$? No\\
Magma command: \texttt{HyperbolicCoxeterGraph(28)} or \texttt{HyperbolicCoxeterMatrix(28)} \\
\newpage
\subsubsection{\EHNC{15}}
\setcounter{mycnt}{15}
\[\Gamma(W,S) = \quad\csname CoxGrHNC\Roman{mycnt}\endcsname\]

\vskip12pt 
\[ M=\begin{pmatrix} 
1 & 6 & 2 & 2 \\
6 & 1 & 3 & 2 \\
2 & 3 & 1 & 4 \\
2 & 2 & 4 & 1
\end{pmatrix}\]
Numerator of $p_{(W,S)}(t)$:
\begin{dmath*} n_{(W,S)}(t) = 
1+3t+5t^{2}+7t^{3}+8t^{4}+8t^{5}+7t^{6}+5t^{7}+3t^{8}+t^{9}
\end{dmath*}
Denominator of $p_{(W,S)}(t)$:
\begin{dmath*} d_{(W,S)}(t) = 
1-t-t^{3}-t^{5}+t^{6}+t^{8}
\end{dmath*}
Initial values of $a_k$:
\begin{dmath*}\breakingcomma
 (a_k)_{k=0}^{9}=(1
, 4
, 9
, 17
, 29
, 47
, 74
, 113
, 170
, 253
).
 \end{dmath*}
Exponential growth rate:
\[\omega(W,S) =
1.46557123187676802665673122522
\dots\]
Cocompact? No\\
$(W,S)\in \CM$? No\\
Magma command: \texttt{HyperbolicCoxeterGraph(29)} or \texttt{HyperbolicCoxeterMatrix(29)} \\
\newpage
\subsubsection{\EHNC{16}}
\setcounter{mycnt}{16}
\[\Gamma(W,S) = \quad\csname CoxGrHNC\Roman{mycnt}\endcsname\]

\vskip12pt 
\[ M=\begin{pmatrix} 
1 & 6 & 2 & 2 \\
6 & 1 & 3 & 2 \\
2 & 3 & 1 & 5 \\
2 & 2 & 5 & 1
\end{pmatrix}\]
Numerator of $p_{(W,S)}(t)$:
\begin{dmath*} n_{(W,S)}(t) = 
1+3t+5t^{2}+7t^{3}+9t^{4}+11t^{5}+12t^{6}+12t^{7}+12t^{8}+12t^{9}+11t^{10}+9t^{11}+7t^{12}+5t^{13}+3t^{14}+t^{15}
\end{dmath*}
Denominator of $p_{(W,S)}(t)$:
\begin{dmath*} d_{(W,S)}(t) = 
1-t-t^{3}-t^{5}+t^{12}+t^{14}
\end{dmath*}
Initial values of $a_k$:
\begin{dmath*}\breakingcomma
 (a_k)_{k=0}^{15}=(1
, 4
, 9
, 17
, 30
, 51
, 84
, 135
, 215
, 341
, 538
, 846
, 1328
, 2082
, 3262
, 5108
).
 \end{dmath*}
Exponential growth rate:
\[\omega(W,S) =
1.56488512068504857917661865311
\dots\]
Cocompact? No\\
$(W,S)\in \CM$? No\\
Magma command: \texttt{HyperbolicCoxeterGraph(30)} or \texttt{HyperbolicCoxeterMatrix(30)} \\
\newpage
\subsubsection{\EHNC{17}}
\setcounter{mycnt}{17}
\[\Gamma(W,S) = \quad\csname CoxGrHNC\Roman{mycnt}\endcsname\]

\vskip12pt 
\[ M=\begin{pmatrix} 
1 & 6 & 2 & 2 \\
6 & 1 & 3 & 2 \\
2 & 3 & 1 & 3 \\
2 & 2 & 3 & 1
\end{pmatrix}\]
Numerator of $p_{(W,S)}(t)$:
\begin{dmath*} n_{(W,S)}(t) = 
1+3t+5t^{2}+7t^{3}+8t^{4}+8t^{5}+7t^{6}+5t^{7}+3t^{8}+t^{9}
\end{dmath*}
Denominator of $p_{(W,S)}(t)$:
\begin{dmath*} d_{(W,S)}(t) = 
1-t-t^{4}+t^{8}
\end{dmath*}
Initial values of $a_k$:
\begin{dmath*}\breakingcomma
 (a_k)_{k=0}^{9}=(1
, 4
, 9
, 16
, 25
, 37
, 53
, 74
, 101
, 135
).
 \end{dmath*}
Exponential growth rate:
\[\omega(W,S) =
1.29646623878370094099709854474
\dots\]
Cocompact? No\\
$(W,S)\in \CM$? Yes\\
Magma command: \texttt{HyperbolicCoxeterGraph(31)} or \texttt{HyperbolicCoxeterMatrix(31)} \\
\newpage
\subsubsection{\EHNC{18}}
\setcounter{mycnt}{18}
\[\Gamma(W,S) = \quad\csname CoxGrHNC\Roman{mycnt}\endcsname\]

\vskip12pt 
\[ M=\begin{pmatrix} 
1 & 3 & 2 & 2 \\
3 & 1 & 6 & 2 \\
2 & 6 & 1 & 3 \\
2 & 2 & 3 & 1
\end{pmatrix}\]
Numerator of $p_{(W,S)}(t)$:
\begin{dmath*} n_{(W,S)}(t) = 
1+2t+2t^{2}+2t^{3}+2t^{4}+2t^{5}+t^{6}
\end{dmath*}
Denominator of $p_{(W,S)}(t)$:
\begin{dmath*} d_{(W,S)}(t) = 
1-2t+t^{2}-t^{4}+t^{6}
\end{dmath*}
Initial values of $a_k$:
\begin{dmath*}\breakingcomma
 (a_k)_{k=0}^{6}=(1
, 4
, 9
, 16
, 26
, 42
, 67
).
 \end{dmath*}
Exponential growth rate:
\[\omega(W,S) =
1.49709404876279664895121309733
\dots\]
Cocompact? No\\
$(W,S)\in \CM$? No\\
Magma command: \texttt{HyperbolicCoxeterGraph(32)} or \texttt{HyperbolicCoxeterMatrix(32)} \\
\newpage
\subsubsection{\EHNC{19}}
\setcounter{mycnt}{19}
\[\Gamma(W,S) = \quad\csname CoxGrHNC\Roman{mycnt}\endcsname\]

\vskip12pt 
\[ M=\begin{pmatrix} 
1 & 6 & 2 & 2 \\
6 & 1 & 3 & 2 \\
2 & 3 & 1 & 6 \\
2 & 2 & 6 & 1
\end{pmatrix}\]
Numerator of $p_{(W,S)}(t)$:
\begin{dmath*} n_{(W,S)}(t) = 
1+3t+4t^{2}+4t^{3}+4t^{4}+4t^{5}+3t^{6}+t^{7}
\end{dmath*}
Denominator of $p_{(W,S)}(t)$:
\begin{dmath*} d_{(W,S)}(t) = 
1-t-t^{2}-t^{5}+t^{6}+t^{7}
\end{dmath*}
Initial values of $a_k$:
\begin{dmath*}\breakingcomma
 (a_k)_{k=0}^{7}=(1
, 4
, 9
, 17
, 30
, 52
, 88
, 145
).
 \end{dmath*}
Exponential growth rate:
\[\omega(W,S) =
1.61803398874989484820458683437
\dots\]
Cocompact? No\\
$(W,S)\in \CM$? No\\
Magma command: \texttt{HyperbolicCoxeterGraph(33)} or \texttt{HyperbolicCoxeterMatrix(33)} \\
\newpage
\subsubsection{\EHNC{20}}
\setcounter{mycnt}{20}
\[\Gamma(W,S) = \quad\csname CoxGrHNC\Roman{mycnt}\endcsname\]

\vskip12pt 
\[ M=\begin{pmatrix} 
1 & 6 & 2 & 2 \\
6 & 1 & 3 & 3 \\
2 & 3 & 1 & 2 \\
2 & 3 & 2 & 1
\end{pmatrix}\]
Numerator of $p_{(W,S)}(t)$:
\begin{dmath*} n_{(W,S)}(t) = 
1+3t+5t^{2}+7t^{3}+8t^{4}+8t^{5}+7t^{6}+5t^{7}+3t^{8}+t^{9}
\end{dmath*}
Denominator of $p_{(W,S)}(t)$:
\begin{dmath*} d_{(W,S)}(t) = 
1-t-t^{3}-2t^{5}+t^{6}+t^{8}+t^{9}
\end{dmath*}
Initial values of $a_k$:
\begin{dmath*}\breakingcomma
 (a_k)_{k=0}^{9}=(1
, 4
, 9
, 17
, 29
, 48
, 79
, 127
, 202
, 318
).
 \end{dmath*}
Exponential growth rate:
\[\omega(W,S) =
1.56435053512350983063450400972
\dots\]
Cocompact? No\\
$(W,S)\in \CM$? No\\
Magma command: \texttt{HyperbolicCoxeterGraph(34)} or \texttt{HyperbolicCoxeterMatrix(34)} \\
\newpage
\subsubsection{\EHNC{21}}
\setcounter{mycnt}{21}
\[\Gamma(W,S) = \quad\csname CoxGrHNC\Roman{mycnt}\endcsname\]

\vskip12pt 
\[ M=\begin{pmatrix} 
1 & 4 & 2 & 2 \\
4 & 1 & 4 & 3 \\
2 & 4 & 1 & 2 \\
2 & 3 & 2 & 1
\end{pmatrix}\]
Numerator of $p_{(W,S)}(t)$:
\begin{dmath*} n_{(W,S)}(t) = 
1+2t+2t^{2}+3t^{3}+3t^{4}+2t^{5}+2t^{6}+t^{7}
\end{dmath*}
Denominator of $p_{(W,S)}(t)$:
\begin{dmath*} d_{(W,S)}(t) = 
1-2t+t^{2}-t^{3}+t^{4}-t^{5}+t^{6}
\end{dmath*}
Initial values of $a_k$:
\begin{dmath*}\breakingcomma
 (a_k)_{k=0}^{7}=(1
, 4
, 9
, 18
, 33
, 56
, 93
, 151
).
 \end{dmath*}
Exponential growth rate:
\[\omega(W,S) =
1.57014731219605436291066543514
\dots\]
Cocompact? No\\
$(W,S)\in \CM$? Yes\\
Magma command: \texttt{HyperbolicCoxeterGraph(35)} or \texttt{HyperbolicCoxeterMatrix(35)} \\
\newpage
\subsubsection{\EHNC{22}}
\setcounter{mycnt}{22}
\[\Gamma(W,S) = \quad\csname CoxGrHNC\Roman{mycnt}\endcsname\]

\vskip12pt 
\[ M=\begin{pmatrix} 
1 & 4 & 2 & 2 \\
4 & 1 & 4 & 4 \\
2 & 4 & 1 & 2 \\
2 & 4 & 2 & 1
\end{pmatrix}\]
Numerator of $p_{(W,S)}(t)$:
\begin{dmath*} n_{(W,S)}(t) = 
1+3t+4t^{2}+4t^{3}+3t^{4}+t^{5}
\end{dmath*}
Denominator of $p_{(W,S)}(t)$:
\begin{dmath*} d_{(W,S)}(t) = 
1-t-t^{2}-2t^{3}+t^{4}+2t^{5}
\end{dmath*}
Initial values of $a_k$:
\begin{dmath*}\breakingcomma
 (a_k)_{k=0}^{5}=(1
, 4
, 9
, 19
, 38
, 70
).
 \end{dmath*}
Exponential growth rate:
\[\omega(W,S) =
1.80843400251150161450701496589
\dots\]
Cocompact? No\\
$(W,S)\in \CM$? No\\
Magma command: \texttt{HyperbolicCoxeterGraph(36)} or \texttt{HyperbolicCoxeterMatrix(36)} \\
\newpage
\subsubsection{\EHNC{23}}
\setcounter{mycnt}{23}
\[\Gamma(W,S) = \quad\csname CoxGrHNC\Roman{mycnt}\endcsname\]

\vskip12pt 
\[ M=\begin{pmatrix} 
1 & 3 & 3 & 3 \\
3 & 1 & 3 & 3 \\
3 & 3 & 1 & 3 \\
3 & 3 & 3 & 1
\end{pmatrix}\]
Numerator of $p_{(W,S)}(t)$:
\begin{dmath*} n_{(W,S)}(t) = 
1+2t+2t^{2}+t^{3}
\end{dmath*}
Denominator of $p_{(W,S)}(t)$:
\begin{dmath*} d_{(W,S)}(t) = 
1-2t-2t^{2}+3t^{3}
\end{dmath*}
Initial values of $a_k$:
\begin{dmath*}\breakingcomma
 (a_k)_{k=0}^{3}=(1
, 4
, 12
, 30
).
 \end{dmath*}
Exponential growth rate:
\[\omega(W,S) =
2.30277563773199464655961063374
\dots\]
Cocompact? No\\
$(W,S)\in \CM$? No\\
Magma command: \texttt{HyperbolicCoxeterGraph(37)} or \texttt{HyperbolicCoxeterMatrix(37)} \\
\newpage
\subsubsection{\EHNC{24}}
\setcounter{mycnt}{24}
\[\Gamma(W,S) = \quad\csname CoxGrHNC\Roman{mycnt}\endcsname\]

\vskip12pt 
\[ M=\begin{pmatrix} 
1 & 4 & 2 & 2 & 2 \\
4 & 1 & 3 & 2 & 2 \\
2 & 3 & 1 & 4 & 2 \\
2 & 2 & 4 & 1 & 3 \\
2 & 2 & 2 & 3 & 1
\end{pmatrix}\]
Numerator of $p_{(W,S)}(t)$:
\begin{dmath*} n_{(W,S)}(t) = 
-1-4t-9t^{2}-16t^{3}-25t^{4}-36t^{5}-48t^{6}-60t^{7}-71t^{8}-80t^{9}-87t^{10}-92t^{11}-94t^{12}-92t^{13}-87t^{14}-80t^{15}-71t^{16}-60t^{17}-48t^{18}-36t^{19}-25t^{20}-16t^{21}-9t^{22}-4t^{23}-t^{24}
\end{dmath*}
Denominator of $p_{(W,S)}(t)$:
\begin{dmath*} d_{(W,S)}(t) = 
-1+t+t^{3}-t^{4}+2t^{5}-t^{6}+t^{7}-2t^{8}+2t^{9}-2t^{10}+t^{11}-2t^{12}+t^{13}-2t^{14}+t^{15}-t^{16}-t^{18}+t^{19}+t^{23}
\end{dmath*}
Initial values of $a_k$:
\begin{dmath*}\breakingcomma
 (a_k)_{k=0}^{24}=(1
, 5
, 14
, 31
, 60
, 107
, 181
, 294
, 463
, 712
, 1077
, 1610
, 2384
, 3503
, 5116
, 7438
, 10777
, 15573
, 22454
, 32319
, 46457
, 66713
, 95726
, 137270
, 196745
).
 \end{dmath*}
Exponential growth rate:
\[\omega(W,S) =
1.43070541706112540116746869789
\dots\]
Cocompact? No\\
$(W,S)\in \CM$? Yes\\
Magma command: \texttt{HyperbolicCoxeterGraph(38)} or \texttt{HyperbolicCoxeterMatrix(38)} \\
\newpage
\subsubsection{\EHNC{25}}
\setcounter{mycnt}{25}
\[\Gamma(W,S) = \quad\csname CoxGrHNC\Roman{mycnt}\endcsname\]

\vskip12pt 
\[ M=\begin{pmatrix} 
1 & 4 & 2 & 2 & 2 \\
4 & 1 & 3 & 3 & 2 \\
2 & 3 & 1 & 2 & 2 \\
2 & 3 & 2 & 1 & 3 \\
2 & 2 & 2 & 3 & 1
\end{pmatrix}\]
Numerator of $p_{(W,S)}(t)$:
\begin{dmath*} n_{(W,S)}(t) = 
-1-5t-14t^{2}-30t^{3}-54t^{4}-85t^{5}-120t^{6}-155t^{7}-185t^{8}-205t^{9}-212t^{10}-205t^{11}-185t^{12}-155t^{13}-120t^{14}-85t^{15}-54t^{16}-30t^{17}-14t^{18}-5t^{19}-t^{20}
\end{dmath*}
Denominator of $p_{(W,S)}(t)$:
\begin{dmath*} d_{(W,S)}(t) = 
-1+t^{3}+t^{4}+2t^{5}+t^{6}+t^{7}-t^{8}-t^{9}-2t^{10}-2t^{11}-2t^{12}-t^{13}-t^{14}+t^{16}+t^{17}+t^{18}+t^{19}
\end{dmath*}
Initial values of $a_k$:
\begin{dmath*}\breakingcomma
 (a_k)_{k=0}^{20}=(1
, 5
, 14
, 31
, 60
, 106
, 176
, 280
, 431
, 646
, 950
, 1377
, 1973
, 2802
, 3952
, 5543
, 7740
, 10771
, 14947
, 20695
, 28602
).
 \end{dmath*}
Exponential growth rate:
\[\omega(W,S) =
1.37176048915656028424231212289
\dots\]
Cocompact? No\\
$(W,S)\in \CM$? Yes\\
Magma command: \texttt{HyperbolicCoxeterGraph(39)} or \texttt{HyperbolicCoxeterMatrix(39)} \\
\newpage
\subsubsection{\EHNC{26}}
\setcounter{mycnt}{26}
\[\Gamma(W,S) = \quad\csname CoxGrHNC\Roman{mycnt}\endcsname\]

\vskip12pt 
\[ M=\begin{pmatrix} 
1 & 4 & 2 & 2 & 3 \\
4 & 1 & 3 & 3 & 2 \\
2 & 3 & 1 & 2 & 2 \\
2 & 3 & 2 & 1 & 2 \\
3 & 2 & 2 & 2 & 1
\end{pmatrix}\]
Numerator of $p_{(W,S)}(t)$:
\begin{dmath*} n_{(W,S)}(t) = 
-1-4t-9t^{2}-16t^{3}-25t^{4}-36t^{5}-48t^{6}-60t^{7}-71t^{8}-80t^{9}-87t^{10}-92t^{11}-94t^{12}-92t^{13}-87t^{14}-80t^{15}-71t^{16}-60t^{17}-48t^{18}-36t^{19}-25t^{20}-16t^{21}-9t^{22}-4t^{23}-t^{24}
\end{dmath*}
Denominator of $p_{(W,S)}(t)$:
\begin{dmath*} d_{(W,S)}(t) = 
-1+t+t^{3}+t^{5}-t^{8}-t^{11}-t^{12}-t^{13}-t^{15}-t^{17}+t^{20}+t^{23}
\end{dmath*}
Initial values of $a_k$:
\begin{dmath*}\breakingcomma
 (a_k)_{k=0}^{24}=(1
, 5
, 14
, 31
, 61
, 112
, 196
, 331
, 544
, 876
, 1392
, 2192
, 3426
, 5322
, 8231
, 12691
, 19523
, 29980
, 45974
, 70426
, 107803
, 164930
, 252225
, 385597
, 589350
).
 \end{dmath*}
Exponential growth rate:
\[\omega(W,S) =
1.52732217619767494331020421320
\dots\]
Cocompact? No\\
$(W,S)\in \CM$? Yes\\
Magma command: \texttt{HyperbolicCoxeterGraph(40)} or \texttt{HyperbolicCoxeterMatrix(40)} \\
\newpage
\subsubsection{\EHNC{27}}
\setcounter{mycnt}{27}
\[\Gamma(W,S) = \quad\csname CoxGrHNC\Roman{mycnt}\endcsname\]

\vskip12pt 
\[ M=\begin{pmatrix} 
1 & 4 & 2 & 2 & 2 \\
4 & 1 & 3 & 3 & 2 \\
2 & 3 & 1 & 2 & 2 \\
2 & 3 & 2 & 1 & 4 \\
2 & 2 & 2 & 4 & 1
\end{pmatrix}\]
Numerator of $p_{(W,S)}(t)$:
\begin{dmath*} n_{(W,S)}(t) = 
-1-4t-9t^{2}-16t^{3}-24t^{4}-32t^{5}-39t^{6}-44t^{7}-46t^{8}-44t^{9}-39t^{10}-32t^{11}-24t^{12}-16t^{13}-9t^{14}-4t^{15}-t^{16}
\end{dmath*}
Denominator of $p_{(W,S)}(t)$:
\begin{dmath*} d_{(W,S)}(t) = 
-1+t+2t^{3}-t^{4}+2t^{5}-2t^{6}+t^{7}-3t^{8}+t^{9}-2t^{10}-t^{12}+t^{15}+t^{16}
\end{dmath*}
Initial values of $a_k$:
\begin{dmath*}\breakingcomma
 (a_k)_{k=0}^{16}=(1
, 5
, 14
, 32
, 65
, 122
, 219
, 380
, 643
, 1069
, 1756
, 2861
, 4632
, 7466
, 11994
, 19222
, 30756
).
 \end{dmath*}
Exponential growth rate:
\[\omega(W,S) =
1.59270121882598805544129340934
\dots\]
Cocompact? No\\
$(W,S)\in \CM$? No\\
Magma command: \texttt{HyperbolicCoxeterGraph(41)} or \texttt{HyperbolicCoxeterMatrix(41)} \\
\newpage
\subsubsection{\EHNC{28}}
\setcounter{mycnt}{28}
\[\Gamma(W,S) = \quad\csname CoxGrHNC\Roman{mycnt}\endcsname\]

\vskip12pt 
\[ M=\begin{pmatrix} 
1 & 3 & 2 & 3 & 2 \\
3 & 1 & 3 & 2 & 2 \\
2 & 3 & 1 & 3 & 2 \\
3 & 2 & 3 & 1 & 3 \\
2 & 2 & 2 & 3 & 1
\end{pmatrix}\]
Numerator of $p_{(W,S)}(t)$:
\begin{dmath*} n_{(W,S)}(t) = 
-1-5t-14t^{2}-30t^{3}-53t^{4}-80t^{5}-106t^{6}-125t^{7}-132t^{8}-125t^{9}-106t^{10}-80t^{11}-53t^{12}-30t^{13}-14t^{14}-5t^{15}-t^{16}
\end{dmath*}
Denominator of $p_{(W,S)}(t)$:
\begin{dmath*} d_{(W,S)}(t) = 
-1+t^{2}+t^{3}+2t^{4}+t^{5}-t^{7}-3t^{8}-2t^{9}-2t^{10}-t^{11}+t^{12}+t^{13}+t^{14}+t^{15}
\end{dmath*}
Initial values of $a_k$:
\begin{dmath*}\breakingcomma
 (a_k)_{k=0}^{16}=(1
, 5
, 15
, 36
, 75
, 142
, 252
, 428
, 704
, 1132
, 1791
, 2800
, 4339
, 6680
, 10234
, 15621
, 23778
).
 \end{dmath*}
Exponential growth rate:
\[\omega(W,S) =
1.50928431170332633862815772930
\dots\]
Cocompact? No\\
$(W,S)\in \CM$? Yes\\
Magma command: \texttt{HyperbolicCoxeterGraph(42)} or \texttt{HyperbolicCoxeterMatrix(42)} \\
\newpage
\subsubsection{\EHNC{29}}
\setcounter{mycnt}{29}
\[\Gamma(W,S) = \quad\csname CoxGrHNC\Roman{mycnt}\endcsname\]

\vskip12pt 
\[ M=\begin{pmatrix} 
1 & 3 & 2 & 3 & 2 \\
3 & 1 & 3 & 2 & 2 \\
2 & 3 & 1 & 3 & 2 \\
3 & 2 & 3 & 1 & 4 \\
2 & 2 & 2 & 4 & 1
\end{pmatrix}\]
Numerator of $p_{(W,S)}(t)$:
\begin{dmath*} n_{(W,S)}(t) = 
-1-4t-9t^{2}-16t^{3}-24t^{4}-32t^{5}-39t^{6}-44t^{7}-46t^{8}-44t^{9}-39t^{10}-32t^{11}-24t^{12}-16t^{13}-9t^{14}-4t^{15}-t^{16}
\end{dmath*}
Denominator of $p_{(W,S)}(t)$:
\begin{dmath*} d_{(W,S)}(t) = 
-1+t+t^{2}+t^{3}-t^{7}-2t^{8}-t^{9}-t^{10}-t^{11}+t^{15}+t^{16}
\end{dmath*}
Initial values of $a_k$:
\begin{dmath*}\breakingcomma
 (a_k)_{k=0}^{16}=(1
, 5
, 15
, 37
, 81
, 165
, 322
, 611
, 1137
, 2088
, 3802
, 6883
, 12413
, 22329
, 40096
, 71916
, 128888
).
 \end{dmath*}
Exponential growth rate:
\[\omega(W,S) =
1.78932392876860973173394906469
\dots\]
Cocompact? No\\
$(W,S)\in \CM$? No\\
Magma command: \texttt{HyperbolicCoxeterGraph(43)} or \texttt{HyperbolicCoxeterMatrix(43)} \\
\newpage
\subsubsection{\EHNC{30}}
\setcounter{mycnt}{30}
\[\Gamma(W,S) = \quad\csname CoxGrHNC\Roman{mycnt}\endcsname\]

\vskip12pt 
\[ M=\begin{pmatrix} 
1 & 4 & 2 & 2 & 2 \\
4 & 1 & 3 & 3 & 3 \\
2 & 3 & 1 & 2 & 2 \\
2 & 3 & 2 & 1 & 2 \\
2 & 3 & 2 & 2 & 1
\end{pmatrix}\]
Numerator of $p_{(W,S)}(t)$:
\begin{dmath*} n_{(W,S)}(t) = 
-1-4t-9t^{2}-16t^{3}-23t^{4}-28t^{5}-30t^{6}-28t^{7}-23t^{8}-16t^{9}-9t^{10}-4t^{11}-t^{12}
\end{dmath*}
Denominator of $p_{(W,S)}(t)$:
\begin{dmath*} d_{(W,S)}(t) = 
-1+t+3t^{3}-t^{4}+2t^{5}-3t^{6}-t^{7}-4t^{8}+t^{11}+2t^{12}
\end{dmath*}
Initial values of $a_k$:
\begin{dmath*}\breakingcomma
 (a_k)_{k=0}^{12}=(1
, 5
, 14
, 33
, 70
, 137
, 259
, 476
, 855
, 1518
, 2671
, 4670
, 8135
).
 \end{dmath*}
Exponential growth rate:
\[\omega(W,S) =
1.72583201939633721748849591766
\dots\]
Cocompact? No\\
$(W,S)\in \CM$? Yes\\
Magma command: \texttt{HyperbolicCoxeterGraph(44)} or \texttt{HyperbolicCoxeterMatrix(44)} \\
\newpage
\subsubsection{\EHNC{31}}
\setcounter{mycnt}{31}
\[\Gamma(W,S) = \quad\csname CoxGrHNC\Roman{mycnt}\endcsname\]

\vskip12pt 
\[ M=\begin{pmatrix} 
1 & 3 & 2 & 3 & 2 \\
3 & 1 & 3 & 2 & 3 \\
2 & 3 & 1 & 3 & 2 \\
3 & 2 & 3 & 1 & 3 \\
2 & 3 & 2 & 3 & 1
\end{pmatrix}\]
Numerator of $p_{(W,S)}(t)$:
\begin{dmath*} n_{(W,S)}(t) = 
-1-4t-9t^{2}-16t^{3}-23t^{4}-28t^{5}-30t^{6}-28t^{7}-23t^{8}-16t^{9}-9t^{10}-4t^{11}-t^{12}
\end{dmath*}
Denominator of $p_{(W,S)}(t)$:
\begin{dmath*} d_{(W,S)}(t) = 
-1+t+2t^{2}+t^{3}-2t^{5}-t^{6}-2t^{7}-3t^{8}+t^{9}-t^{10}+t^{11}+2t^{12}
\end{dmath*}
Initial values of $a_k$:
\begin{dmath*}\breakingcomma
 (a_k)_{k=0}^{12}=(1
, 5
, 16
, 43
, 103
, 231
, 499
, 1053
, 2190
, 4516
, 9263
, 18937
, 38638
).
 \end{dmath*}
Exponential growth rate:
\[\omega(W,S) =
2.03378256690775761901140546096
\dots\]
Cocompact? No\\
$(W,S)\in \CM$? No\\
Magma command: \texttt{HyperbolicCoxeterGraph(45)} or \texttt{HyperbolicCoxeterMatrix(45)} \\
\newpage
\subsubsection{\EHNC{32}}
\setcounter{mycnt}{32}
\[\Gamma(W,S) = \quad\csname CoxGrHNC\Roman{mycnt}\endcsname\]

\vskip12pt 
\[ M=\begin{pmatrix} 
1 & 4 & 2 & 2 & 3 \\
4 & 1 & 3 & 2 & 2 \\
2 & 3 & 1 & 4 & 2 \\
2 & 2 & 4 & 1 & 3 \\
3 & 2 & 2 & 3 & 1
\end{pmatrix}\]
Numerator of $p_{(W,S)}(t)$:
\begin{dmath*} n_{(W,S)}(t) = 
-1-4t-9t^{2}-16t^{3}-25t^{4}-36t^{5}-48t^{6}-60t^{7}-71t^{8}-80t^{9}-87t^{10}-92t^{11}-94t^{12}-92t^{13}-87t^{14}-80t^{15}-71t^{16}-60t^{17}-48t^{18}-36t^{19}-25t^{20}-16t^{21}-9t^{22}-4t^{23}-t^{24}
\end{dmath*}
Denominator of $p_{(W,S)}(t)$:
\begin{dmath*} d_{(W,S)}(t) = 
-1+t+t^{2}+t^{3}+t^{5}-3t^{8}-t^{10}-t^{11}-3t^{12}-t^{13}-t^{14}-t^{15}-2t^{16}-t^{17}+t^{20}+t^{22}+t^{23}
\end{dmath*}
Initial values of $a_k$:
\begin{dmath*}\breakingcomma
 (a_k)_{k=0}^{24}=(1
, 5
, 15
, 37
, 82
, 171
, 343
, 671
, 1290
, 2451
, 4624
, 8683
, 16254
, 30362
, 56638
, 105560
, 196625
, 366112
, 681520
, 1268443
, 2360568
, 4392703
, 8173859
, 15209302
, 28299764
).
 \end{dmath*}
Exponential growth rate:
\[\omega(W,S) =
1.86061798901663192520827350833
\dots\]
Cocompact? No\\
$(W,S)\in \CM$? No\\
Magma command: \texttt{HyperbolicCoxeterGraph(46)} or \texttt{HyperbolicCoxeterMatrix(46)} \\
\newpage
% produced with RUNfortables.magma
\subsubsection{\EHNC{33}}
\setcounter{mycnt}{33}
\[\Gamma(W,S) = \quad\csname CoxGrHNC\Roman{mycnt}\endcsname\]

\vskip12pt 
\[ M=\begin{pmatrix} 
1 & 3 & 2 & 2 & 2 & 2 \\
3 & 1 & 4 & 2 & 2 & 2 \\
2 & 4 & 1 & 3 & 2 & 2 \\
2 & 2 & 3 & 1 & 3 & 2 \\
2 & 2 & 2 & 3 & 1 & 3 \\
2 & 2 & 2 & 2 & 3 & 1
\end{pmatrix}\]
Numerator of $p_{(W,S)}(t)$:
\begin{dmath*} n_{(W,S)}(t) = 
1+5t+14t^{2}+30t^{3}+55t^{4}+91t^{5}+139t^{6}+199t^{7}+270t^{8}+350t^{9}+436t^{10}+524t^{11}+609t^{12}+685t^{13}+747t^{14}+791t^{15}+814t^{16}+814t^{17}+791t^{18}+747t^{19}+685t^{20}+609t^{21}+524t^{22}+436t^{23}+350t^{24}+270t^{25}+199t^{26}+139t^{27}+91t^{28}+55t^{29}+30t^{30}+14t^{31}+5t^{32}+t^{33}
\end{dmath*}
Denominator of $p_{(W,S)}(t)$:
\begin{dmath*} d_{(W,S)}(t) = 
1-t-t^{5}+t^{12}+t^{14}+t^{16}+t^{18}-t^{19}-t^{22}-t^{23}-t^{26}+t^{32}
\end{dmath*}
Initial values of $a_k$:
\begin{dmath*}\breakingcomma
 (a_k)_{k=0}^{33}=(1
, 6
, 20
, 50
, 105
, 197
, 342
, 561
, 881
, 1336
, 1969
, 2835
, 4004
, 5564
, 7626
, 10330
, 13853
, 18418
, 24305
, 31865
, 41538
, 53876
, 69572
, 89496
, 114738
, 146662
, 186974
, 237809
, 301839
, 382407
, 483693
, 610922
, 770624
, 970960
).
 \end{dmath*}
Exponential growth rate:
\[\omega(W,S) =
1.24813167325321011601939347627
\dots\]
Cocompact? No\\
$(W,S)\in \CM$? Yes\\
Magma command: \texttt{HyperbolicCoxeterGraph(47)} or \texttt{HyperbolicCoxeterMatrix(47)} \\
\newpage
\subsubsection{\EHNC{34}}
\setcounter{mycnt}{34}
\[\Gamma(W,S) = \quad\csname CoxGrHNC\Roman{mycnt}\endcsname\]

\vskip12pt 
\[ M=\begin{pmatrix} 
1 & 3 & 2 & 2 & 2 & 2 \\
3 & 1 & 3 & 2 & 2 & 2 \\
2 & 3 & 1 & 4 & 2 & 2 \\
2 & 2 & 4 & 1 & 3 & 2 \\
2 & 2 & 2 & 3 & 1 & 3 \\
2 & 2 & 2 & 2 & 3 & 1
\end{pmatrix}\]
Numerator of $p_{(W,S)}(t)$:
\begin{dmath*} n_{(W,S)}(t) = 
1+5t+13t^{2}+25t^{3}+41t^{4}+61t^{5}+84t^{6}+108t^{7}+131t^{8}+151t^{9}+167t^{10}+179t^{11}+186t^{12}+186t^{13}+179t^{14}+167t^{15}+151t^{16}+131t^{17}+108t^{18}+84t^{19}+61t^{20}+41t^{21}+25t^{22}+13t^{23}+5t^{24}+t^{25}
\end{dmath*}
Denominator of $p_{(W,S)}(t)$:
\begin{dmath*} d_{(W,S)}(t) = 
1-t-t^{2}+t^{3}-2t^{5}+t^{6}+t^{7}+t^{10}-t^{11}+t^{12}+t^{13}-t^{14}-t^{15}+2t^{16}-t^{18}-t^{19}-t^{22}-t^{23}+t^{24}+t^{25}
\end{dmath*}
Initial values of $a_k$:
\begin{dmath*}\breakingcomma
 (a_k)_{k=0}^{25}=(1
, 6
, 20
, 50
, 105
, 198
, 348
, 582
, 937
, 1462
, 2224
, 3316
, 4867
, 7052
, 10107
, 14352
, 20223
, 28314
, 39429
, 54654
, 75456
, 103822
, 142444
, 194964
, 266301
, 363096
).
 \end{dmath*}
Exponential growth rate:
\[\omega(W,S) =
1.34911552067238295535561274062
\dots\]
Cocompact? No\\
$(W,S)\in \CM$? Yes\\
Magma command: \texttt{HyperbolicCoxeterGraph(48)} or \texttt{HyperbolicCoxeterMatrix(48)} \\
\newpage
\subsubsection{\EHNC{35}}
\setcounter{mycnt}{35}
\[\Gamma(W,S) = \quad\csname CoxGrHNC\Roman{mycnt}\endcsname\]

\vskip12pt 
\[ M=\begin{pmatrix} 
1 & 3 & 2 & 2 & 2 & 2 \\
3 & 1 & 4 & 2 & 2 & 2 \\
2 & 4 & 1 & 3 & 2 & 2 \\
2 & 2 & 3 & 1 & 3 & 2 \\
2 & 2 & 2 & 3 & 1 & 4 \\
2 & 2 & 2 & 2 & 4 & 1
\end{pmatrix}\]
Numerator of $p_{(W,S)}(t)$:
\begin{dmath*} n_{(W,S)}(t) = 
1+5t+13t^{2}+25t^{3}+41t^{4}+61t^{5}+84t^{6}+108t^{7}+131t^{8}+151t^{9}+167t^{10}+179t^{11}+186t^{12}+186t^{13}+179t^{14}+167t^{15}+151t^{16}+131t^{17}+108t^{18}+84t^{19}+61t^{20}+41t^{21}+25t^{22}+13t^{23}+5t^{24}+t^{25}
\end{dmath*}
Denominator of $p_{(W,S)}(t)$:
\begin{dmath*} d_{(W,S)}(t) = 
1-t-t^{2}+t^{4}-t^{5}+2t^{8}-t^{11}+2t^{12}-t^{14}-t^{15}+2t^{16}-t^{18}-t^{19}-t^{22}-t^{23}+t^{24}+t^{25}
\end{dmath*}
Initial values of $a_k$:
\begin{dmath*}\breakingcomma
 (a_k)_{k=0}^{25}=(1
, 6
, 20
, 51
, 111
, 218
, 399
, 694
, 1162
, 1888
, 2996
, 4667
, 7163
, 10862
, 16310
, 24296
, 35957
, 52931
, 77575
, 113278
, 164912
, 239481
, 347050
, 502076
, 725321
, 1046597
).
 \end{dmath*}
Exponential growth rate:
\[\omega(W,S) =
1.43464552105983588116454981990
\dots\]
Cocompact? No\\
$(W,S)\in \CM$? No\\
Magma command: \texttt{HyperbolicCoxeterGraph(49)} or \texttt{HyperbolicCoxeterMatrix(49)} \\
\newpage
\subsubsection{\EHNC{36}}
\setcounter{mycnt}{36}
\[\Gamma(W,S) = \quad\csname CoxGrHNC\Roman{mycnt}\endcsname\]

\vskip12pt 
\[ M=\begin{pmatrix} 
1 & 3 & 2 & 2 & 2 & 2 \\
3 & 1 & 4 & 2 & 2 & 2 \\
2 & 4 & 1 & 3 & 2 & 2 \\
2 & 2 & 3 & 1 & 3 & 3 \\
2 & 2 & 2 & 3 & 1 & 2 \\
2 & 2 & 2 & 3 & 2 & 1
\end{pmatrix}\]
Numerator of $p_{(W,S)}(t)$:
\begin{dmath*} n_{(W,S)}(t) = 
1+5t+13t^{2}+25t^{3}+41t^{4}+61t^{5}+84t^{6}+108t^{7}+131t^{8}+151t^{9}+167t^{10}+179t^{11}+186t^{12}+186t^{13}+179t^{14}+167t^{15}+151t^{16}+131t^{17}+108t^{18}+84t^{19}+61t^{20}+41t^{21}+25t^{22}+13t^{23}+5t^{24}+t^{25}
\end{dmath*}
Denominator of $p_{(W,S)}(t)$:
\begin{dmath*} d_{(W,S)}(t) = 
1-t-t^{2}+t^{4}-2t^{5}+t^{6}+2t^{8}+t^{10}-2t^{11}+3t^{12}-t^{14}-2t^{15}+3t^{16}-t^{18}-2t^{19}-t^{22}-2t^{23}+t^{24}+2t^{25}
\end{dmath*}
Initial values of $a_k$:
\begin{dmath*}\breakingcomma
 (a_k)_{k=0}^{25}=(1
, 6
, 20
, 51
, 111
, 219
, 405
, 715
, 1220
, 2026
, 3294
, 5269
, 8321
, 13008
, 20168
, 31064
, 47595
, 72614
, 110405
, 167394
, 253223
, 382354
, 576478
, 868116
, 1306014
, 1963236
).
 \end{dmath*}
Exponential growth rate:
\[\omega(W,S) =
1.49807706548831508637056306574
\dots\]
Cocompact? No\\
$(W,S)\in \CM$? Yes\\
Magma command: \texttt{HyperbolicCoxeterGraph(50)} or \texttt{HyperbolicCoxeterMatrix(50)} \\
\newpage
\subsubsection{\EHNC{37}}
\setcounter{mycnt}{37}
\[\Gamma(W,S) = \quad\csname CoxGrHNC\Roman{mycnt}\endcsname\]

\vskip12pt 
\[ M=\begin{pmatrix} 
1 & 4 & 2 & 2 & 2 & 2 \\
4 & 1 & 3 & 2 & 2 & 2 \\
2 & 3 & 1 & 3 & 3 & 2 \\
2 & 2 & 3 & 1 & 2 & 2 \\
2 & 2 & 3 & 2 & 1 & 3 \\
2 & 2 & 2 & 2 & 3 & 1
\end{pmatrix}\]
Numerator of $p_{(W,S)}(t)$:
\begin{dmath*} n_{(W,S)}(t) = 
1+4t+9t^{2}+17t^{3}+28t^{4}+41t^{5}+56t^{6}+72t^{7}+87t^{8}+100t^{9}+110t^{10}+115t^{11}+115t^{12}+110t^{13}+100t^{14}+87t^{15}+72t^{16}+56t^{17}+41t^{18}+28t^{19}+17t^{20}+9t^{21}+4t^{22}+t^{23}
\end{dmath*}
Denominator of $p_{(W,S)}(t)$:
\begin{dmath*} d_{(W,S)}(t) = 
1-2t+t^{2}-t^{4}+t^{5}-t^{7}+2t^{8}-t^{9}+t^{11}-t^{13}+t^{14}-t^{15}+t^{17}-2t^{18}+t^{19}-t^{21}+t^{22}
\end{dmath*}
Initial values of $a_k$:
\begin{dmath*}\breakingcomma
 (a_k)_{k=0}^{23}=(1
, 6
, 20
, 51
, 111
, 217
, 393
, 673
, 1104
, 1750
, 2699
, 4071
, 6029
, 8794
, 12666
, 18052
, 25503
, 35764
, 49841
, 69092
, 95348
, 131078
, 179610
, 245427
).
 \end{dmath*}
Exponential growth rate:
\[\omega(W,S) =
1.34370894458074694707436312719
\dots\]
Cocompact? No\\
$(W,S)\in \CM$? Yes\\
Magma command: \texttt{HyperbolicCoxeterGraph(51)} or \texttt{HyperbolicCoxeterMatrix(51)} \\
\newpage
\subsubsection{\EHNC{38}}
\setcounter{mycnt}{38}
\[\Gamma(W,S) = \quad\csname CoxGrHNC\Roman{mycnt}\endcsname\]

\vskip12pt 
\[ M=\begin{pmatrix} 
1 & 4 & 2 & 2 & 2 & 2 \\
4 & 1 & 3 & 2 & 2 & 2 \\
2 & 3 & 1 & 3 & 3 & 2 \\
2 & 2 & 3 & 1 & 2 & 2 \\
2 & 2 & 3 & 2 & 1 & 4 \\
2 & 2 & 2 & 2 & 4 & 1
\end{pmatrix}\]
Numerator of $p_{(W,S)}(t)$:
\begin{dmath*} n_{(W,S)}(t) = 
1+5t+13t^{2}+25t^{3}+40t^{4}+56t^{5}+71t^{6}+83t^{7}+90t^{8}+90t^{9}+83t^{10}+71t^{11}+56t^{12}+40t^{13}+25t^{14}+13t^{15}+5t^{16}+t^{17}
\end{dmath*}
Denominator of $p_{(W,S)}(t)$:
\begin{dmath*} d_{(W,S)}(t) = 
1-t-t^{2}-t^{3}+t^{4}+t^{6}+3t^{8}-t^{11}-t^{13}-2t^{14}-2t^{15}+t^{16}+2t^{17}
\end{dmath*}
Initial values of $a_k$:
\begin{dmath*}\breakingcomma
 (a_k)_{k=0}^{17}=(1
, 6
, 20
, 52
, 117
, 239
, 458
, 839
, 1486
, 2564
, 4337
, 7225
, 11893
, 19396
, 31402
, 50549
, 81006
, 129356
).
 \end{dmath*}
Exponential growth rate:
\[\omega(W,S) =
1.57569318209261794675800072515
\dots\]
Cocompact? No\\
$(W,S)\in \CM$? No\\
Magma command: \texttt{HyperbolicCoxeterGraph(52)} or \texttt{HyperbolicCoxeterMatrix(52)} \\
\newpage
\subsubsection{\EHNC{39}}
\setcounter{mycnt}{39}
\[\Gamma(W,S) = \quad\csname CoxGrHNC\Roman{mycnt}\endcsname\]

\vskip12pt 
\[ M=\begin{pmatrix} 
1 & 3 & 2 & 2 & 2 & 2 \\
3 & 1 & 3 & 3 & 3 & 2 \\
2 & 3 & 1 & 2 & 2 & 3 \\
2 & 3 & 2 & 1 & 2 & 2 \\
2 & 3 & 2 & 2 & 1 & 2 \\
2 & 2 & 3 & 2 & 2 & 1
\end{pmatrix}\]
Numerator of $p_{(W,S)}(t)$:
\begin{dmath*} n_{(W,S)}(t) = 
1+5t+13t^{2}+25t^{3}+40t^{4}+56t^{5}+71t^{6}+83t^{7}+90t^{8}+90t^{9}+83t^{10}+71t^{11}+56t^{12}+40t^{13}+25t^{14}+13t^{15}+5t^{16}+t^{17}
\end{dmath*}
Denominator of $p_{(W,S)}(t)$:
\begin{dmath*} d_{(W,S)}(t) = 
1-t-t^{2}-t^{3}+t^{4}+2t^{5}-t^{7}+t^{8}-t^{12}-t^{14}+t^{16}
\end{dmath*}
Initial values of $a_k$:
\begin{dmath*}\breakingcomma
 (a_k)_{k=0}^{17}=(1
, 6
, 20
, 52
, 117
, 237
, 445
, 791
, 1347
, 2216
, 3550
, 5568
, 8582
, 13044
, 19604
, 29189
, 43129
, 63332
).
 \end{dmath*}
Exponential growth rate:
\[\omega(W,S) =
1.42400358256305844195776866838
\dots\]
Cocompact? No\\
$(W,S)\in \CM$? Yes\\
Magma command: \texttt{HyperbolicCoxeterGraph(53)} or \texttt{HyperbolicCoxeterMatrix(53)} \\
\newpage
\subsubsection{\EHNC{40}}
\setcounter{mycnt}{40}
\[\Gamma(W,S) = \quad\csname CoxGrHNC\Roman{mycnt}\endcsname\]

\vskip12pt 
\[ M=\begin{pmatrix} 
1 & 3 & 2 & 2 & 2 & 2 \\
3 & 1 & 3 & 3 & 3 & 2 \\
2 & 3 & 1 & 2 & 2 & 4 \\
2 & 3 & 2 & 1 & 2 & 2 \\
2 & 3 & 2 & 2 & 1 & 2 \\
2 & 2 & 4 & 2 & 2 & 1
\end{pmatrix}\]
Numerator of $p_{(W,S)}(t)$:
\begin{dmath*} n_{(W,S)}(t) = 
1+5t+13t^{2}+25t^{3}+40t^{4}+56t^{5}+71t^{6}+83t^{7}+90t^{8}+90t^{9}+83t^{10}+71t^{11}+56t^{12}+40t^{13}+25t^{14}+13t^{15}+5t^{16}+t^{17}
\end{dmath*}
Denominator of $p_{(W,S)}(t)$:
\begin{dmath*} d_{(W,S)}(t) = 
1-t-t^{2}-2t^{3}+2t^{4}+2t^{6}-t^{7}+4t^{8}-2t^{11}+t^{12}-t^{13}-3t^{14}-3t^{15}+t^{16}+3t^{17}
\end{dmath*}
Initial values of $a_k$:
\begin{dmath*}\breakingcomma
 (a_k)_{k=0}^{17}=(1
, 6
, 20
, 53
, 123
, 260
, 518
, 990
, 1834
, 3320
, 5908
, 10380
, 18059
, 31188
, 53561
, 91588
, 156107
, 265419
).
 \end{dmath*}
Exponential growth rate:
\[\omega(W,S) =
1.68670898797603657664798272337
\dots\]
Cocompact? No\\
$(W,S)\in \CM$? No\\
Magma command: \texttt{HyperbolicCoxeterGraph(54)} or \texttt{HyperbolicCoxeterMatrix(54)} \\
\newpage
\subsubsection{\EHNC{41}}
\setcounter{mycnt}{41}
\[\Gamma(W,S) = \quad\csname CoxGrHNC\Roman{mycnt}\endcsname\]

\vskip12pt 
\[ M=\begin{pmatrix} 
1 & 3 & 3 & 3 & 3 & 3 \\
3 & 1 & 2 & 2 & 2 & 2 \\
3 & 2 & 1 & 2 & 2 & 2 \\
3 & 2 & 2 & 1 & 2 & 2 \\
3 & 2 & 2 & 2 & 1 & 2 \\
3 & 2 & 2 & 2 & 2 & 1
\end{pmatrix}\]
Numerator of $p_{(W,S)}(t)$:
\begin{dmath*} n_{(W,S)}(t) = 
1+4t+8t^{2}+13t^{3}+18t^{4}+20t^{5}+20t^{6}+18t^{7}+13t^{8}+8t^{9}+4t^{10}+t^{11}
\end{dmath*}
Denominator of $p_{(W,S)}(t)$:
\begin{dmath*} d_{(W,S)}(t) = 
1-2t-2t^{3}+5t^{4}-t^{5}+t^{6}+2t^{8}-5t^{9}-3t^{10}+4t^{11}
\end{dmath*}
Initial values of $a_k$:
\begin{dmath*}\breakingcomma
 (a_k)_{k=0}^{11}=(1
, 6
, 20
, 55
, 135
, 301
, 637
, 1301
, 2575
, 5000
, 9580
, 18146
).
 \end{dmath*}
Exponential growth rate:
\[\omega(W,S) =
1.84299223582523685359162961980
\dots\]
Cocompact? No\\
$(W,S)\in \CM$? Yes\\
Magma command: \texttt{HyperbolicCoxeterGraph(55)} or \texttt{HyperbolicCoxeterMatrix(55)} \\
\newpage
\subsubsection{\EHNC{42}}
\setcounter{mycnt}{42}
\[\Gamma(W,S) = \quad\csname CoxGrHNC\Roman{mycnt}\endcsname\]

\vskip12pt 
\[ M=\begin{pmatrix} 
1 & 4 & 2 & 3 & 2 & 2 \\
4 & 1 & 3 & 2 & 2 & 2 \\
2 & 3 & 1 & 2 & 2 & 3 \\
3 & 2 & 2 & 1 & 3 & 2 \\
2 & 2 & 2 & 3 & 1 & 3 \\
2 & 2 & 3 & 2 & 3 & 1
\end{pmatrix}\]
Numerator of $p_{(W,S)}(t)$:
\begin{dmath*} n_{(W,S)}(t) = 
1+4t+9t^{2}+16t^{3}+25t^{4}+37t^{5}+52t^{6}+69t^{7}+87t^{8}+105t^{9}+123t^{10}+140t^{11}+154t^{12}+163t^{13}+167t^{14}+167t^{15}+163t^{16}+154t^{17}+140t^{18}+123t^{19}+105t^{20}+87t^{21}+69t^{22}+52t^{23}+37t^{24}+25t^{25}+16t^{26}+9t^{27}+4t^{28}+t^{29}
\end{dmath*}
Denominator of $p_{(W,S)}(t)$:
\begin{dmath*} d_{(W,S)}(t) = 
1-2t+t^{3}-t^{6}+2t^{8}-t^{11}+2t^{12}+t^{14}-t^{15}+t^{18}-2t^{19}-2t^{22}+t^{23}-t^{26}+t^{29}
\end{dmath*}
Initial values of $a_k$:
\begin{dmath*}\breakingcomma
 (a_k)_{k=0}^{29}=(1
, 6
, 21
, 57
, 133
, 282
, 560
, 1062
, 1948
, 3486
, 6124
, 10609
, 18184
, 30914
, 52228
, 87815
, 147108
, 245739
, 409604
, 681590
, 1132709
, 1880523
, 3119635
, 5172145
, 8571140
, 14198848
, 23515238
, 38936283
, 64459855
, 106701370
).
 \end{dmath*}
Exponential growth rate:
\[\omega(W,S) =
1.65462148869357522321191538305
\dots\]
Cocompact? No\\
$(W,S)\in \CM$? No\\
Magma command: \texttt{HyperbolicCoxeterGraph(56)} or \texttt{HyperbolicCoxeterMatrix(56)} \\
\newpage
\subsubsection{\EHNC{43}}
\setcounter{mycnt}{43}
\[\Gamma(W,S) = \quad\csname CoxGrHNC\Roman{mycnt}\endcsname\]

\vskip12pt 
\[ M=\begin{pmatrix} 
1 & 4 & 2 & 2 & 2 & 3 \\
4 & 1 & 3 & 2 & 2 & 2 \\
2 & 3 & 1 & 3 & 2 & 2 \\
2 & 2 & 3 & 1 & 4 & 2 \\
2 & 2 & 2 & 4 & 1 & 3 \\
3 & 2 & 2 & 2 & 3 & 1
\end{pmatrix}\]
Numerator of $p_{(W,S)}(t)$:
\begin{dmath*} n_{(W,S)}(t) = 
1+4t+9t^{2}+16t^{3}+25t^{4}+36t^{5}+48t^{6}+60t^{7}+71t^{8}+80t^{9}+87t^{10}+92t^{11}+94t^{12}+92t^{13}+87t^{14}+80t^{15}+71t^{16}+60t^{17}+48t^{18}+36t^{19}+25t^{20}+16t^{21}+9t^{22}+4t^{23}+t^{24}
\end{dmath*}
Denominator of $p_{(W,S)}(t)$:
\begin{dmath*} d_{(W,S)}(t) = 
1-2t+t^{4}-2t^{5}+3t^{6}-2t^{7}+5t^{8}-2t^{9}+2t^{10}-4t^{11}+8t^{12}-6t^{13}+2t^{14}-4t^{15}+7t^{16}-6t^{17}+3t^{18}-6t^{19}+3t^{20}-2t^{21}-4t^{23}+5t^{24}
\end{dmath*}
Initial values of $a_k$:
\begin{dmath*}\breakingcomma
 (a_k)_{k=0}^{24}=(1
, 6
, 21
, 58
, 140
, 312
, 660
, 1348
, 2687
, 5262
, 10176
, 19504
, 37145
, 70416
, 133044
, 250774
, 471870
, 886798
, 1665097
, 3124462
, 5860164
, 10987520
, 20596135
, 38600866
, 72336018
).
 \end{dmath*}
Exponential growth rate:
\[\omega(W,S) =
1.87335449595753946905901409963
\dots\]
Cocompact? No\\
$(W,S)\in \CM$? No\\
Magma command: \texttt{HyperbolicCoxeterGraph(57)} or \texttt{HyperbolicCoxeterMatrix(57)} \\
\newpage
\subsubsection{\EHNC{44}}
\setcounter{mycnt}{44}
\[\Gamma(W,S) = \quad\csname CoxGrHNC\Roman{mycnt}\endcsname\]

\vskip12pt 
\[ M=\begin{pmatrix} 
1 & 3 & 2 & 2 & 2 & 2 \\
3 & 1 & 3 & 2 & 2 & 3 \\
2 & 3 & 1 & 3 & 2 & 2 \\
2 & 2 & 3 & 1 & 3 & 2 \\
2 & 2 & 2 & 3 & 1 & 3 \\
2 & 3 & 2 & 2 & 3 & 1
\end{pmatrix}\]
Numerator of $p_{(W,S)}(t)$:
\begin{dmath*} n_{(W,S)}(t) = 
1+4t+9t^{2}+16t^{3}+25t^{4}+35t^{5}+44t^{6}+51t^{7}+55t^{8}+55t^{9}+51t^{10}+44t^{11}+35t^{12}+25t^{13}+16t^{14}+9t^{15}+4t^{16}+t^{17}
\end{dmath*}
Denominator of $p_{(W,S)}(t)$:
\begin{dmath*} d_{(W,S)}(t) = 
1-2t+t^{3}+t^{7}-t^{8}+t^{9}-t^{11}-t^{15}+t^{16}
\end{dmath*}
Initial values of $a_k$:
\begin{dmath*}\breakingcomma
 (a_k)_{k=0}^{17}=(1
, 6
, 21
, 57
, 133
, 280
, 547
, 1011
, 1792
, 3076
, 5150
, 8456
, 13673
, 21842
, 34557
, 54256
, 84664
, 131468
).
 \end{dmath*}
Exponential growth rate:
\[\omega(W,S) =
1.52179436592917726350419855685
\dots\]
Cocompact? No\\
$(W,S)\in \CM$? Yes\\
Magma command: \texttt{HyperbolicCoxeterGraph(58)} or \texttt{HyperbolicCoxeterMatrix(58)} \\
\newpage
\subsubsection{\EHNC{45}}
\setcounter{mycnt}{45}
\[\Gamma(W,S) = \quad\csname CoxGrHNC\Roman{mycnt}\endcsname\]

\vskip12pt 
\[ M=\begin{pmatrix} 
1 & 3 & 2 & 2 & 2 & 2 & 2 \\
3 & 1 & 3 & 3 & 2 & 2 & 2 \\
2 & 3 & 1 & 2 & 2 & 2 & 3 \\
2 & 3 & 2 & 1 & 3 & 2 & 2 \\
2 & 2 & 2 & 3 & 1 & 4 & 2 \\
2 & 2 & 2 & 2 & 4 & 1 & 2 \\
2 & 2 & 3 & 2 & 2 & 2 & 1
\end{pmatrix}\]
Numerator of $p_{(W,S)}(t)$:
\begin{dmath*} n_{(W,S)}(t) = 
-1-7t-27t^{2}-77t^{3}-181t^{4}-371t^{5}-686t^{6}-1170t^{7}-1869t^{8}-2826t^{9}-4075t^{10}-5635t^{11}-7504t^{12}-9654t^{13}-12029t^{14}-14546t^{15}-17099t^{16}-19566t^{17}-21818t^{18}-23729t^{19}-25186t^{20}-26099t^{21}-26410t^{22}-26099t^{23}-25186t^{24}-23729t^{25}-21818t^{26}-19566t^{27}-17099t^{28}-14546t^{29}-12029t^{30}-9654t^{31}-7504t^{32}-5635t^{33}-4075t^{34}-2826t^{35}-1869t^{36}-1170t^{37}-686t^{38}-371t^{39}-181t^{40}-77t^{41}-27t^{42}-7t^{43}-t^{44}
\end{dmath*}
Denominator of $p_{(W,S)}(t)$:
\begin{dmath*} d_{(W,S)}(t) = 
-1+t^{3}+t^{4}+t^{5}+t^{6}+t^{7}+t^{9}-t^{10}-t^{11}-3t^{12}-3t^{13}-4t^{14}-3t^{15}-3t^{16}-t^{17}-t^{18}+t^{19}+2t^{20}+4t^{21}+5t^{22}+5t^{23}+4t^{24}+4t^{25}+3t^{26}+2t^{27}+t^{28}-t^{29}-2t^{30}-3t^{31}-3t^{32}-3t^{33}-3t^{34}-2t^{35}-t^{36}-t^{37}+t^{40}+t^{41}+t^{42}+t^{43}
\end{dmath*}
Initial values of $a_k$:
\begin{dmath*}\breakingcomma
 (a_k)_{k=0}^{44}=(1
, 7
, 27
, 78
, 189
, 406
, 799
, 1472
, 2576
, 4326
, 7025
, 11096
, 17124
, 25912
, 38557
, 56551
, 81916
, 117385
, 166642
, 234641
, 328031
, 455722
, 629637
, 865712
, 1185222
, 1616541
, 2197479
, 2978382
, 4026244
, 5430161
, 7308561
, 9818787
, 13169802
, 17639024
, 23594631
, 31525115
, 42078438
, 56113898
, 74770831
, 99559608
, 132482147
, 176191512
, 234203270
, 311175373
, 413278763
).
 \end{dmath*}
Exponential growth rate:
\[\omega(W,S) =
1.32374816773374595778541093503
\dots\]
Cocompact? No\\
$(W,S)\in \CM$? Yes\\
Magma command: \texttt{HyperbolicCoxeterGraph(59)} or \texttt{HyperbolicCoxeterMatrix(59)} \\
\newpage
\subsubsection{\EHNC{46}}
\setcounter{mycnt}{46}
\[\Gamma(W,S) = \quad\csname CoxGrHNC\Roman{mycnt}\endcsname\]

\vskip12pt 
\[ M=\begin{pmatrix} 
1 & 3 & 2 & 2 & 2 & 2 & 2 \\
3 & 1 & 3 & 3 & 2 & 2 & 2 \\
2 & 3 & 1 & 2 & 2 & 2 & 2 \\
2 & 3 & 2 & 1 & 3 & 3 & 2 \\
2 & 2 & 2 & 3 & 1 & 2 & 2 \\
2 & 2 & 2 & 3 & 2 & 1 & 3 \\
2 & 2 & 2 & 2 & 2 & 3 & 1
\end{pmatrix}\]
Numerator of $p_{(W,S)}(t)$:
\begin{dmath*} n_{(W,S)}(t) = 
-1-7t-26t^{2}-70t^{3}-155t^{4}-301t^{5}-531t^{6}-869t^{7}-1338t^{8}-1957t^{9}-2737t^{10}-3678t^{11}-4767t^{12}-5976t^{13}-7262t^{14}-8570t^{15}-9837t^{16}-10996t^{17}-11981t^{18}-12733t^{19}-13205t^{20}-13366t^{21}-13205t^{22}-12733t^{23}-11981t^{24}-10996t^{25}-9837t^{26}-8570t^{27}-7262t^{28}-5976t^{29}-4767t^{30}-3678t^{31}-2737t^{32}-1957t^{33}-1338t^{34}-869t^{35}-531t^{36}-301t^{37}-155t^{38}-70t^{39}-26t^{40}-7t^{41}-t^{42}
\end{dmath*}
Denominator of $p_{(W,S)}(t)$:
\begin{dmath*} d_{(W,S)}(t) = 
-1+t^{2}+t^{3}+t^{4}-t^{10}-2t^{11}-2t^{12}-2t^{13}-2t^{14}-t^{15}+t^{19}+3t^{20}+4t^{21}+3t^{22}+2t^{23}+t^{24}+2t^{25}+2t^{26}-t^{28}-t^{29}-2t^{30}-2t^{31}-2t^{32}-t^{33}-t^{34}-t^{35}+t^{40}+t^{41}
\end{dmath*}
Initial values of $a_k$:
\begin{dmath*}\breakingcomma
 (a_k)_{k=0}^{42}=(1
, 7
, 27
, 78
, 190
, 413
, 826
, 1550
, 2767
, 4746
, 7879
, 12732
, 20116
, 31185
, 47573
, 71585
, 106462
, 156748
, 228798
, 331481
, 477153
, 683003
, 972913
, 1380027
, 1950298
, 2747385
, 3859412
, 5408292
, 7562588
, 10555253
, 14708093
, 20465493
, 28440911
, 39480969
, 54753788
, 75870727
, 105054152
, 145368622
, 201039439
, 277891557
, 383954305
, 530294533
, 732164421
).
 \end{dmath*}
Exponential growth rate:
\[\omega(W,S) =
1.37735023351086082704676893243
\dots\]
Cocompact? No\\
$(W,S)\in \CM$? Yes\\
Magma command: \texttt{HyperbolicCoxeterGraph(60)} or \texttt{HyperbolicCoxeterMatrix(60)} \\
\newpage
% produced with RUNfortables.magma
\subsubsection{\EHNC{47}}
\setcounter{mycnt}{47}
\[\Gamma(W,S) = \quad\csname CoxGrHNC\Roman{mycnt}\endcsname\]

\vskip12pt 
\[ M=\begin{pmatrix} 
1 & 3 & 2 & 2 & 2 & 3 & 2 \\
3 & 1 & 3 & 2 & 2 & 2 & 2 \\
2 & 3 & 1 & 3 & 2 & 2 & 2 \\
2 & 2 & 3 & 1 & 3 & 2 & 2 \\
2 & 2 & 2 & 3 & 1 & 3 & 2 \\
3 & 2 & 2 & 2 & 3 & 1 & 3 \\
2 & 2 & 2 & 2 & 2 & 3 & 1
\end{pmatrix}\]
Numerator of $p_{(W,S)}(t)$:
\begin{dmath*} n_{(W,S)}(t) = 
-1-8t-34t^{2}-104t^{3}-259t^{4}-560t^{5}-1091t^{6}-1959t^{7}-3290t^{8}-5221t^{9}-7888t^{10}-11411t^{11}-15877t^{12}-21322t^{13}-27715t^{14}-34947t^{15}-42827t^{16}-51086t^{17}-59389t^{18}-67355t^{19}-74584t^{20}-80688t^{21}-85323t^{22}-88219t^{23}-89204t^{24}-88219t^{25}-85323t^{26}-80688t^{27}-74584t^{28}-67355t^{29}-59389t^{30}-51086t^{31}-42827t^{32}-34947t^{33}-27715t^{34}-21322t^{35}-15877t^{36}-11411t^{37}-7888t^{38}-5221t^{39}-3290t^{40}-1959t^{41}-1091t^{42}-560t^{43}-259t^{44}-104t^{45}-34t^{46}-8t^{47}-t^{48}
\end{dmath*}
Denominator of $p_{(W,S)}(t)$:
\begin{dmath*} d_{(W,S)}(t) = 
-1-t+t^{2}+2t^{3}+2t^{4}+2t^{5}+2t^{6}+2t^{7}+t^{8}-2t^{10}-6t^{11}-9t^{12}-10t^{13}-12t^{14}-12t^{15}-10t^{16}-8t^{17}-5t^{18}+6t^{20}+13t^{21}+17t^{22}+19t^{23}+21t^{24}+22t^{25}+21t^{26}+18t^{27}+12t^{28}+6t^{29}+t^{30}-3t^{31}-7t^{32}-11t^{33}-13t^{34}-12t^{35}-11t^{36}-9t^{37}-6t^{38}-4t^{39}-2t^{40}+t^{41}+2t^{42}+2t^{43}+2t^{44}+2t^{45}+2t^{46}+t^{47}
\end{dmath*}
Initial values of $a_k$:
\begin{dmath*}\breakingcomma
 (a_k)_{k=0}^{48}=(1
, 7
, 28
, 85
, 218
, 499
, 1052
, 2084
, 3933
, 7143
, 12580
, 21611
, 36379
, 60226
, 98342
, 158760
, 253880
, 402800
, 634880
, 995187
, 1552811
, 2413559
, 3739327
, 5777651
, 8906778
, 13704392
, 21052395
, 32296630
, 49490329
, 75765133
, 115896490
, 177165202
, 270670159
, 413328443
, 630922591
, 962743095
, 1468661022
, 2239902537
, 3415462591
, 5207108762
, 7937470340
, 12098059829
, 18437656740
, 28096940887
, 42813575064
, 65234599093
, 99392284324
, 151428976166
, 230701201950
).
 \end{dmath*}
Exponential growth rate:
\[\omega(W,S) =
1.52320775283771810350512463006
\dots\]
Cocompact? No\\
$(W,S)\in \CM$? Yes\\
Magma command: \texttt{HyperbolicCoxeterGraph(61)} or \texttt{HyperbolicCoxeterMatrix(61)} \\
\newpage
\subsubsection{\EHNC{48}}
\setcounter{mycnt}{48}
\[\Gamma(W,S) = \quad\csname CoxGrHNC\Roman{mycnt}\endcsname\]

\vskip12pt 
\[ M=\begin{pmatrix} 
1 & 3 & 2 & 2 & 2 & 2 & 2 & 2 \\
3 & 1 & 3 & 2 & 2 & 2 & 2 & 2 \\
2 & 3 & 1 & 3 & 3 & 2 & 2 & 2 \\
2 & 2 & 3 & 1 & 2 & 2 & 2 & 2 \\
2 & 2 & 3 & 2 & 1 & 3 & 2 & 2 \\
2 & 2 & 2 & 2 & 3 & 1 & 3 & 2 \\
2 & 2 & 2 & 2 & 2 & 3 & 1 & 4 \\
2 & 2 & 2 & 2 & 2 & 2 & 4 & 1
\end{pmatrix}\]
Numerator of $p_{(W,S)}(t)$:
\begin{dmath*} n_{(W,S)}(t) = 
1+5t+15t^{2}+35t^{3}+69t^{4}+122t^{5}+199t^{6}+306t^{7}+449t^{8}+633t^{9}+863t^{10}+1142t^{11}+1472t^{12}+1852t^{13}+2279t^{14}+2748t^{15}+3252t^{16}+3783t^{17}+4330t^{18}+4880t^{19}+5418t^{20}+5929t^{21}+6399t^{22}+6815t^{23}+7165t^{24}+7438t^{25}+7625t^{26}+7720t^{27}+7720t^{28}+7625t^{29}+7438t^{30}+7165t^{31}+6815t^{32}+6399t^{33}+5929t^{34}+5418t^{35}+4880t^{36}+4330t^{37}+3783t^{38}+3252t^{39}+2748t^{40}+2279t^{41}+1852t^{42}+1472t^{43}+1142t^{44}+863t^{45}+633t^{46}+449t^{47}+306t^{48}+199t^{49}+122t^{50}+69t^{51}+35t^{52}+15t^{53}+5t^{54}+t^{55}
\end{dmath*}
Denominator of $p_{(W,S)}(t)$:
\begin{dmath*} d_{(W,S)}(t) = 
1-3t+4t^{2}-5t^{3}+6t^{4}-5t^{5}+3t^{6}-t^{7}-t^{8}+2t^{9}-t^{10}-t^{11}+4t^{12}-7t^{13}+10t^{14}-12t^{15}+12t^{16}-11t^{17}+8t^{18}-4t^{19}+2t^{20}-t^{21}-t^{22}+3t^{24}-5t^{25}+6t^{26}-7t^{27}+6t^{28}-3t^{29}+3t^{31}-6t^{32}+7t^{33}-6t^{34}+5t^{35}-2t^{36}-t^{37}+2t^{38}-2t^{39}+2t^{40}-3t^{42}+5t^{43}-6t^{44}+7t^{45}-7t^{46}+5t^{47}-3t^{48}+t^{49}+t^{50}-2t^{51}+2t^{52}-2t^{53}+t^{54}
\end{dmath*}
Initial values of $a_k$:
\begin{dmath*}\breakingcomma
 (a_k)_{k=0}^{55}=(1
, 8
, 35
, 113
, 302
, 708
, 1507
, 2979
, 5555
, 9881
, 16906
, 28002
, 45126
, 71038
, 109594
, 166138
, 248026
, 365325
, 531743
, 765865
, 1092793
, 1546319
, 2171802
, 3029974
, 4201967
, 5795948
, 7955871
, 10873012
, 14801163
, 20076637
, 27144591
, 36593645
, 49201396
, 65994232
, 88325911
, 117980762
, 157309185
, 209405510
, 278341407
, 369472130
, 489838244
, 648692517
, 858190864
, 1134298286
, 1497976553
, 1976741069
, 2606701445
, 3435235806
, 4524495335
, 5955996418
, 7836637467
, 10306581890
, 13549585350
, 17806524457
, 23393118426
, 30723142138
).
 \end{dmath*}
Exponential growth rate:
\[\omega(W,S) =
1.30923602686968771526540931654
\dots\]
Cocompact? No\\
$(W,S)\in \CM$? Yes\\
Magma command: \texttt{HyperbolicCoxeterGraph(62)} or \texttt{HyperbolicCoxeterMatrix(62)} \\
\newpage
\subsubsection{\EHNC{49}}
\setcounter{mycnt}{49}
\[\Gamma(W,S) = \quad\csname CoxGrHNC\Roman{mycnt}\endcsname\]

\vskip12pt 
\[ M=\begin{pmatrix} 
1 & 3 & 2 & 2 & 2 & 2 & 2 & 2 \\
3 & 1 & 3 & 3 & 2 & 2 & 2 & 2 \\
2 & 3 & 1 & 2 & 2 & 2 & 2 & 2 \\
2 & 3 & 2 & 1 & 3 & 2 & 2 & 2 \\
2 & 2 & 2 & 3 & 1 & 3 & 3 & 2 \\
2 & 2 & 2 & 2 & 3 & 1 & 2 & 2 \\
2 & 2 & 2 & 2 & 3 & 2 & 1 & 3 \\
2 & 2 & 2 & 2 & 2 & 2 & 3 & 1
\end{pmatrix}\]
Numerator of $p_{(W,S)}(t)$:
\begin{dmath*} n_{(W,S)}(t) = 
1+5t+15t^{2}+35t^{3}+69t^{4}+121t^{5}+194t^{6}+291t^{7}+414t^{8}+565t^{9}+746t^{10}+958t^{11}+1201t^{12}+1472t^{13}+1767t^{14}+2079t^{15}+2401t^{16}+2725t^{17}+3043t^{18}+3347t^{19}+3628t^{20}+3878t^{21}+4088t^{22}+4252t^{23}+4364t^{24}+4421t^{25}+4421t^{26}+4364t^{27}+4252t^{28}+4088t^{29}+3878t^{30}+3628t^{31}+3347t^{32}+3043t^{33}+2725t^{34}+2401t^{35}+2079t^{36}+1767t^{37}+1472t^{38}+1201t^{39}+958t^{40}+746t^{41}+565t^{42}+414t^{43}+291t^{44}+194t^{45}+121t^{46}+69t^{47}+35t^{48}+15t^{49}+5t^{50}+t^{51}
\end{dmath*}
Denominator of $p_{(W,S)}(t)$:
\begin{dmath*} d_{(W,S)}(t) = 
1-3t+4t^{2}-5t^{3}+6t^{4}-7t^{5}+9t^{6}-9t^{7}+9t^{8}-10t^{9}+10t^{10}-11t^{11}+12t^{12}-12t^{13}+14t^{14}-16t^{15}+17t^{16}-18t^{17}+17t^{18}-16t^{19}+15t^{20}-13t^{21}+10t^{22}-6t^{23}+2t^{24}-2t^{26}+4t^{27}-5t^{28}+6t^{29}-6t^{30}+6t^{31}-7t^{32}+8t^{33}-9t^{34}+10t^{35}-10t^{36}+11t^{37}-11t^{38}+10t^{39}-9t^{40}+7t^{41}-5t^{42}+3t^{43}-t^{44}+t^{46}-2t^{47}+2t^{48}-2t^{49}+t^{50}
\end{dmath*}
Initial values of $a_k$:
\begin{dmath*}\breakingcomma
 (a_k)_{k=0}^{51}=(1
, 8
, 35
, 113
, 302
, 709
, 1515
, 3014
, 5668
, 10184
, 17623
, 29553
, 48262
, 77052
, 120644
, 185736
, 281771
, 421991
, 624878
, 916122
, 1331305
, 1919560
, 2748553
, 3911257
, 5535151
, 7794704
, 10928310
, 15261246
, 21236776
, 29458262
, 40746153
, 56215081
, 77378126
, 106287771
, 145726392
, 199463621
, 272603980
, 372056358
, 507167904
, 690579753
, 939382020
, 1276672513
, 1733660019
, 2352502092
, 3190133420
, 4323430038
, 5856174918
, 7928452581
, 10729318900
, 14513886794
, 19626365574
, 26531126928
).
 \end{dmath*}
Exponential growth rate:
\[\omega(W,S) =
1.34774184227273844060669340702
\dots\]
Cocompact? No\\
$(W,S)\in \CM$? Yes\\
Magma command: \texttt{HyperbolicCoxeterGraph(63)} or \texttt{HyperbolicCoxeterMatrix(63)} \\
\newpage
\subsubsection{\EHNC{50}}
\setcounter{mycnt}{50}
\[\Gamma(W,S) = \quad\csname CoxGrHNC\Roman{mycnt}\endcsname\]

\vskip12pt 
\[ M=\begin{pmatrix} 
1 & 3 & 2 & 2 & 2 & 2 & 3 & 2 \\
3 & 1 & 3 & 3 & 2 & 2 & 2 & 2 \\
2 & 3 & 1 & 2 & 2 & 3 & 2 & 2 \\
2 & 3 & 2 & 1 & 3 & 2 & 2 & 2 \\
2 & 2 & 2 & 3 & 1 & 2 & 2 & 2 \\
2 & 2 & 3 & 2 & 2 & 1 & 2 & 2 \\
3 & 2 & 2 & 2 & 2 & 2 & 1 & 3 \\
2 & 2 & 2 & 2 & 2 & 2 & 3 & 1
\end{pmatrix}\]
Numerator of $p_{(W,S)}(t)$:
\begin{dmath*} n_{(W,S)}(t) = 
1+7t+27t^{2}+77t^{3}+182t^{4}+378t^{5}+713t^{6}+1247t^{7}+2051t^{8}+3205t^{9}+4795t^{10}+6909t^{11}+9632t^{12}+13040t^{13}+17194t^{14}+22134t^{15}+27874t^{16}+34398t^{17}+41657t^{18}+49567t^{19}+58009t^{20}+66831t^{21}+75852t^{22}+84868t^{23}+93659t^{24}+101997t^{25}+109655t^{26}+116417t^{27}+122087t^{28}+126497t^{29}+129514t^{30}+131046t^{31}+131046t^{32}+129514t^{33}+126497t^{34}+122087t^{35}+116417t^{36}+109655t^{37}+101997t^{38}+93659t^{39}+84868t^{40}+75852t^{41}+66831t^{42}+58009t^{43}+49567t^{44}+41657t^{45}+34398t^{46}+27874t^{47}+22134t^{48}+17194t^{49}+13040t^{50}+9632t^{51}+6909t^{52}+4795t^{53}+3205t^{54}+2051t^{55}+1247t^{56}+713t^{57}+378t^{58}+182t^{59}+77t^{60}+27t^{61}+7t^{62}+t^{63}
\end{dmath*}
Denominator of $p_{(W,S)}(t)$:
\begin{dmath*} d_{(W,S)}(t) = 
1-t-t^{4}-t^{8}+2t^{9}-t^{10}+t^{11}+t^{12}+t^{14}+t^{15}-t^{17}+3t^{18}-2t^{19}-t^{22}-2t^{23}+t^{24}-2t^{25}-2t^{26}+t^{27}-2t^{28}-t^{29}-2t^{32}+2t^{33}-t^{34}+2t^{36}+t^{37}+2t^{39}+t^{40}+2t^{42}+t^{44}+t^{46}-t^{47}+t^{48}-t^{49}-t^{50}-t^{52}-t^{54}-t^{56}+t^{62}
\end{dmath*}
Initial values of $a_k$:
\begin{dmath*}\breakingcomma
 (a_k)_{k=0}^{63}=(1
, 8
, 35
, 112
, 295
, 681
, 1429
, 2788
, 5135
, 9027
, 15271
, 25017
, 39881
, 62108
, 94787
, 142133
, 209857
, 305652
, 439829
, 626145
, 882880
, 1234240
, 1712183
, 2358789
, 3229332
, 4396262
, 5954364
, 8027434
, 10776910
, 14413026
, 19209222
, 25520749
, 33808676
, 44670858
, 58881873
, 77444504
, 101656076
, 133193914
, 174225414
, 227549778
, 296780471
, 386580055
, 502962393
, 653681486
, 848731688
, 1100991108
, 1427049094
, 1848270346
, 2392163162
, 3094138558
, 3999771737
, 5167709130
, 6673404991
, 8613923915
, 11114112975
, 14334533628
, 18481654550
, 23820949180
, 30693724996
, 39538746888
, 50920019216
, 65562479340
, 84397854087
, 108623571115
).
 \end{dmath*}
Exponential growth rate:
\[\omega(W,S) =
1.28425783453546273020101777611
\dots\]
Cocompact? No\\
$(W,S)\in \CM$? Yes\\
Magma command: \texttt{HyperbolicCoxeterGraph(64)} or \texttt{HyperbolicCoxeterMatrix(64)} \\
\newpage
\subsubsection{\EHNC{51}}
\setcounter{mycnt}{51}
\[\Gamma(W,S) = \quad\csname CoxGrHNC\Roman{mycnt}\endcsname\]

\vskip12pt 
\[ M=\begin{pmatrix} 
1 & 3 & 2 & 2 & 2 & 2 & 2 & 2 \\
3 & 1 & 3 & 3 & 2 & 2 & 2 & 2 \\
2 & 3 & 1 & 2 & 2 & 3 & 2 & 2 \\
2 & 3 & 2 & 1 & 3 & 2 & 2 & 2 \\
2 & 2 & 2 & 3 & 1 & 2 & 2 & 3 \\
2 & 2 & 3 & 2 & 2 & 1 & 3 & 2 \\
2 & 2 & 2 & 2 & 2 & 3 & 1 & 3 \\
2 & 2 & 2 & 2 & 3 & 2 & 3 & 1
\end{pmatrix}\]
Numerator of $p_{(W,S)}(t)$:
\begin{dmath*} n_{(W,S)}(t) = 
1+6t+21t^{2}+56t^{3}+125t^{4}+246t^{5}+440t^{6}+730t^{7}+1139t^{8}+1689t^{9}+2400t^{10}+3289t^{11}+4369t^{12}+5647t^{13}+7123t^{14}+8788t^{15}+10624t^{16}+12603t^{17}+14688t^{18}+16834t^{19}+18990t^{20}+21101t^{21}+23110t^{22}+24961t^{23}+26600t^{24}+27978t^{25}+29052t^{26}+29788t^{27}+30162t^{28}+30162t^{29}+29788t^{30}+29052t^{31}+27978t^{32}+26600t^{33}+24961t^{34}+23110t^{35}+21101t^{36}+18990t^{37}+16834t^{38}+14688t^{39}+12603t^{40}+10624t^{41}+8788t^{42}+7123t^{43}+5647t^{44}+4369t^{45}+3289t^{46}+2400t^{47}+1689t^{48}+1139t^{49}+730t^{50}+440t^{51}+246t^{52}+125t^{53}+56t^{54}+21t^{55}+6t^{56}+t^{57}
\end{dmath*}
Denominator of $p_{(W,S)}(t)$:
\begin{dmath*} d_{(W,S)}(t) = 
1-2t+t^{2}-t^{3}+t^{5}-t^{6}+3t^{7}-2t^{8}+3t^{9}-3t^{10}+2t^{11}-2t^{12}+t^{13}+2t^{14}-4t^{15}+6t^{16}-7t^{17}+5t^{18}-6t^{19}+3t^{20}-3t^{21}+2t^{23}-5t^{24}+5t^{25}-7t^{26}+6t^{27}-6t^{28}+5t^{29}-3t^{30}+2t^{31}+t^{32}-2t^{33}+4t^{34}-4t^{35}+6t^{36}-4t^{37}+5t^{38}-3t^{39}+3t^{40}-t^{41}+2t^{43}-3t^{44}+3t^{45}-3t^{46}+2t^{47}-t^{48}-t^{50}-t^{53}+t^{54}-t^{55}+t^{56}
\end{dmath*}
Initial values of $a_k$:
\begin{dmath*}\breakingcomma
 (a_k)_{k=0}^{57}=(1
, 8
, 36
, 121
, 339
, 838
, 1891
, 3982
, 7943
, 15171
, 27970
, 50087
, 87548
, 149957
, 252507
, 419086
, 687064
, 1114656
, 1792229
, 2859642
, 4532806
, 7144326
, 11205639
, 17501946
, 27237152
, 42255037
, 65376593
, 100914333
, 155456142
, 239059571
, 367071019
, 562896143
, 862218066
, 1319418919
, 2017354196
, 3082228636
, 4706233833
, 7181994176
, 10954976418
, 16703225397
, 25458666374
, 38791633275
, 59091565192
, 89993973354
, 137029878038
, 208613604045
, 317545182799
, 483295403300
, 735481153159
, 1119150989516
, 1702823736286
, 2590713891827
, 3941324303089
, 5995722063918
, 9120540404928
, 13873375375548
, 21102239631352
, 32096808721498
).
 \end{dmath*}
Exponential growth rate:
\[\omega(W,S) =
1.52071977182142647854923187854
\dots\]
Cocompact? No\\
$(W,S)\in \CM$? Yes\\
Magma command: \texttt{HyperbolicCoxeterGraph(65)} or \texttt{HyperbolicCoxeterMatrix(65)} \\
\newpage
% produced with RUNfortables.magma
\subsubsection{\EHNC{52}}
\setcounter{mycnt}{52}
\[\Gamma(W,S) = \quad\csname CoxGrHNC\Roman{mycnt}\endcsname\]

\vskip12pt 
\[ M=\begin{pmatrix} 
1 & 3 & 2 & 2 & 2 & 2 & 2 & 2 & 2 \\
3 & 1 & 3 & 3 & 2 & 2 & 2 & 2 & 2 \\
2 & 3 & 1 & 2 & 2 & 3 & 2 & 2 & 2 \\
2 & 3 & 2 & 1 & 3 & 2 & 2 & 2 & 2 \\
2 & 2 & 2 & 3 & 1 & 2 & 2 & 3 & 2 \\
2 & 2 & 3 & 2 & 2 & 1 & 3 & 2 & 2 \\
2 & 2 & 2 & 2 & 2 & 3 & 1 & 2 & 2 \\
2 & 2 & 2 & 2 & 3 & 2 & 2 & 1 & 3 \\
2 & 2 & 2 & 2 & 2 & 2 & 2 & 3 & 1
\end{pmatrix}\]
Numerator of $p_{(W,S)}(t)$:
\begin{dmath*} n_{(W,S)}(t) = 
-1-8t-35t^{2}-112t^{3}-294t^{4}-672t^{5}-1386t^{6}-2640t^{7}-4718t^{8}-8000t^{9}-12978t^{10}-20272t^{11}-30645t^{12}-45016t^{13}-64470t^{14}-90264t^{15}-123829t^{16}-166768t^{17}-220849t^{18}-287992t^{19}-370250t^{20}-469784t^{21}-588833t^{22}-729680t^{23}-894613t^{24}-1085880t^{25}-1305640t^{26}-1555912t^{27}-1838523t^{28}-2155056t^{29}-2506798t^{30}-2894688t^{31}-3319268t^{32}-3780640t^{33}-4278429t^{34}-4811752t^{35}-5379194t^{36}-5978792t^{37}-6608029t^{38}-7263840t^{39}-7942628t^{40}-8640288t^{41}-9352240t^{42}-10073472t^{43}-10798593t^{44}-11521896t^{45}-12237428t^{46}-12939064t^{47}-13620586t^{48}-14275768t^{49}-14898464t^{50}-15482696t^{51}-16022740t^{52}-16513208t^{53}-16949127t^{54}-17326016t^{55}-17639957t^{56}-17887656t^{57}-18066494t^{58}-18174568t^{59}-18210722t^{60}-18174568t^{61}-18066494t^{62}-17887656t^{63}-17639957t^{64}-17326016t^{65}-16949127t^{66}-16513208t^{67}-16022740t^{68}-15482696t^{69}-14898464t^{70}-14275768t^{71}-13620586t^{72}-12939064t^{73}-12237428t^{74}-11521896t^{75}-10798593t^{76}-10073472t^{77}-9352240t^{78}-8640288t^{79}-7942628t^{80}-7263840t^{81}-6608029t^{82}-5978792t^{83}-5379194t^{84}-4811752t^{85}-4278429t^{86}-3780640t^{87}-3319268t^{88}-2894688t^{89}-2506798t^{90}-2155056t^{91}-1838523t^{92}-1555912t^{93}-1305640t^{94}-1085880t^{95}-894613t^{96}-729680t^{97}-588833t^{98}-469784t^{99}-370250t^{100}-287992t^{101}-220849t^{102}-166768t^{103}-123829t^{104}-90264t^{105}-64470t^{106}-45016t^{107}-30645t^{108}-20272t^{109}-12978t^{110}-8000t^{111}-4718t^{112}-2640t^{113}-1386t^{114}-672t^{115}-294t^{116}-112t^{117}-35t^{118}-8t^{119}-t^{120}
\end{dmath*}
Denominator of $p_{(W,S)}(t)$:
\begin{dmath*} d_{(W,S)}(t) = 
-1+t+t^{5}-t^{6}+t^{7}-t^{10}+2t^{11}-2t^{12}+t^{13}-t^{14}+t^{15}-3t^{16}+3t^{17}-3t^{18}+t^{19}-2t^{20}+3t^{21}-5t^{22}+4t^{23}-3t^{24}+2t^{25}-4t^{26}+6t^{27}-6t^{28}+5t^{29}-3t^{30}+4t^{31}-6t^{32}+9t^{33}-6t^{34}+5t^{35}-3t^{36}+7t^{37}-8t^{38}+10t^{39}-5t^{40}+4t^{41}-4t^{42}+9t^{43}-9t^{44}+8t^{45}-3t^{46}+3t^{47}-6t^{48}+9t^{49}-8t^{50}+4t^{51}-2t^{52}+2t^{53}-7t^{54}+7t^{55}-5t^{56}-t^{57}-t^{58}+t^{59}-6t^{60}+2t^{61}-t^{62}-4t^{63}+t^{65}-3t^{66}-2t^{67}+3t^{68}-5t^{69}+t^{70}+t^{72}-5t^{73}+6t^{74}-4t^{75}+2t^{76}-t^{77}+5t^{78}-6t^{79}+6t^{80}-2t^{81}+3t^{82}-2t^{83}+6t^{84}-5t^{85}+5t^{86}-t^{87}+3t^{88}-3t^{89}+5t^{90}-3t^{91}+3t^{92}-t^{93}+2t^{94}-3t^{95}+3t^{96}-2t^{97}+t^{98}-t^{99}+t^{100}-2t^{101}+t^{102}-t^{103}-t^{105}-t^{107}-t^{111}+t^{119}
\end{dmath*}
Initial values of $a_k$:
\begin{dmath*}\breakingcomma
 (a_k)_{k=0}^{120}=(1
, 9
, 44
, 156
, 450
, 1123
, 2517
, 5193
, 10032
, 18370
, 32176
, 54285
, 88701
, 140988
, 218771
, 332374
, 495629
, 726898
, 1050360
, 1497626
, 2109759
, 2939794
, 4055876
, 5545161
, 7518656
, 10117212
, 13518932
, 17948316
, 23687537
, 31090329
, 40599072
, 52765787
, 68277910
, 87989906
, 112962016
, 144507711
, 184251767
, 234201290
, 296832525
, 375196897
, 473050476
, 595011959
, 746755357
, 935244905
, 1169021330
, 1458550571
, 1816648421
, 2258997441
, 2804775992
, 3477423470
, 4305570972
, 5324172854
, 6575882199
, 8112722371
, 9998117935
, 12309361684
, 15140610829
, 18606525178
, 22846684084
, 28030947969
, 34365965403
, 42103069336
, 51547857712
, 63071816245
, 77126416907
, 94260217469
, 115139598622
, 140573909887
, 171545958647
, 209248974207
, 255131418070
, 310951301418
, 378842021757
, 461392155724
, 561742159771
, 683701553741
, 831890917132
, 1011913941820
, 1230565891741
, 1496086160121
, 1818464237475
, 2209810368294
, 2684804553198
, 3261240433597
, 3960684083221
, 4809271953106
, 5838677328577
, 7087280845841
, 8601588108887
, 10437946519118
, 12664624412986
, 15364328899320
, 18637254885238
, 22604777267307
, 27413921856917
, 33242779169937
, 40307059786812
, 48868031846809
, 59242131912162
, 71812601780226
, 87043578080842
, 105497151391296
, 127854020424211
, 154938498579099
, 187748789619442
, 227493642282737
, 275636727318469
, 333950363333666
, 404580560263369
, 490125763809582
, 593732185981030
, 719209214285061
, 871169127405874
, 1055196235262979
, 1278051638757815
, 1547921108718555
, 1874715162276619
, 2270432325932566
, 2749598887833535
, 3329801241948252
, 4032330316365198
).
 \end{dmath*}
Exponential growth rate:
\[\omega(W,S) =
1.21043510584084537675231410608
\dots\]
Cocompact? No\\
$(W,S)\in \CM$? Yes\\
Magma command: \texttt{HyperbolicCoxeterGraph(66)} or \texttt{HyperbolicCoxeterMatrix(66)} \\
\newpage
% produced with RUNfortables.magma
\subsubsection{\EHNC{53}}
\setcounter{mycnt}{53}
\[\Gamma(W,S) = \quad\csname CoxGrHNC\Roman{mycnt}\endcsname\]

\vskip12pt 
\[ M=\begin{pmatrix} 
1 & 3 & 2 & 2 & 2 & 2 & 2 & 2 & 2 \\
3 & 1 & 3 & 3 & 2 & 2 & 2 & 2 & 2 \\
2 & 3 & 1 & 2 & 2 & 3 & 2 & 2 & 2 \\
2 & 3 & 2 & 1 & 3 & 2 & 2 & 2 & 2 \\
2 & 2 & 2 & 3 & 1 & 2 & 2 & 2 & 2 \\
2 & 2 & 3 & 2 & 2 & 1 & 3 & 2 & 2 \\
2 & 2 & 2 & 2 & 2 & 3 & 1 & 3 & 2 \\
2 & 2 & 2 & 2 & 2 & 2 & 3 & 1 & 4 \\
2 & 2 & 2 & 2 & 2 & 2 & 2 & 4 & 1
\end{pmatrix}\]
Numerator of $p_{(W,S)}(t)$:
\begin{dmath*} n_{(W,S)}(t) = 
-1-8t-35t^{2}-112t^{3}-294t^{4}-672t^{5}-1386t^{6}-2640t^{7}-4719t^{8}-8008t^{9}-13013t^{10}-20384t^{11}-30939t^{12}-45688t^{13}-65856t^{14}-92904t^{15}-128547t^{16}-174768t^{17}-233827t^{18}-308264t^{19}-400895t^{20}-514800t^{21}-653303t^{22}-819944t^{23}-1018442t^{24}-1252648t^{25}-1526489t^{26}-1843904t^{27}-2208773t^{28}-2624840t^{29}-3095631t^{30}-3624368t^{31}-4213881t^{32}-4866520t^{33}-5584069t^{34}-6367664t^{35}-7217717t^{36}-8133848t^{37}-9114827t^{38}-10158528t^{39}-11261896t^{40}-12420928t^{41}-13630669t^{42}-14885224t^{43}-16177787t^{44}-17500688t^{45}-18845457t^{46}-20202904t^{47}-21563214t^{48}-22916056t^{49}-24250704t^{50}-25556168t^{51}-26821333t^{52}-28035104t^{53}-29186555t^{54}-30265080t^{55}-31260543t^{56}-32163424t^{57}-32964958t^{58}-33657264t^{59}-34233462t^{60}-34687776t^{61}-35015621t^{62}-35213672t^{63}-35279914t^{64}-35213672t^{65}-35015621t^{66}-34687776t^{67}-34233462t^{68}-33657264t^{69}-32964958t^{70}-32163424t^{71}-31260543t^{72}-30265080t^{73}-29186555t^{74}-28035104t^{75}-26821333t^{76}-25556168t^{77}-24250704t^{78}-22916056t^{79}-21563214t^{80}-20202904t^{81}-18845457t^{82}-17500688t^{83}-16177787t^{84}-14885224t^{85}-13630669t^{86}-12420928t^{87}-11261896t^{88}-10158528t^{89}-9114827t^{90}-8133848t^{91}-7217717t^{92}-6367664t^{93}-5584069t^{94}-4866520t^{95}-4213881t^{96}-3624368t^{97}-3095631t^{98}-2624840t^{99}-2208773t^{100}-1843904t^{101}-1526489t^{102}-1252648t^{103}-1018442t^{104}-819944t^{105}-653303t^{106}-514800t^{107}-400895t^{108}-308264t^{109}-233827t^{110}-174768t^{111}-128547t^{112}-92904t^{113}-65856t^{114}-45688t^{115}-30939t^{116}-20384t^{117}-13013t^{118}-8008t^{119}-4719t^{120}-2640t^{121}-1386t^{122}-672t^{123}-294t^{124}-112t^{125}-35t^{126}-8t^{127}-t^{128}
\end{dmath*}
Denominator of $p_{(W,S)}(t)$:
\begin{dmath*} d_{(W,S)}(t) = 
-1+t+t^{3}-t^{4}+t^{5}-t^{6}+2t^{7}-3t^{8}+3t^{9}-3t^{10}+4t^{11}-5t^{12}+5t^{13}-6t^{14}+7t^{15}-10t^{16}+9t^{17}-11t^{18}+12t^{19}-15t^{20}+15t^{21}-18t^{22}+19t^{23}-23t^{24}+24t^{25}-26t^{26}+29t^{27}-32t^{28}+35t^{29}-35t^{30}+40t^{31}-42t^{32}+47t^{33}-45t^{34}+53t^{35}-52t^{36}+59t^{37}-55t^{38}+63t^{39}-61t^{40}+68t^{41}-63t^{42}+71t^{43}-68t^{44}+73t^{45}-68t^{46}+73t^{47}-72t^{48}+72t^{49}-69t^{50}+69t^{51}-69t^{52}+65t^{53}-65t^{54}+59t^{55}-62t^{56}+52t^{57}-55t^{58}+44t^{59}-49t^{60}+35t^{61}-40t^{62}+26t^{63}-32t^{64}+17t^{65}-22t^{66}+7t^{67}-13t^{68}-t^{69}-3t^{70}-9t^{71}+6t^{72}-16t^{73}+14t^{74}-22t^{75}+22t^{76}-26t^{77}+28t^{78}-30t^{79}+32t^{80}-32t^{81}+36t^{82}-33t^{83}+38t^{84}-33t^{85}+38t^{86}-32t^{87}+38t^{88}-31t^{89}+35t^{90}-28t^{91}+33t^{92}-26t^{93}+29t^{94}-23t^{95}+26t^{96}-21t^{97}+21t^{98}-17t^{99}+18t^{100}-15t^{101}+13t^{102}-12t^{103}+11t^{104}-10t^{105}+7t^{106}-8t^{107}+6t^{108}-6t^{109}+3t^{110}-4t^{111}+2t^{112}-3t^{113}+t^{114}-2t^{115}+t^{116}-t^{117}+t^{127}
\end{dmath*}
Initial values of $a_k$:
\begin{dmath*}\breakingcomma
 (a_k)_{k=0}^{128}=(1
, 9
, 44
, 157
, 459
, 1167
, 2674
, 5653
, 11208
, 21089
, 37995
, 65997
, 111123
, 182161
, 291755
, 457893
, 705918
, 1071235
, 1602942
, 2368685
, 3461132
, 5006586
, 7176420
, 10202230
, 14395879
, 20175969
, 28102753
, 38924119
, 53636087
, 73562316
, 100458493
, 136649269
, 185207745
, 250190556
, 336945566
, 452514352
, 606158380
, 810046527
, 1080152993
, 1437429467
, 1909334694
, 2531829671
, 3351979324
, 4431343944
, 5850398831
, 7714292326
, 10160345670
, 13367819366
, 17570628311
, 23073892797
, 30275478701
, 39694026167
, 52005415693
, 68090204792
, 89095327570
, 116514335985
, 152291743179
, 198958694404
, 259809354347
, 339130210052
, 442498139585
, 577167839357
, 752575363948
, 980992535215
, 1278377372904
, 1665479201873
, 2169274629918
, 2824833370994
, 3677742477525
, 4787255976494
, 6230386818680
, 8107222879329
, 10547832946993
, 13721237991138
, 17847065022438
, 23210685310567
, 30182878275271
, 39245373474627
, 51024027151611
, 66331914520801
, 86225300424048
, 112076335963556
, 145667478016854
, 189314121085900
, 246023869238781
, 319703393084125
, 415427085642637
, 539785976134320
, 701340873604478
, 911210871536687
, 1183837641758625
, 1537978019444895
, 1997993059895361
, 2595522108601397
, 3371656866952006
, 4379764772271412
, 5689155600307147
, 7389843102345769
, 9598728683511850
, 12467631776291851
, 16193718368252024
, 21033043811450286
, 27318139878878764
, 35480853721837394
, 46082006986499405
, 59849911625682080
, 77730387042784396
, 100951712878179705
, 131108977208806458
, 170273611581569922
, 221135633563548015
, 287188363087169296
, 372968294964346399
, 484366596739304009
, 629033618520299751
, 816904187173487870
, 1060879749675404442
, 1377714198797812744
, 1789164197912778502
, 2323482980560599246
, 3017360181158529376
, 3918440874703318769
, 5088596767750047628
, 6608174120222506986
, 8581510032698709984
, 11144095810139660266
, 14471879188487321288
, 18793344047908898668
, 24405196916250410177
).
 \end{dmath*}
Exponential growth rate:
\[\omega(W,S) =
1.29857549835704593429648428209
\dots\]
Cocompact? No\\
$(W,S)\in \CM$? Yes\\
Magma command: \texttt{HyperbolicCoxeterGraph(67)} or \texttt{HyperbolicCoxeterMatrix(67)} \\
\newpage
% produced with printtables.magma
\subsubsection{\EHNC{54}}
\setcounter{mycnt}{54}
\[\Gamma(W,S) = \quad\csname CoxGrHNC\Roman{mycnt}\endcsname\]

\vskip12pt 
\[ M=\begin{pmatrix} 
1 & 3 & 2 & 2 & 2 & 2 & 2 & 2 & 2 \\
3 & 1 & 3 & 3 & 2 & 2 & 2 & 2 & 2 \\
2 & 3 & 1 & 2 & 2 & 3 & 2 & 2 & 2 \\
2 & 3 & 2 & 1 & 3 & 2 & 2 & 2 & 2 \\
2 & 2 & 2 & 3 & 1 & 2 & 2 & 2 & 2 \\
2 & 2 & 3 & 2 & 2 & 1 & 3 & 2 & 2 \\
2 & 2 & 2 & 2 & 2 & 3 & 1 & 3 & 3 \\
2 & 2 & 2 & 2 & 2 & 2 & 3 & 1 & 2 \\
2 & 2 & 2 & 2 & 2 & 2 & 3 & 2 & 1
\end{pmatrix}\]
Numerator of $p_{(W,S)}(t)$:
\begin{dmath*} n_{(W,S)}(t) = 
-1-8t-35t^{2}-112t^{3}-294t^{4}-672t^{5}-1386t^{6}-2640t^{7}-4718t^{8}-8000t^{9}-12978t^{10}-20272t^{11}-30645t^{12}-45016t^{13}-64470t^{14}-90264t^{15}-123829t^{16}-166768t^{17}-220849t^{18}-287992t^{19}-370250t^{20}-469784t^{21}-588833t^{22}-729680t^{23}-894613t^{24}-1085880t^{25}-1305640t^{26}-1555912t^{27}-1838523t^{28}-2155056t^{29}-2506798t^{30}-2894688t^{31}-3319268t^{32}-3780640t^{33}-4278429t^{34}-4811752t^{35}-5379194t^{36}-5978792t^{37}-6608029t^{38}-7263840t^{39}-7942628t^{40}-8640288t^{41}-9352240t^{42}-10073472t^{43}-10798593t^{44}-11521896t^{45}-12237428t^{46}-12939064t^{47}-13620586t^{48}-14275768t^{49}-14898464t^{50}-15482696t^{51}-16022740t^{52}-16513208t^{53}-16949127t^{54}-17326016t^{55}-17639957t^{56}-17887656t^{57}-18066494t^{58}-18174568t^{59}-18210722t^{60}-18174568t^{61}-18066494t^{62}-17887656t^{63}-17639957t^{64}-17326016t^{65}-16949127t^{66}-16513208t^{67}-16022740t^{68}-15482696t^{69}-14898464t^{70}-14275768t^{71}-13620586t^{72}-12939064t^{73}-12237428t^{74}-11521896t^{75}-10798593t^{76}-10073472t^{77}-9352240t^{78}-8640288t^{79}-7942628t^{80}-7263840t^{81}-6608029t^{82}-5978792t^{83}-5379194t^{84}-4811752t^{85}-4278429t^{86}-3780640t^{87}-3319268t^{88}-2894688t^{89}-2506798t^{90}-2155056t^{91}-1838523t^{92}-1555912t^{93}-1305640t^{94}-1085880t^{95}-894613t^{96}-729680t^{97}-588833t^{98}-469784t^{99}-370250t^{100}-287992t^{101}-220849t^{102}-166768t^{103}-123829t^{104}-90264t^{105}-64470t^{106}-45016t^{107}-30645t^{108}-20272t^{109}-12978t^{110}-8000t^{111}-4718t^{112}-2640t^{113}-1386t^{114}-672t^{115}-294t^{116}-112t^{117}-35t^{118}-8t^{119}-t^{120}
\end{dmath*}
Denominator of $p_{(W,S)}(t)$:
\begin{dmath*} d_{(W,S)}(t) = 
-1+t+t^{3}-t^{4}+t^{5}+t^{7}-2t^{8}+t^{9}-t^{10}+t^{11}-2t^{12}+t^{13}-2t^{14}+t^{15}-3t^{16}+t^{17}-t^{18}+t^{19}-3t^{20}+t^{21}+t^{23}-2t^{24}+t^{25}+t^{26}+t^{27}+t^{28}+t^{29}+3t^{30}+3t^{32}+6t^{34}+4t^{36}-t^{37}+7t^{38}-t^{39}+4t^{40}-2t^{41}+6t^{42}-2t^{43}+4t^{44}-2t^{45}+3t^{46}-3t^{47}-2t^{49}-2t^{51}-4t^{52}-2t^{53}-4t^{54}-8t^{56}-9t^{58}+t^{59}-11t^{60}+3t^{61}-12t^{62}+3t^{63}-12t^{64}+5t^{65}-12t^{66}+6t^{67}-12t^{68}+6t^{69}-10t^{70}+8t^{71}-10t^{72}+7t^{73}-8t^{74}+7t^{75}-5t^{76}+8t^{77}-5t^{78}+6t^{79}-2t^{80}+6t^{81}-t^{82}+5t^{83}+t^{84}+3t^{85}+2t^{86}+3t^{87}+3t^{88}+t^{89}+2t^{90}+4t^{92}+2t^{94}-2t^{95}+3t^{96}-2t^{97}+2t^{98}-2t^{99}+2t^{100}-3t^{101}+t^{102}-2t^{103}+t^{104}-2t^{105}-2t^{107}+t^{108}-t^{109}-t^{111}+t^{119}
\end{dmath*}
Initial values of $a_k$:
\begin{dmath*}\breakingcomma
 (a_k)_{k=0}^{120}=(1
, 9
, 44
, 157
, 459
, 1167
, 2675
, 5662
, 11252
, 21246
, 38455
, 67173
, 113843
, 187990
, 303520
, 480517
, 747719
, 1145917
, 1732604
, 2588317
, 3825267
, 5599052
, 8124523
, 11697239
, 16722428
, 23754022
, 33547203
, 47129045
, 65893379
, 91728068
, 127185605
, 175711588
, 241950487
, 332154567
, 454730421
, 620969020
, 846020396
, 1150194303
, 1560695140
, 2113935213
, 2858617990
, 3859846301
, 5204594546
, 7008995734
, 9428042769
, 12668500839
, 17006090055
, 22808346059
, 30565029415
, 40928569712
, 54767847512
, 73239702815
, 97884000448
, 130749997773
, 174564303550
, 232954094510
, 310743741989
, 414348958224
, 552300482616
, 735939833010
, 980343596314
, 1305551255020
, 1738196140293
, 2313671758498
, 3079009097595
, 4096698090037
, 5449762846614
, 7248501762278
, 9639438337006
, 12817207448420
, 17040339321990
, 22652218771092
, 30108915916644
, 40016140386775
, 53178308867880
, 70663695474810
, 93890934909392
, 124743874895194
, 165724066456187
, 220153223494245
, 292442022731648
, 388446977884164
, 515944241304424
, 685258637606732
, 910098780707438
, 1208665782182757
, 1605125170803392
, 2131560997310281
, 2830570066882411
, 3758705973057042
, 4991051281061408
, 6627287374519898
, 8799752503483275
, 11684139232271074
, 15513695763782814
, 20598078745008616
, 27348381012241251
, 36310356681461428
, 48208528342188718
, 64004740386520904
, 84975889739589802
, 112817114758860012
, 149778780027452800
, 198848325369412222
, 263991672276043582
, 350473692914286353
, 465283634761974874
, 617699873770619830
, 820039625897382703
, 1088654190341044596
, 1445250135168415186
, 1918643169912321727
, 2547086407703971283
, 3381361125858834484
, 4488879737771917988
, 5959132467517148269
, 7910917778617442198
, 10501940720669299565
, 13941554664738054174
, 18507675856438148565
, 24569237341673508883
).
 \end{dmath*}
Exponential growth rate:
\[\omega(W,S) =
1.32748480370989959113896848419
\dots\]
Cocompact? No\\
$(W,S)\in \CM$? Yes\\
Magma command: \texttt{HyperbolicCoxeterGraph(68)} or \texttt{HyperbolicCoxeterMatrix(68)} \\
\newpage
\subsubsection{\EHNC{55}}
\setcounter{mycnt}{55}
\[\Gamma(W,S) = \quad\csname CoxGrHNC\Roman{mycnt}\endcsname\]

\vskip12pt 
\[ M=\begin{pmatrix} 
1 & 3 & 2 & 2 & 2 & 2 & 2 & 2 & 2 \\
3 & 1 & 3 & 3 & 2 & 2 & 2 & 2 & 2 \\
2 & 3 & 1 & 2 & 2 & 3 & 2 & 2 & 2 \\
2 & 3 & 2 & 1 & 3 & 2 & 2 & 2 & 2 \\
2 & 2 & 2 & 3 & 1 & 2 & 2 & 3 & 2 \\
2 & 2 & 3 & 2 & 2 & 1 & 3 & 2 & 2 \\
2 & 2 & 2 & 2 & 2 & 3 & 1 & 2 & 3 \\
2 & 2 & 2 & 2 & 3 & 2 & 2 & 1 & 3 \\
2 & 2 & 2 & 2 & 2 & 2 & 3 & 3 & 1
\end{pmatrix}\]
Numerator of $p_{(W,S)}(t)$:
\begin{dmath*} n_{(W,S)}(t) = 
-1-8t-35t^{2}-112t^{3}-294t^{4}-672t^{5}-1386t^{6}-2640t^{7}-4718t^{8}-8000t^{9}-12978t^{10}-20272t^{11}-30645t^{12}-45016t^{13}-64470t^{14}-90264t^{15}-123829t^{16}-166768t^{17}-220849t^{18}-287992t^{19}-370250t^{20}-469784t^{21}-588833t^{22}-729680t^{23}-894613t^{24}-1085880t^{25}-1305640t^{26}-1555912t^{27}-1838523t^{28}-2155056t^{29}-2506798t^{30}-2894688t^{31}-3319268t^{32}-3780640t^{33}-4278429t^{34}-4811752t^{35}-5379194t^{36}-5978792t^{37}-6608029t^{38}-7263840t^{39}-7942628t^{40}-8640288t^{41}-9352240t^{42}-10073472t^{43}-10798593t^{44}-11521896t^{45}-12237428t^{46}-12939064t^{47}-13620586t^{48}-14275768t^{49}-14898464t^{50}-15482696t^{51}-16022740t^{52}-16513208t^{53}-16949127t^{54}-17326016t^{55}-17639957t^{56}-17887656t^{57}-18066494t^{58}-18174568t^{59}-18210722t^{60}-18174568t^{61}-18066494t^{62}-17887656t^{63}-17639957t^{64}-17326016t^{65}-16949127t^{66}-16513208t^{67}-16022740t^{68}-15482696t^{69}-14898464t^{70}-14275768t^{71}-13620586t^{72}-12939064t^{73}-12237428t^{74}-11521896t^{75}-10798593t^{76}-10073472t^{77}-9352240t^{78}-8640288t^{79}-7942628t^{80}-7263840t^{81}-6608029t^{82}-5978792t^{83}-5379194t^{84}-4811752t^{85}-4278429t^{86}-3780640t^{87}-3319268t^{88}-2894688t^{89}-2506798t^{90}-2155056t^{91}-1838523t^{92}-1555912t^{93}-1305640t^{94}-1085880t^{95}-894613t^{96}-729680t^{97}-588833t^{98}-469784t^{99}-370250t^{100}-287992t^{101}-220849t^{102}-166768t^{103}-123829t^{104}-90264t^{105}-64470t^{106}-45016t^{107}-30645t^{108}-20272t^{109}-12978t^{110}-8000t^{111}-4718t^{112}-2640t^{113}-1386t^{114}-672t^{115}-294t^{116}-112t^{117}-35t^{118}-8t^{119}-t^{120}
\end{dmath*}
Denominator of $p_{(W,S)}(t)$:
\begin{dmath*} d_{(W,S)}(t) = 
-1+t+t^{2}-t^{8}-t^{9}-t^{10}-t^{12}-t^{13}-2t^{14}-3t^{16}+t^{17}-t^{18}+t^{19}-t^{20}+4t^{21}-t^{22}+5t^{23}+t^{24}+6t^{25}+t^{26}+10t^{27}+3t^{28}+11t^{29}+4t^{30}+12t^{31}+t^{32}+16t^{33}+4t^{34}+14t^{35}+2t^{36}+15t^{37}-t^{38}+15t^{39}-3t^{40}+10t^{41}-6t^{42}+11t^{43}-12t^{44}+7t^{45}-12t^{46}+t^{47}-22t^{48}+t^{49}-23t^{50}-5t^{51}-27t^{52}-10t^{53}-31t^{54}-9t^{55}-33t^{56}-16t^{57}-33t^{58}-17t^{59}-36t^{60}-16t^{61}-31t^{62}-19t^{63}-31t^{64}-15t^{65}-26t^{66}-14t^{67}-22t^{68}-14t^{69}-16t^{70}-8t^{71}-13t^{72}-7t^{73}-5t^{74}-5t^{75}-2t^{76}+t^{77}+5t^{78}+6t^{80}+3t^{81}+11t^{82}+6t^{83}+11t^{84}+4t^{85}+14t^{86}+6t^{87}+12t^{88}+5t^{89}+12t^{90}+3t^{91}+11t^{92}+3t^{93}+9t^{94}+t^{95}+6t^{96}-t^{97}+6t^{98}+2t^{100}-2t^{101}+2t^{102}-3t^{103}+t^{104}-2t^{105}-t^{106}-3t^{107}-2t^{109}-2t^{111}-t^{112}-t^{113}+t^{114}-t^{115}+t^{119}+t^{120}
\end{dmath*}
Initial values of $a_k$:
\begin{dmath*}\breakingcomma
 (a_k)_{k=0}^{120}=(1
, 9
, 45
, 166
, 505
, 1343
, 3234
, 7217
, 15168
, 30375
, 58466
, 108893
, 197287
, 349172
, 605791
, 1033204
, 1736441
, 2881447
, 4729014
, 7687121
, 12391500
, 19829426
, 31529624
, 49853209
, 78438880
, 122883443
, 191781113
, 298309483
, 462648018
, 715663835
, 1104525809
, 1701251829
, 2615716264
, 4015437886
, 6155673070
, 9425168189
, 14415702273
, 22027767206
, 33631133187
, 51308762771
, 78227283722
, 119199613476
, 181539298185
, 276357683935
, 420533276349
, 639701374864
, 972792240121
, 1478919456150
, 2247835016912
, 3415797183752
, 5189652367208
, 7883381661728
, 11973561828816
, 18183527326042
, 27611082863951
, 41922297747970
, 63645568511303
, 96617820257528
, 146661547854871
, 222612106515078
, 337876424991867
, 512798023521969
, 778245394718662
, 1181056515904121
, 1792299537093436
, 2719806230388852
, 4127188135794855
, 6262688312130295
, 9502955706302135
, 14419460110563698
, 21879257317490451
, 33197870119989258
, 50371238894068495
, 76427635462169755
, 115961583304206266
, 175943914469035511
, 266950775520151048
, 405028253819305037
, 614521352201274367
, 932366049358825991
, 1414601116301627307
, 2146247687519843765
, 3256298296783853198
, 4940457847917914017
, 7495645596348870793
, 11372340587760162444
, 17253998651673746573
, 26177541528857663848
, 39716159364852851959
, 60256653084779642426
, 91420209675024650357
, 138700789277133467477
, 210433665898111389093
, 319264866455640633349
, 484380556396274614712
, 734889391269675048187
, 1114954110168802585175
, 1691577020914623779147
, 2566411881460668090005
, 3893684495198719938940
, 5907381242453099997900
, 8962498613673968190698
, 13597625621418707776018
, 20629891737060758578778
, 31299018708989254024584
, 47485871129662715216698
, 72044033982430746326228
, 109302870891385539921751
, 165830747076201024426208
, 251592965363165470348195
, 381708547431558123419099
, 579115562686445812000681
, 878614874674833307159971
, 1333005175669404951239257
, 2022390852426875259872545
, 3068303563286387450586147
, 4655127014019199979582618
, 7062601897784866092320734
, 10715141282127376381903008
, 16256650232673061019323763
, 24664039594390624031805157
).
 \end{dmath*}
Exponential growth rate:
\[\omega(W,S) =
1.51716581622718753951807447475
\dots\]
Cocompact? No\\
$(W,S)\in \CM$? No\\
Magma command: \texttt{HyperbolicCoxeterGraph(69)} or \texttt{HyperbolicCoxeterMatrix(69)} \\
\newpage
\subsubsection{\EHNC{56}}
\setcounter{mycnt}{56}
\[\Gamma(W,S) = \quad\csname CoxGrHNC\Roman{mycnt}\endcsname\]

\vskip12pt 
\[ M=\begin{pmatrix} 
1 & 3 & 2 & 2 & 2 & 2 & 2 & 2 & 2 & 2 \\
3 & 1 & 3 & 3 & 2 & 2 & 2 & 2 & 2 & 2 \\
2 & 3 & 1 & 2 & 2 & 3 & 2 & 2 & 2 & 2 \\
2 & 3 & 2 & 1 & 3 & 2 & 2 & 2 & 2 & 2 \\
2 & 2 & 2 & 3 & 1 & 2 & 2 & 2 & 2 & 2 \\
2 & 2 & 3 & 2 & 2 & 1 & 3 & 2 & 2 & 2 \\
2 & 2 & 2 & 2 & 2 & 3 & 1 & 3 & 2 & 2 \\
2 & 2 & 2 & 2 & 2 & 2 & 3 & 1 & 3 & 2 \\
2 & 2 & 2 & 2 & 2 & 2 & 2 & 3 & 1 & 3 \\
2 & 2 & 2 & 2 & 2 & 2 & 2 & 2 & 3 & 1
\end{pmatrix}\]
Numerator of $p_{(W,S)}(t)$:
\begin{dmath*} n_{(W,S)}(t) = 
1+9t+43t^{2}+147t^{3}+406t^{4}+966t^{5}+2058t^{6}+4026t^{7}+7359t^{8}+12727t^{9}+21021t^{10}+33397t^{11}+51323t^{12}+76627t^{13}+111544t^{14}+158760t^{15}+221451t^{16}+303315t^{17}+408595t^{18}+542091t^{19}+709159t^{20}+915695t^{21}+1168103t^{22}+1473247t^{23}+1838386t^{24}+2271090t^{25}+2779137t^{26}+3370393t^{27}+4052677t^{28}+4833613t^{29}+5720471t^{30}+6719999t^{31}+7838249t^{32}+9080401t^{33}+10450589t^{34}+11951733t^{35}+13585381t^{36}+15351565t^{37}+17248675t^{38}+19273355t^{39}+21420424t^{40}+23682824t^{41}+26051597t^{42}+28515893t^{43}+31063011t^{44}+33678475t^{45}+36346145t^{46}+39048361t^{47}+41766118t^{48}+44479270t^{49}+47166760t^{50}+49806872t^{51}+52377501t^{52}+54856437t^{53}+57221659t^{54}+59451635t^{55}+61525623t^{56}+63423967t^{57}+65128382t^{58}+66622222t^{59}+67890726t^{60}+68921238t^{61}+69703397t^{62}+70229293t^{63}+70493586t^{64}+70493586t^{65}+70229293t^{66}+69703397t^{67}+68921238t^{68}+67890726t^{69}+66622222t^{70}+65128382t^{71}+63423967t^{72}+61525623t^{73}+59451635t^{74}+57221659t^{75}+54856437t^{76}+52377501t^{77}+49806872t^{78}+47166760t^{79}+44479270t^{80}+41766118t^{81}+39048361t^{82}+36346145t^{83}+33678475t^{84}+31063011t^{85}+28515893t^{86}+26051597t^{87}+23682824t^{88}+21420424t^{89}+19273355t^{90}+17248675t^{91}+15351565t^{92}+13585381t^{93}+11951733t^{94}+10450589t^{95}+9080401t^{96}+7838249t^{97}+6719999t^{98}+5720471t^{99}+4833613t^{100}+4052677t^{101}+3370393t^{102}+2779137t^{103}+2271090t^{104}+1838386t^{105}+1473247t^{106}+1168103t^{107}+915695t^{108}+709159t^{109}+542091t^{110}+408595t^{111}+303315t^{112}+221451t^{113}+158760t^{114}+111544t^{115}+76627t^{116}+51323t^{117}+33397t^{118}+21021t^{119}+12727t^{120}+7359t^{121}+4026t^{122}+2058t^{123}+966t^{124}+406t^{125}+147t^{126}+43t^{127}+9t^{128}+t^{129}
\end{dmath*}
Denominator of $p_{(W,S)}(t)$:
\begin{dmath*} d_{(W,S)}(t) = 
1-t-t^{2}+t^{3}-t^{7}+t^{8}+t^{9}-t^{10}-t^{11}+t^{13}+t^{18}-t^{19}+t^{21}-t^{22}+t^{24}+t^{26}-t^{27}-t^{28}+t^{30}-t^{32}+t^{36}-2t^{37}-t^{38}+t^{39}-t^{43}+t^{45}-t^{46}-t^{47}+2t^{48}-t^{50}+t^{51}+t^{54}-t^{55}+2t^{57}-t^{59}+2t^{60}-t^{62}+t^{63}+t^{66}-t^{67}+t^{69}-t^{71}-t^{74}+t^{75}-2t^{77}+t^{81}-2t^{82}-t^{83}+t^{84}-t^{86}+t^{93}-t^{94}-t^{95}+2t^{96}+t^{97}-t^{98}+t^{101}+t^{102}-t^{103}+t^{105}-t^{116}-t^{126}+t^{128}
\end{dmath*}
Initial values of $a_k$:
\begin{dmath*}\breakingcomma
 (a_k)_{k=0}^{129}=(1
, 10
, 54
, 210
, 660
, 1782
, 4290
, 9439
, 19315
, 37234
, 68278
, 120001
, 203344
, 333805
, 532919
, 830114
, 1265021
, 1890329
, 2775291
, 4010005
, 5710615
, 8025601
, 11143356
, 15301282
, 20796675
, 27999712
, 37368904
, 49469438
, 64994898
, 84792933
, 109895531
, 141554661
, 181284162
, 230908894
, 292622324
, 369053901
, 463347782
, 579254711
, 721239128
, 894603899
, 1105635418
, 1361772247
, 1671800936
, 2046083209
, 2496819326
, 3038353148
, 3687525252
, 4464081381
, 5391144590
, 6495760682
, 7809527936
, 9369323740
, 11218142588
, 13406062012
, 15991355433
, 19041773675
, 22636020043
, 26865447474
, 31836010393
, 37670508620
, 44511166062
, 52522593077
, 61895188422
, 72849044722
, 85638430564
, 100556932786
, 117943354482
, 138188477891
, 161742816915
, 189125501786
, 220934458690
, 257858070311
, 300688529679
, 350337129850
, 407851766338
, 474436968451
, 551476820432
, 640561184346
, 743515694870
, 862436062536
, 999727297688
, 1158148553744
, 1340864386789
, 1551503340752
, 1794224895363
, 2073795959943
, 2395678262351
, 2766128171937
, 3192310711372
, 3682429758441
, 4245876719495
, 4893400276052
, 5637300170476
, 6491648411961
, 7472541757274
, 8598389859929
, 9890244095851
, 11372172773596
, 13071689234721
, 15020240258507
, 17253763220362
, 19813321632446
, 22745830038376
, 26104880764079
, 29951686769914
, 34356156834627
, 39398121563175
, 45168731286271
, 51772049853271
, 59326871661399
, 67968793069831
, 77852573681215
, 89154827909069
, 102077092870613
, 116849325045960
, 133733885434364
, 153030081239518
, 175079341568896
, 200271115396429
, 229049592295061
, 261921360402763
, 299464131977844
, 342336684994929
, 391290189836894
, 447181113597707
, 510985921221938
, 583817823118798
, 666945853514324
, 761816603228029
, 870078975443968
).
 \end{dmath*}
Exponential growth rate:
\[\omega(W,S) =
%% addded by hand
\tau =
1.13807874335567263390849987408
\dots\]
Cocompact? No\\
$(W,S)\in \CM$? Yes\\
Magma command: \texttt{HyperbolicCoxeterGraph(70)} or \texttt{HyperbolicCoxeterMatrix(70)} \\
\newpage
\subsubsection{\EHNC{57}}
\setcounter{mycnt}{57}
\[\Gamma(W,S) = \quad\csname CoxGrHNC\Roman{mycnt}\endcsname\]

\vskip12pt 
\[ M=\begin{pmatrix} 
1 & 3 & 2 & 2 & 2 & 2 & 2 & 2 & 2 & 2 \\
3 & 1 & 3 & 3 & 2 & 2 & 2 & 2 & 2 & 2 \\
2 & 3 & 1 & 2 & 2 & 3 & 2 & 2 & 2 & 2 \\
2 & 3 & 2 & 1 & 3 & 2 & 2 & 2 & 2 & 2 \\
2 & 2 & 2 & 3 & 1 & 2 & 2 & 2 & 2 & 2 \\
2 & 2 & 3 & 2 & 2 & 1 & 3 & 2 & 2 & 2 \\
2 & 2 & 2 & 2 & 2 & 3 & 1 & 3 & 2 & 2 \\
2 & 2 & 2 & 2 & 2 & 2 & 3 & 1 & 3 & 2 \\
2 & 2 & 2 & 2 & 2 & 2 & 2 & 3 & 1 & 4 \\
2 & 2 & 2 & 2 & 2 & 2 & 2 & 2 & 4 & 1
\end{pmatrix}\]
Numerator of $p_{(W,S)}(t)$:
\begin{dmath*} n_{(W,S)}(t) = 
1+9t+43t^{2}+147t^{3}+406t^{4}+966t^{5}+2058t^{6}+4026t^{7}+7359t^{8}+12727t^{9}+21021t^{10}+33397t^{11}+51323t^{12}+76627t^{13}+111544t^{14}+158760t^{15}+221451t^{16}+303315t^{17}+408595t^{18}+542091t^{19}+709159t^{20}+915695t^{21}+1168103t^{22}+1473247t^{23}+1838386t^{24}+2271090t^{25}+2779137t^{26}+3370393t^{27}+4052677t^{28}+4833613t^{29}+5720471t^{30}+6719999t^{31}+7838249t^{32}+9080401t^{33}+10450589t^{34}+11951733t^{35}+13585381t^{36}+15351565t^{37}+17248675t^{38}+19273355t^{39}+21420424t^{40}+23682824t^{41}+26051597t^{42}+28515893t^{43}+31063011t^{44}+33678475t^{45}+36346145t^{46}+39048361t^{47}+41766118t^{48}+44479270t^{49}+47166760t^{50}+49806872t^{51}+52377501t^{52}+54856437t^{53}+57221659t^{54}+59451635t^{55}+61525623t^{56}+63423967t^{57}+65128382t^{58}+66622222t^{59}+67890726t^{60}+68921238t^{61}+69703397t^{62}+70229293t^{63}+70493586t^{64}+70493586t^{65}+70229293t^{66}+69703397t^{67}+68921238t^{68}+67890726t^{69}+66622222t^{70}+65128382t^{71}+63423967t^{72}+61525623t^{73}+59451635t^{74}+57221659t^{75}+54856437t^{76}+52377501t^{77}+49806872t^{78}+47166760t^{79}+44479270t^{80}+41766118t^{81}+39048361t^{82}+36346145t^{83}+33678475t^{84}+31063011t^{85}+28515893t^{86}+26051597t^{87}+23682824t^{88}+21420424t^{89}+19273355t^{90}+17248675t^{91}+15351565t^{92}+13585381t^{93}+11951733t^{94}+10450589t^{95}+9080401t^{96}+7838249t^{97}+6719999t^{98}+5720471t^{99}+4833613t^{100}+4052677t^{101}+3370393t^{102}+2779137t^{103}+2271090t^{104}+1838386t^{105}+1473247t^{106}+1168103t^{107}+915695t^{108}+709159t^{109}+542091t^{110}+408595t^{111}+303315t^{112}+221451t^{113}+158760t^{114}+111544t^{115}+76627t^{116}+51323t^{117}+33397t^{118}+21021t^{119}+12727t^{120}+7359t^{121}+4026t^{122}+2058t^{123}+966t^{124}+406t^{125}+147t^{126}+43t^{127}+9t^{128}+t^{129}
\end{dmath*}
Denominator of $p_{(W,S)}(t)$:
\begin{dmath*} d_{(W,S)}(t) = 
1-t-t^{2}+t^{4}-t^{7}+2t^{8}-t^{9}-t^{11}+2t^{12}-t^{13}+t^{14}-2t^{15}+4t^{16}-2t^{17}+2t^{18}-4t^{19}+4t^{20}-3t^{21}+3t^{22}-5t^{23}+7t^{24}-6t^{25}+4t^{26}-8t^{27}+8t^{28}-7t^{29}+6t^{30}-11t^{31}+11t^{32}-9t^{33}+7t^{34}-13t^{35}+14t^{36}-12t^{37}+9t^{38}-14t^{39}+17t^{40}-12t^{41}+11t^{42}-17t^{43}+18t^{44}-11t^{45}+12t^{46}-17t^{47}+22t^{48}-12t^{49}+12t^{50}-15t^{51}+18t^{52}-10t^{53}+13t^{54}-15t^{55}+18t^{56}-7t^{57}+9t^{58}-13t^{59}+15t^{60}-5t^{61}+4t^{62}-9t^{63}+11t^{64}-2t^{65}+t^{66}-8t^{67}+5t^{68}+3t^{69}-5t^{70}-5t^{71}+2t^{72}+5t^{73}-9t^{74}-4t^{76}+6t^{77}-11t^{78}+2t^{79}-4t^{80}+11t^{81}-15t^{82}+3t^{83}-6t^{84}+11t^{85}-13t^{86}+6t^{87}-7t^{88}+11t^{89}-10t^{90}+6t^{91}-7t^{92}+12t^{93}-11t^{94}+5t^{95}-2t^{96}+10t^{97}-8t^{98}+5t^{99}-4t^{100}+8t^{101}-5t^{102}+3t^{103}-2t^{104}+7t^{105}-4t^{106}+t^{107}-t^{108}+4t^{109}-4t^{110}+t^{111}+2t^{113}-2t^{114}-t^{116}+t^{117}-2t^{118}-t^{119}+t^{120}+t^{121}-t^{122}-t^{126}-t^{127}+t^{128}+t^{129}
\end{dmath*}
Initial values of $a_k$:
\begin{dmath*}\breakingcomma
 (a_k)_{k=0}^{129}=(1
, 10
, 54
, 211
, 670
, 1837
, 4511
, 10164
, 21372
, 42461
, 80456
, 146453
, 257576
, 439737
, 731492
, 1189385
, 1895303
, 2966538
, 4569479
, 6938155
, 10399242
, 15405661
, 22581568
, 32782420
, 47174953
, 67343412
, 95430337
, 134322770
, 187898084
, 261347981
, 361604860
, 497902108
, 682509426
, 931696722
, 1266996231
, 1716853460
, 2318784727
, 3122194308
, 4192049912
, 5613674455
, 7498988890
, 9994640337
, 13292578626
, 17643811249
, 23376282788
, 30918104583
, 40827722389
, 53833078175
, 70882428234
, 93210263700
, 122422793492
, 160608760934
, 210483060795
, 275572815611
, 360458404054
, 471085597317
, 615169694219
, 802718664905
, 1046710220764
, 1363967946511
, 1776294833455
, 2311939611487
, 3007493315690
, 3910341993799
, 5081838238162
, 6601401731762
, 8571820358954
, 11126102685163
, 14436334968866
, 18725128064688
, 24280410299034
, 31474542870055
, 40789019029006
, 52846375944171
, 68451422879386
, 88644502312624
, 114770292124740
, 148566678965260
, 192279552431940
, 248811073383127
, 321911169332696
, 416424849778674
, 538611600798466
, 696557851816559
, 900709618591768
, 1164560315974483
, 1505538919229874
, 1946156801585217
, 2515488550277946
, 3251083976659119
, 4201436824441640
, 5429172198402497
, 7015161877527333
, 9063837531692819
, 11710050418142653
, 15127927541843153
, 19542305168241197
, 25243489556668348
, 32606312911092122
, 42114734115797707
, 54393597287231435
, 70249630345345213
, 90724371429561216
, 117162492736480076
, 151300000457959614
, 195378092062406966
, 252290133527009368
, 325771389446942888
, 420643940402394984
, 543132838082343219
, 701274216337578165
, 905442101338822690
, 1169028441034311320
, 1509320912583409902
, 1948636023969676239
, 2515781751352945290
, 3247945542489153507
, 4193131382870789444
, 5413305590445056478
, 6988457432914535367
, 9021840589835510028
, 11646738835362918745
, 15035199163250762845
, 19409304453500869488
, 25055724131296354852
, 32344495989415000564
, 41753269499105833978
, 53898598674530470926
, 69576334313558198456
, 89813761456979011678
).
 \end{dmath*}
Exponential growth rate:
\[\omega(W,S) =
1.29074823652170787022549809467
\dots\]
Cocompact? No\\
$(W,S)\in \CM$? No\\
Magma command: \texttt{HyperbolicCoxeterGraph(71)} or \texttt{HyperbolicCoxeterMatrix(71)} \\
\newpage
\subsubsection{\EHNC{58}}
\setcounter{mycnt}{58}
\[\Gamma(W,S) = \quad\csname CoxGrHNC\Roman{mycnt}\endcsname\]

\vskip12pt 
\[ M=\begin{pmatrix} 
1 & 3 & 2 & 2 & 2 & 2 & 2 & 2 & 2 & 2 \\
3 & 1 & 3 & 3 & 2 & 2 & 2 & 2 & 2 & 2 \\
2 & 3 & 1 & 2 & 2 & 2 & 2 & 2 & 2 & 2 \\
2 & 3 & 2 & 1 & 3 & 2 & 2 & 2 & 2 & 2 \\
2 & 2 & 2 & 3 & 1 & 3 & 2 & 2 & 2 & 2 \\
2 & 2 & 2 & 2 & 3 & 1 & 3 & 2 & 2 & 2 \\
2 & 2 & 2 & 2 & 2 & 3 & 1 & 3 & 3 & 2 \\
2 & 2 & 2 & 2 & 2 & 2 & 3 & 1 & 2 & 2 \\
2 & 2 & 2 & 2 & 2 & 2 & 3 & 2 & 1 & 3 \\
2 & 2 & 2 & 2 & 2 & 2 & 2 & 2 & 3 & 1
\end{pmatrix}\]
Numerator of $p_{(W,S)}(t)$:
\begin{dmath*} n_{(W,S)}(t) = 
1+9t+43t^{2}+147t^{3}+406t^{4}+966t^{5}+2058t^{6}+4026t^{7}+7359t^{8}+12727t^{9}+21021t^{10}+33397t^{11}+51323t^{12}+76627t^{13}+111544t^{14}+158760t^{15}+221451t^{16}+303315t^{17}+408595t^{18}+542091t^{19}+709159t^{20}+915695t^{21}+1168103t^{22}+1473247t^{23}+1838386t^{24}+2271090t^{25}+2779137t^{26}+3370393t^{27}+4052677t^{28}+4833613t^{29}+5720471t^{30}+6719999t^{31}+7838249t^{32}+9080401t^{33}+10450589t^{34}+11951733t^{35}+13585381t^{36}+15351565t^{37}+17248675t^{38}+19273355t^{39}+21420424t^{40}+23682824t^{41}+26051597t^{42}+28515893t^{43}+31063011t^{44}+33678475t^{45}+36346145t^{46}+39048361t^{47}+41766118t^{48}+44479270t^{49}+47166760t^{50}+49806872t^{51}+52377501t^{52}+54856437t^{53}+57221659t^{54}+59451635t^{55}+61525623t^{56}+63423967t^{57}+65128382t^{58}+66622222t^{59}+67890726t^{60}+68921238t^{61}+69703397t^{62}+70229293t^{63}+70493586t^{64}+70493586t^{65}+70229293t^{66}+69703397t^{67}+68921238t^{68}+67890726t^{69}+66622222t^{70}+65128382t^{71}+63423967t^{72}+61525623t^{73}+59451635t^{74}+57221659t^{75}+54856437t^{76}+52377501t^{77}+49806872t^{78}+47166760t^{79}+44479270t^{80}+41766118t^{81}+39048361t^{82}+36346145t^{83}+33678475t^{84}+31063011t^{85}+28515893t^{86}+26051597t^{87}+23682824t^{88}+21420424t^{89}+19273355t^{90}+17248675t^{91}+15351565t^{92}+13585381t^{93}+11951733t^{94}+10450589t^{95}+9080401t^{96}+7838249t^{97}+6719999t^{98}+5720471t^{99}+4833613t^{100}+4052677t^{101}+3370393t^{102}+2779137t^{103}+2271090t^{104}+1838386t^{105}+1473247t^{106}+1168103t^{107}+915695t^{108}+709159t^{109}+542091t^{110}+408595t^{111}+303315t^{112}+221451t^{113}+158760t^{114}+111544t^{115}+76627t^{116}+51323t^{117}+33397t^{118}+21021t^{119}+12727t^{120}+7359t^{121}+4026t^{122}+2058t^{123}+966t^{124}+406t^{125}+147t^{126}+43t^{127}+9t^{128}+t^{129}
\end{dmath*}
Denominator of $p_{(W,S)}(t)$:
\begin{dmath*} d_{(W,S)}(t) = 
1-t-t^{2}+t^{4}-2t^{7}+3t^{8}-3t^{11}+3t^{12}-t^{13}+2t^{14}-4t^{15}+6t^{16}-3t^{17}+4t^{18}-8t^{19}+8t^{20}-5t^{21}+5t^{22}-11t^{23}+13t^{24}-9t^{25}+9t^{26}-17t^{27}+16t^{28}-14t^{29}+13t^{30}-22t^{31}+22t^{32}-17t^{33}+16t^{34}-28t^{35}+29t^{36}-24t^{37}+21t^{38}-31t^{39}+34t^{40}-24t^{41}+25t^{42}-37t^{43}+38t^{44}-24t^{45}+26t^{46}-38t^{47}+44t^{48}-24t^{49}+26t^{50}-34t^{51}+37t^{52}-21t^{53}+26t^{54}-33t^{55}+36t^{56}-13t^{57}+16t^{58}-27t^{59}+30t^{60}-9t^{61}+6t^{62}-17t^{63}+19t^{64}-t^{65}-13t^{67}+7t^{68}+9t^{69}-13t^{70}-5t^{71}+t^{72}+14t^{73}-22t^{74}+5t^{75}-10t^{76}+16t^{77}-24t^{78}+8t^{79}-13t^{80}+26t^{81}-31t^{82}+10t^{83}-13t^{84}+23t^{85}-29t^{86}+14t^{87}-14t^{88}+24t^{89}-22t^{90}+12t^{91}-14t^{92}+23t^{93}-21t^{94}+10t^{95}-5t^{96}+18t^{97}-16t^{98}+10t^{99}-6t^{100}+13t^{101}-10t^{102}+5t^{103}-3t^{104}+12t^{105}-7t^{106}+2t^{107}-t^{108}+5t^{109}-6t^{110}+2t^{111}+3t^{113}-3t^{114}-t^{116}+t^{117}-3t^{118}-2t^{119}+2t^{120}+2t^{121}-2t^{122}-t^{126}-2t^{127}+t^{128}+2t^{129}
\end{dmath*}
Initial values of $a_k$:
\begin{dmath*}\breakingcomma
 (a_k)_{k=0}^{129}=(1
, 10
, 54
, 211
, 670
, 1837
, 4511
, 10165
, 21382
, 42515
, 80667
, 147124
, 259423
, 444303
, 741879
, 1211502
, 1939931
, 3052663
, 4729562
, 7226320
, 10903806
, 16268060
, 24024622
, 35152125
, 51001727
, 73431144
, 104984944
, 149136603
, 210612915
, 295828090
, 413463788
, 575243121
, 796962233
, 1099863643
, 1512462700
, 2072974354
, 2832534761
, 3859474634
, 5244983516
, 7110612589
, 9618206530
, 12983043171
, 17491207666
, 23522554382
, 31581039603
, 42334773979
, 56668888359
, 75755286603
, 101144648267
, 134887740135
, 179695326493
, 239148902098
, 317978330141
, 422427540819
, 560736116358
, 743773357308
, 985872952056
, 1305931523169
, 1728854238613
, 2287456848404
, 3024967902312
, 3998320102490
, 5282479136264
, 6976136371713
, 9209194331304
, 12152608569608
, 16031326561296
, 21141296709594
, 27871826021838
, 36734966232055
, 48404135209992
, 63764872820880
, 83981539787007
, 110584962575702
, 145587596235031
, 191634837741854
, 252203828892897
, 331864642397744
, 436623413341938
, 574373109525902
, 755485686590926
, 993589949076104
, 1306593326831873
, 1718024015248241
, 2258793880215417
, 2969513984034111
, 3903535896404259
, 5130946200736404
, 6743812842160051
, 8863075510312859
, 11647595093860410
, 15306038558942342
, 20112487437192127
, 26426936279463903
, 34722212710744410
, 45619330385485035
, 59933916006281147
, 78737178661680897
, 103435975814670560
, 135877956420279891
, 178489634344674489
, 234457704310892605
, 307967141592444727
, 404513866676272551
, 531315323652770340
, 697849631792749168
, 916563569392237446
, 1203802254118631874
, 1581029935738905963
, 2076433050801478199
, 2727025226937244844
, 3581411397207223937
, 4703417389530468832
, 6176855965547140239
, 8111785120017481938
, 10652725848543309566
, 13989452863002237310
, 18371163798670491230
, 24125084650642000378
, 31680900323919205347
, 41602833998570946012
, 54631769953609256971
, 71740564175957243762
, 94206671460240699882
, 123707510254695723394
, 162445683703060369660
, 213313403848598492340
, 280108392137724381436
, 367817371591649486224
, 482988311081921174114
).
 \end{dmath*}
Exponential growth rate:
\[\omega(W,S) =
1.31304825690637874842910420169
\dots\]
Cocompact? No\\
$(W,S)\in \CM$? Yes\\
Magma command: \texttt{HyperbolicCoxeterGraph(72)} or \texttt{HyperbolicCoxeterMatrix(72)} \\

%%%%%%%%%%%%%%%%%%%%%%%%%%%%%%%%%%% 
%% References
\bibliographystyle{amsalpha}
\bibliography{all-references,terragni-references}

\bigskip
\addtocontents{toc}{\protect\end{multicols}}
\end{document}